\documentclass{informs-stsy}

\usepackage{endnotes}
\let\footnote=\endnote
\let\enotesize=\normalsize
\def\notesname{Endnotes}%
\def\makeenmark{\hbox to1.275em{\theenmark.\enskip\hss}}
\def\enoteformat{\rightskip0pt\leftskip0pt\parindent=1.275em
  \leavevmode\llap{\makeenmark}}

\usepackage{natbib}
 \bibpunct[, ]{(}{)}{,}{a}{}{,}%
 \def\bibfont{\small}%
\usepackage[bookmarksnumbered]{hyperref}
\usepackage[english]{babel}
\usepackage{amsmath,bbm}
\usepackage{graphicx}
\usepackage{amsfonts}
\usepackage{tikz}
\usepackage{mathtools}
\usepackage{float}
\usepackage{multirow}
\usepackage{multicol}
\usepackage{lscape}
\usepackage{amssymb}
\usepackage{chngpage}
\usepackage[T1]{fontenc}

\RequirePackage[normalem]{ulem} 
\RequirePackage{color}\definecolor{RED}{rgb}{1,0,0}\definecolor{BLUE}{rgb}{0,0,1} 

\newtheorem{thm}{Theorem}[section]
\newtheorem{cor}[thm]{Corollary}
\newtheorem{lem}[thm]{Lemma}
\newtheorem{prop}[thm]{Proposition}
\newtheorem{defn}[thm]{Definition}

\newcommand{\HidePart}[1]{ \iffalse #1 \fi  }

\newcommand{\probability}{\operatorname{\mathbb{P}}\probarg}
\DeclarePairedDelimiterX{\probarg}[1]{(}{)}{%
  \ifnum\currentgrouptype=16 \else\begingroup\fi
  \activatebar#1
  \ifnum\currentgrouptype=16 \else\endgroup\fi
}
\newcommand{\expect}{\operatorname{\mathbb{E}}\expectarg}
\DeclarePairedDelimiterX{\expectarg}[1]{[}{]}{%
  \ifnum\currentgrouptype=16 \else\begingroup\fi
  \activatebar#1
  \ifnum\currentgrouptype=16 \else\endgroup\fi
}

\newcommand{\activatebar}{%
  \begingroup\lccode`\~=`\|
  \lowercase{\endgroup\let~}\innermid 
  \mathcode`|=\string"8000
}

\newcommand{\InventoryExpCosts}[1]{\mathbb{E}\left[Nh^{(N)} \left(#1-Q_i+\left(\max_{j\leq N}Q_j-#1\right)^+\right)+b^{(N)}\left(\max_{j\leq N}Q_j-#1\right)^+\right]}
\newcommand{\OptimalWAIN}{I_N^A}
\newcommand{\OptimalWABetaN}{\beta_N^A}
\newcommand{\ApproxWAIN}{\hat{I}_N^A}
\newcommand{\ApproxWABetaN}{\hat{\beta}_N^A}
\newcommand{\ApproxWAC}[1]{\hat{C}^A_N(#1)}
\newcommand{\ApproxWACosts}{\hat{F}^A_N}

\newcommand{\ApproxC}[1]{\hat{C}_N(#1)}
\newcommand{\ApproxIN}{\hat{I}_N}
\newcommand{\ApproxBetaN}{\hat{\beta}_N}
\newcommand{\ApproxCosts}{\hat{F}_N}
\newcommand{\ApproxQuantile}{\hat{P}^{-1}_N}

\newcommand{\LimitN}{\overset{N\to\infty}{\longrightarrow}}
\newcommand{\LimitD}{\overset{d}{\longrightarrow}}
\newcommand{\LimitP}{\overset{\mathbb{P}}{\longrightarrow}}

\newcommand{\LimitL}{\overset{L_1}{\longrightarrow}}

\DeclareRobustCommand{\VAN}[3]{#2}

\begin{document}


\RUNAUTHOR{Meijer, Schol, van Jaarsveld, Vlasiou, and Zwart}

\RUNTITLE{Inventory and capacity in large-scale assembly systems}

\TITLE{Optimization of inventory and capacity in large-scale assembly systems using extreme-value theory}

\ARTICLEAUTHORS{%
\AUTHOR{Mirjam S. Meijer}
\AFF{Kühne Logistics University \EMAIL{mirjam.meijer@the-klu.org} \URL{}}
\AUTHOR{Dennis Schol}
\AFF{Eindhoven University of Technology \EMAIL{c.schol@tue.nl} \URL{}}
\AUTHOR{Willem van Jaarsveld}
\AFF{Eindhoven University of Technology \EMAIL{w.l.v.jaarsveld@tue.nl} \URL{}}
\AUTHOR{Maria Vlasiou}
\AFF{ University of Twente, Eindhoven University of Technology \EMAIL{m.vlasiou@utwente.nl} \URL{}}
\AUTHOR{Bert Zwart}
\AFF{Eindhoven University of Technology, CWI \EMAIL{bert.zwart@cwi.nl} \URL{}}
} 

\ABSTRACT{%
High-tech systems are typically produced in two stages: 1) Production of components using specialized equipment and staff; 2) System assembly/integration. Component production capacity is subject to fluctuations, causing a high risk of shortages of at least one component, which results in costly delays. Companies hedge this risk by strategic investments in excess production capacity and in buffer inventories of components. To optimize these, it is crucial to characterize the relation between component shortage risk and capacity and inventory investments. We suppose that component production capacity and produce demand are normally distributed over finite time intervals, and we accordingly model the production system as a {symmetric} fork-join queueing network with $N$ statistically identical queues with a common arrival process and independent service processes. {Assuming a symmetric cost structure,} we subsequently apply extreme value theory to gain analytic insights into this optimization problem. We derive several new results for this queueing network, notably that the scaled maximum of $N$ steady-state queue lengths converges in distribution to a Gaussian random variable. These results translate into asymptotically optimal methods to dimension the system. Tests on a range of problems reveal that these methods typically work well for systems of moderate size.
}%



\maketitle

%

\section{Introduction}
Delivery reliability is a key performance indicator for high-tech manufacturers such as ASML, Philips, and Airbus. High-tech systems such as wafer steppers, medical imaging equipment, and aircraft are produced by assembling thousands of components, each produced by highly skilled staff using specialized equipment. This production system facilitates modular design and testing, but it is also vulnerable: the shortage of a single component will result in delivery delays that cause customer grievances, a built-up of inventory of other components, and a severe reduction in turnover and cashflow. For example, in 2021 ASML was hit by material shortages in its supply chain, causing it to cut its revenue guidance \citep{barrons}. 
Also in other industries with higher demand volumes, e.g.\ car manufacturing, many components are required to assemble the final product and a single missing item can hinder production of the entire end-product. An example is the shutdown of complete manufacturing lines at several car manufacturers due to shortages of semiconductors \citep{ewing2021}.

Two complementary approaches may contribute to guaranteeing a reliable production system by reducing the risk of component shortages: excess component production capacity and inventory buffers. Production capacity and inventory buffers have a qualitatively different role in the mitigation of component shortages. Excess production capacity implies that the expected maximum number of components that can be produced per quarter exceeds the expected demand per quarter, e.g.\ as a rule, production capacity may be $110\%$ of expected demand. Inventory buffers are components that are produced in anticipation of demand; typically, such anticipative production continues until the inventory buffer reaches a \emph{target}, e.g.\ of 6 weeks of demand. Excess production capacity is always available, while inventory buffers are consumed when used to absorb production or demand fluctuations.

Joint optimization of excess component production capacity and component buffer is the ultimate goal because investments in excess component production capacity and component buffer inventories run into the hundreds of millions of euros \citep{asmlreport}. High-level investment plans for capacity and inventory may be devised for each product line (e.g.\ ASML's TWINSCAN XT range or Philips’ Azurion 7 C range), depending on the role of the product line in the company’s portfolio and other considerations. Despite the strategic importance of these investments, there is a lack of quantitative methods for determining appropriate investments in capacity and inventory to achieve the desired level of delivery reliability. Indeed, despite decades of research in inventory management, the joint optimization of production capacity and inventory remains a considerable challenge \citep{bradley2002}. While the topic has increasingly been studied \citep[see e.g.][]{reed2017}, the focus of analysis has been on problems with a single component. The much more common situation of assembling a system from many components has proved very challenging. 

In this paper, we make a step towards overcoming this challenge. We propose a stylized model capturing key features of high-tech manufacturing that is based on interactions with high-tech manufacturers in the Netherlands, and that yields new insights into the joint optimization of capacity and inventory for large-scale assembly systems.
We focus on a single product line. Typically, a majority of the expensive components used in high-tech products are common to all products in a product line, while being unique to that line, and we consider capacity and inventory optimization for those common components. Component shortages result in delays in the start of the assembly/integration process. Given the tight production planning that is common at high-tech manufacturers, such delays in turn result in costly delivery delays. Component production is \emph{capacitated} and subject to random fluctuations. E.g.\ the production capacity of components may be $\mu \pm \sigma$  items per quarter, and we assume a normal distribution for this per-period production capacity \citep[e.g.][]{bradley2002,wu2014}, which is the most natural assumption, as the stochastic term represents the error around the mean. We adopt a continuous-time model, and we likewise assume that production capacity in every finite interval is linear with normally distributed white noise, i.e.\ cumulative net production is a Brownian motion with drift $-\beta<0$ and variance $\sigma^2$ \citep[cf.][]{bradley2002,harrison2013}. We analyze the steady-state behavior of this system.

To analyze the overall production system, we consider a symmetric fork-join network of $N$ queues driven by a common arrival process and having independent, identical service processes. Due to this common arrival process, total inventory per component including backlogged items is equal for all components. However, as a result of variations in the service times, the number of backlogged items may
vary per component. We express the optimal component production capacity and inventory in this model in terms of the steady-state delay distribution of the slowest component, which has the form of a maximum of $N$ all-time suprema of Brownian motions, and we subsequently focus on analyzing this delay distribution. In particular, in large-scale systems with many components/queues, one can expect that the maximum delay (which is due to stochasticity of demand and service processes) grows without bound as a function of the size of the system. To analyze and quantify this phenomenon, we derive new analytic results for the delays in this fork-join network as $N\rightarrow \infty$. {To do so we make a major assumption, which is that the randomness and cost characteristics of each of the $N$ suppliers are identical, resulting in a symmetric system with identical net service capacities and basestock levels. 
 The symmetry we impose makes a mathematical treatment of our model within reach. While this is a shortcoming of our work, it already reveals useful insights, and we complement our analytic results with simulation experiments for asymmetric systems.}

\paragraph{\textbf{\emph{Extreme value analysis.}}}Original equipment manufacturers (OEMs) typically \emph{level} the demand to smooth the production process. Accordingly, in our base model, we assume that demand is completely leveled, which corresponds to a fork-join queue with a deterministic arrival stream. Extremes for this network as $N\rightarrow \infty$ are obtained using extreme value theory (EVT), and based on those results, in Section \ref{sec: basic model - independent} we derive easy-to-calculate expressions for capacity and inventory that are asymptotically optimal as the number of components grows large. We provide order bounds between the costs under optimal and approximate inventory and capacity.  
In particular, inspired by the literature on call centers \cite{borst2004dimensioning, gans2003telephone} and \cite{van2019economies}, we distinguish three regimes that depend on the growth rates of cost parameters and are determined by the probability $\gamma_N$ of not having enough inventory. Given that $\gamma_N\rightarrow \gamma$, we say that the regime is {\em balanced} if $\gamma \in (0,1)$. Furthermore, we are in the quality-driven regime if $\gamma=0$ and in the efficiency-driven regime if $\gamma=1$. For the base model, we establish asymptotic cost optimality in all three regimes. For the balanced, quality-driven, and efficiency-driven regimes, we have convergence rates of $1/(N\log N), \gamma_N/(N\log(N/\gamma_N))$ and $1/\log N$ respectively. 

\paragraph{\textbf{\emph{Demand fluctuations.}}} Other than the number of produced components being stochastic, despite efforts to level demand, typically some demand variation remains. Thus, a natural choice is that the demand has, apart from a linear term, a white noise term as well, which is normally distributed. Therefore in Section~\ref{sec: model - dependent}, we assume that the cumulative stochastic demand for systems is modeled by a Brownian motion with variance $\sigma_A^2$. \citep[cf.][for a single-component manufacturing system]{bradley2002}. This implies that the demand over any finite time period is a normal variable, which is a standard assumption in literature \citep[e.g.][]{klosterhalfen2014,atan2016}. In high-tech manufacturing, normally distributed demand is a suitable assumption especially when considering longer time periods, but it is also a reasonable approximation for shorter periods.  
As a consequence of these demand variations, component delays become \emph{dependent}, since they face the same stochastic demands from system assembly. The question is now how this affects the maximum delay as the number of queues/components $N\rightarrow\infty$. Most of the work in extreme value theory has been done for independent random variables; cf.\ \cite{de2007extreme,resnick2013extreme}, and suitable results from extreme value theory are absent for our setting, rendering the analysis of extremes in the dependent case challenging.

\paragraph{\textbf{\emph{New extreme-value limit.}}}
Our answer to this challenge is somewhat surprising: in Theorem \ref{thm: second order convergence process}, we prove that the scaled maximum queue length converges to a normally distributed random variable as $N\rightarrow\infty$. In particular, if $Q_i(\infty, \beta)$ is the invariant queue length at node $i$,
\begin{align}\label{eq: limit approximate problem dependent}
\frac{\max_{i\leq N} Q_i(\infty, \beta)-\frac{\sigma^2}{2\beta}\log N}{\sqrt{\log N}}\LimitD\frac{\sigma\sigma_A}{\sqrt{2}\beta}X,
\end{align}
with $X$ standard normal. 
An intuitive explanation of this result is the following. Using Lindley's recursion, we can write the maximum queue length as a maximum of $N$ suprema. By using subadditivity arguments, we can separate the independent and dependent part, the independent part converges using standard extreme-value results, while the dependent part satisfies a central limit theorem. To the best of our knowledge, we are the first who prove a result of this type. A consequence of this convergence result is that, with proper scaling of holding and backorder costs, the optimal inventory for stochastic demand converges to a scaled version of the quantile function of the normal distribution, while this quantile function also appears in the limit of the optimal capacity. 

\paragraph{\textbf{\emph{Numerical experiments.}}} In Section \ref{subsec: numerical results dependent}, numerical experiments show that we typically are most of the times 10\% off the optimum (e.g. when $N$ is in the range from 10 to 100); cf.~Tables \ref{tab: costs asymptotic} and \ref{tab: costs asymptotic b=3N}. Naturally, the difference goes to 0 as $N\to\infty$; cf.~Theorem \ref{thm: order bounds dep}. We give an improvement of this approximation by combining our results for deterministic demand and stochastic demand. Based on this approximation, we optimize the capacity and inventory decisions and we test the quality of these approximations through numerical experiments. It turns out that these approximations perform well already when considering a limited number of components, and are typically less than 2\% off the optimum.

\paragraph{\textbf{\emph{Limitations of simulation.}}} In Section \ref{subsec: numerical results dependent}, we explain the simulation procedure in the case of stochastic demand. We aim to approximate the maximum queue length of the all-time supremum of $N$ dependent Brownian motions. Because the dependence structure between two all-time suprema of Brownian motions is complicated, we cannot resort to an easy simulation procedure, for example by using copulas. We namely need to simulate discretized approximations of all of these $N$ Brownian paths. Subsequently, we need to cut the Brownian path at some finite time point. We then record the largest observations of all of these paths. Subsequently, we compute the maximum of $N$ of these records to obtain one observation of a maximum queue length. Afterwards, we need to repeat this procedure to collect data. Finally, we use the collected data to compute empirical means and to estimate quantile functions. This means that the computation time grows with at least $N$, the size of the fork-join queue. Besides, in this simulation procedure, a lot of discretization and approximation steps are needed, which increase the error. Though the simulation results give a clear indication of the convergence rate of our limit theorem for small fork-join queueing networks, clearly the procedure above is unworkable for a system with a number of servers of the order of thousands, which as a matter of fact shows the usefulness of the limit in Theorem \ref{thm: second order convergence process} as an approximation.

\paragraph{\textbf{\emph{Summary of results.}}}
In this paper, we study an assembly system with $N$ components, where the demand and the number of produced components are deterministic with some random perturbation, which is assumed to be normally distributed. Thus, the total delay for one component in steady state can be modeled by the all-time supremum of a Brownian motion. We model the system as a fork-join queue. We then use results from EVT to estimate the longest queue, and we minimize the total costs in the system using this approximation, cf.\ Theorems \ref{thm: order bounds indep}, \ref{thm: second order convergence process}, and \ref{thm: order bounds dep} for the most important results.

\paragraph{\textbf{\emph{New insights.}}}This paper generates new insights in fork-join queues that lead to new analytical results for an important class of assembly systems. This paper is the first to consider simultaneous optimization of inventory and capacity in a multi-component assembly system with dependent delays. Due to the dependencies in delays, evaluating such a system with fixed capacity and inventory is already a difficult problem.
We provide several asymptotically optimal expressions for capacity and inventory that are either in closed-form or can easily be computed numerically. Our results may help OEMs to optimally allocate budget to capacity and inventory, to cost-efficiently ensure timely deliveries to their customers. 

\paragraph{\textbf{\emph{Overview.}}}The remainder of this paper is organized as follows. In Section \ref{sec: literature review}, we provide an overview of relevant literature. We introduce the general mathematical model in Section \ref{sec: problem formulation} and subsequently present the optimization problem where we need to decide on capacity and inventory to minimize costs. We study the assembly system with deterministic demand in Section \ref{sec: basic model - independent}. We provide explicit expressions and approximations for optimal inventory and capacity. The stochastic demand case, with solutions to the minimization problem and convergence results, is studied in more detail in Section \ref{sec: model - dependent}. A refinement of the approximations from Section \ref{sec: model - dependent} is provided in Section \ref{sec: master formula}, where we combine the lessons learnt in Sections \ref{sec: basic model - independent} and \ref{sec: model - dependent} to obtain better approximations for optimal capacity and inventory. In Section \ref{sec: asymmetry}, we briefly touch upon the case of asymmetric systems and demonstrate that even in these settings our result for symmetric systems remain useful. We give a summary and conclusions in Section \ref{sec: conclusions} and provide most of the proofs in Appendix \ref{sec: appendix}.

\section{Literature Review}\label{sec: literature review}

Simultaneous optimization of capacity and inventory is an important problem in supply chain management, but the literature on this topic is limited due to complexity of the problem \citep{bradley2002}. 
Considering the interaction between a manufacturer and a single supplier, \cite{chaturvedi2016} discuss the trade-off between inventory and capacity and how properly diversifying supply sources can reduce inventory and capacity investments.
\cite{sleptchenko2003} study simultaneous optimization of spare-part inventory and repair capacity. In the last decade, simultaneous optimization of capacity and inventory in a single supplier-manufacturer relationship has been studied increasingly \cite[e.g.][]{reed2017,reddy2020}. \cite{reed2017} show that the square-root staffing rule of \cite{halfin1981heavy} is a valuable tool in optimizing inventory and capacity in a multi-server make-to-stock queue.
\cite{altendorfer2011} study simultaneous optimization of inventory and planned lead-time and \cite{mayorga2011} study the joint optimization of inventory and temporarily available additional capacity. 
Our work differs fundamentally from these studies, as we consider the assembly of multiple components that face the same (stochastic) demand. 

In particular, we derive extreme value results for multi-component assembly systems as the number of components grows large, in order to obtain asymptotically optimal capacity and inventory decisions. We are not aware of related studies of extreme values for inventory and capacity optimization, but the approach is conceptually related to studies that apply asymptotic analysis to analyze inventory control problems, and we next review this literature. Such studies typically analyze inventory models that are inherently high-dimensional: Asymptotic analysis may be used to derive much simpler optimization problems that form an accurate approximation in some relevant asymptotic regime. This approach has led to major progress in the analysis of inventory problems, e.g. for lost-sales models \citep{goldberg2016asymptotic,xin2016}, dual sourcing \citep{xin2018asymptotic}, and assembly-to-order systems \citep{reiman2015,dougru2017} in the presence of large leadtimes. Assemble-to-order systems with high-volume demand are studied by \cite{plambeck2008} and \cite{plambeckward2008}, while \cite{zhang2020simple} study policies for managing perishable inventory when the market size grows large. 
A comprehensive overview of advances using asymptotic analysis can be found in \cite{goldberg2021}. 
While conceptually related, our analysis differs substantially since a queueing model rather than an MDP underlies our problem, and we aim to analyze extremes in the queueing model to optimize certain model parameters. In that sense, our work is related to \cite{glasserman1997}, who provides approximations for setting base-stock levels in single-stage and multi-stage systems that are asymptotically exact as the target service level or the backorder penalty becomes large. 
For single-product lost-sales inventory systems under periodic review, \cite{huh2009} show that order-up-to policies are asymptotically optimal when the lost sales penalty is large compared to the holding cost. \cite{bijvank2014} show the robustness of this result when using the optimal base-stock levels of the corresponding backorder system instead of those of the lost-sales system.
The asymptotic analysis in this paper has also been influenced by related problems for queues with many servers, inspired by agent staffing problems in call centers; we refer to \cite{borst2004dimensioning, gans2003telephone} and \cite{van2019economies} for background.

Brownian motion models are common in the literature on inventory control. Optimal control of inventory that can be described by a Brownian motion is described by \citet[\S7]{harrison2013}, who provides optimality conditions for both discounted and average cost criteria. Closely related to our work is the Brownian Motion Model presented by \citet[\S3]{bradley2002} to study the trade-off between capacity and inventory. They provide closed-form approximations to the optimal capacity and base-stock levels in a system with a single item. 
We consider an assembly system in which multiple components are merged into one end-product. This is an essential difference, since in our model inventory does not only buffer against uncertain demand, but a component may also need to be stored when other components are not yet available.

We note that our study focuses on the common components of a single high-tech system, which is a considerably simpler problem than general assemble-to-order problems \citep[cf.][]{atan2017}. Our focus enables us to obtain results for the key trade-off between capacity, inventory and delivery reliability, while sidestepping the difficulties of inventory control in multi-product assemble-to-order systems with component commonality \cite[see e.g.][]{song1998,lu2005,reiman2015,atan2017}. 

Literature concerning simultaneous optimization of capacity and inventory in single-sourced assembly (or assembly-to-order) systems with multiple components is limited. \cite{zou2004} study how supply chain efficiency can be increased by synchronizing processing times and delivery quantities. 
\cite{pan2016} consider the simultaneous optimization of component prices and production quantities in a two-supplier setting where one supplier has uncertainty in the yield.
Our main contribution compared to the work of \cite{zou2004} and \cite{pan2016} is that we provide approximations of the optimal capacity and base-stock levels that only require two moments.

To analyze the problem at hand, we examine fork-join queueing networks with $N$ servers where the arrival and service streams are almost deterministic with a Brownian component. Our goal is to find and investigate the maximum queue length as $N$ goes to infinity. The queue lengths are dependent random variables due the joint interarrivals. Thus, our paper is related to the convergence of extreme values (maximum queue lengths) of dependent random variables. An overview of early results on extreme value theory for dependent random variables is given in \cite{leadbetter2012extremes}. The authors provide conditions when the sequence of random variables may be treated as a sequence of independent random variables; this is the case when the covariance of random variables $X_i$ and $X_j$ decreases when $i$ and $j$ are further apart from each other. They also present a convergence result for the joint all-time suprema of a finite number of dependent stationary processes, they prove in Theorem~11.2.3 that, under some assumptions, the joint all-time suprema of a finite number of dependent stationary processes are mutually independent. This is somewhat related to the problem that we study; however, we do not investigate stationary processes and we only look at the largest of the $N$ all-time suprema, where $N\to\infty$. 

We investigate the extreme values for a sequence of $N$ Brownian motions. To be precise, we examine the joint all-time suprema of $N$ dependent Brownian motions with a negative and linear drift term, when $N$ is large. A lot of work has been done on joint suprema of Brownian motions. For instance, \cite{kou2016first} give the solution of the Laplace transform of joint first passage times in terms of the solution of a partial differential equation, where the Brownian motions are dependent. \cite{debicki2020exact} analyze the tail asymptotics of the all-time suprema of two dependent Brownian motions. The joint suprema of a finite number of Brownian motions is also studied; cf.\  \cite{debicki2015extremes}, where the authors give tail asymptotics of the joint suprema of independent Gaussian processes over a finite time interval. These are just three examples, but the literature is rich with variations around assumptions on independence and dependence or around whether or not drift terms are linear, with joint suprema of two or more than two processes, with suprema over finite and infinite time intervals, and with extensions to other Gaussian processes. In this paper, we specifically examine the maximum of $N$ all-time suprema of dependent Brownian motions. In this respect, the work of  \cite{brown1977extreme} comes the closest to our work. In that paper, the authors study process convergence of the scaled maximum of $N$ independent Brownian motions to a stationary limiting process whose marginals are Gumbel distributed. However, we add to this by considering the maximum of the all-time suprema of $N$ dependent Brownian motions.

Our work also relates to the literature on fork-join queues. Specifically, we study asymptotic results for a fork-join queueing system with $N$ servers. Most exact results on fork-join queues are limited to systems with two service stations; cf.\ \cite{flatto1984two}, \cite{wright1992two}, \cite{baccelli1985two} and \cite{de1988fredholm}. For fork-join queues with more than two servers only approximations of performance measures are given; cf.\  \cite{ko2004response}, \cite{baccelli1989queueing} and \cite{nelson1988approximate}. Most of these papers focus on fork-join queueing systems where the number of servers is finite, while we investigate a fork-join queue where $N$ goes to infinity. Furthermore, in these papers, the focus lies on steady-state distributions and other one-dimensional performance measures. Work on the heavy-traffic
process limit has also been done. For example, \cite{varma1990heavy} derives a heavy-traffic analysis for fork-join queues, and shows weak convergence of several processes, such as the joint queue lengths in front of each server. Furthermore, \cite{nguyen1993processing} proves that various appearing limiting processes are in fact multi-dimensional reflected Brownian motions. \cite{nguyen1994trouble} extends this result to a fork-join queue with multiple job types. Lu and Pang study fork-join networks in \cite{lu2015gaussian,lu2017heavy,lu2017heavy2}. In \cite{lu2015gaussian}, they investigate a fork-join network where each service station has multiple servers under nonexchangeable synchronization and operates in the quality-driven regime. They derive functional central limit theorems for the number of tasks waiting in the waiting buffers for synchronization and for the number of synchronized jobs. In \cite{lu2017heavy}, they extend this analysis to a fork-join network with a fixed number of service stations, each having many servers, where the system operates in the Halfin-Whitt regime. In \cite{lu2017heavy2}, the authors investigate these heavy-traffic limits for a fixed number of infinite-server stations, where services are dependent and could be disrupted. Finally, we mention \cite{atar2012control}, who investigate the control of a fork-join queue in heavy traffic by using feedback procedures.

\section{Model and preliminaries}\label{sec: problem formulation}
The production system of OEMs such as ASML, Philips, or Airbus consists out of roughly two stages: 1) Component production; and 2) assembly/integration of components. This setup is crucial to enable the modular design, production and testing of components, and substantial value is added in both stages. For these reasons system integration is only initiated after customers have committed to purchasing the system. We consider a manufacturing system in which a manufacturer assembles a final product from $N$ common components, where $N$ is a large number, meaning that all components are required whenever a product is assembled. Each component is produced on a single production line that involves highly skilled staff and specialized equipment. 
In anticipation of uncertain demand, an inventory buffer is built up: production continues until a target inventory position is reached, after which production is switched off until the inventory position drops below this target.
Such base-stock policies are widely used for modeling component inventories \citep[e.g.][]{akccay2004joint,bollapragada2004managing,karsten2012inventory}. Also in a high-tech manufacturing environment, where capacity mainly refers to people working in cleanrooms that can be at work or have a day off instead of expensive machines with high start-up costs, such policies are suitable. 
Despite these inventory buffers, random delays may occur in the production process for each of the components.

\paragraph{\textbf{\emph{Model.}}}
We adopt a symmetric continuous-time model and assume that production capacity in every finite time interval is normally distributed, meaning that cumulative production is a Brownian motion with drift. We then look at this system in equilibrium, and find a trade-off between investing in the base-stock buffer, and investing in capacity. 
To efficiently satisfy demand of the end-product, which may either be deterministic or stochastic, we need to decide how much capacity to establish for each component and how many finished components to keep on inventory as a buffer. Even though it is costly to establish capacity and to hold inventory, not being able to satisfy demand gives rise to backorder costs. Therefore, we need to find capacity and inventory levels that minimize total expected costs.  

To analyze the cost-minimization problem, we model this assembly system by a fork-join network of $N$ statistically identical, but possibly correlated queues. Demand is represented by the common arrival process of jobs going to each server and each server, with independent, identical service processes, represents production of a component. The backlog of each component is represented by a queue of jobs that have not been served yet. After completion of a job, the finished component is stored in a warehouse. 
As demand at each server
is driven by a common arrival process, the total inventory of a component including the number of
backlogged components is equal for all components. However, as the service times vary, the division
between the number of finished components and the number of backlogged components may vary
per server. When all servers have a finished component in their warehouse, the end-product can be assembled. This system is visualized in Figure \ref{fig: tikzpicture fork-join queue}.
\begin{figure}[H]
\centering

\begin{tikzpicture}[level distance=3cm,
level 1/.style={sibling distance=2cm}]

\node[] (Root) [red]  (A) {}
    child[line width=0,grow=right,draw opacity=0] {
    node[]  {}    
    child[draw=black,draw opacity=0,->] { node{} }     
    child[draw=black,draw opacity=0,->] { node{} }
    child[draw=black,draw opacity=0,->] { node {} }
}
;
\node [draw=none, shift={(8.5cm,-0.9cm)}] (A) {\vdots};
\draw (8.5cm,-2.25cm) circle [radius=0.5cm] node {$N$};
\draw (8.5cm,0cm) circle [radius=0.5cm] node {2};
\draw (8.5cm,2.25cm) circle [radius=0.5cm] node{1};
\draw (6,3) -- ++(2cm,0) -- ++(0,-1.5cm) -- ++(-2cm,0)node[above,yshift=1.5cm,xshift=1cm]{\text{Backlog of components}};
\foreach \i in {1,...,4}
  \draw (8cm-\i*10pt,3) -- +(0,-1.5cm) ;
  \draw (6,0.75) -- ++(2cm,0) -- ++(0,-1.5cm) -- ++(-2cm,0);
\foreach \i in {1,...,4}
  \draw (8cm-\i*10pt,0.75) -- +(0,-1.5cm);
  
  \draw (6,-1.5) -- ++(2cm,0) -- ++(0,-1.5cm) -- ++(-2cm,0);
\foreach \i in {1,...,4}
  \draw (8cm-\i*10pt,-1.5) -- +(0,-1.5cm);
  \draw[->] (9,-2.25) -- (10,-2.25) node[above]{};
    \draw[->] (9,0) -- (10,0) node[above]{};
      \draw[->] (9,2.25) -- (10,2.25) node[above]{};
    
      \draw (0,0) circle [radius=0.5cm] node{};
      \draw (0.5,0) -- (3,0) node[above,xshift=-0.75cm]{\text{Arrival stream}} node[below,xshift=-0.75cm]{\text{of demand}};
    \draw[->] (3,0) -- (5.5,0);
    \draw (4,2.25)--(4,-2.25);
    \draw[->] (4,2.25) -- (5.5,2.25);
    \draw[->] (4,-2.25) -- (5.5,-2.25);
   \draw(11.5,-0.5)--(11.5,0.5);
    \draw(11.50,-0.5)--(12.5,-0.5)node[above,xshift=-0.5cm,yshift=0.2cm]{\text{}};
    \draw(11.50,0.5)--(12.5,0.5);
    \draw(12.5,-0.5)--(12.5,0.5);
    \draw[dashed](10,0)--(11.5,0);
    \draw[dashed](10,-2.25)--(11.5,-2.25);    
     \draw(11.5,-2.75)--(11.5,-1.75);
    \draw(11.5,-2.75)--(12.5,-2.75)node[above,xshift=-0.5cm,yshift=0.2cm]{\text{}};
    \draw(11.5,-1.75)--(12.5,-1.75);
    \draw(12.5,-2.75)--(12.5,-1.75);
    \draw[dashed](10,2.25)--(11.5,2.25);
       \draw(11.5,1.75)--(11.5,2.75);
    \draw(11.5,1.75)--(12.5,1.75)node[above,xshift=-0.5cm,yshift=0.2cm]{\text{}};
    \draw(11.5,2.75)--(12.5,2.75)node[above,yshift=0.32cm]{\text{Warehouse}};
    \draw(12.5,1.75)--(12.5,2.75);
    \draw[dashed](12,1.75)--(12,0.5);
    \draw[dashed](12,-1.75)--(12,-0.5);
    
    \draw[dashed,->](12.5,0)--(13.5,0);
    
     \draw[dashed] (14,0) circle [radius=0.5cm] node[above,yshift=1cm,xshift=0cm]{\text{Assembly of}}
node[above,yshift=0.5cm,xshift=0cm]{\text{ components}};
\end{tikzpicture}

\caption{Fork-join queue}
\label{fig: tikzpicture fork-join queue}
\end{figure}

\paragraph{\textbf{\emph{Brownian fork-join queue.}}}
We model queue lengths as reflected Brownian motions, following  \cite{harrison1985brownian, abate1987transient}. Other papers using Brownian queues to analyze assembly systems are for example \cite{plambeck2008} and \cite{plambeckward2008}.
\begin{defn}\label{def: brownian fork-join queue}
For all $i\leq N$, the service process at server $i$ is governed by the Brownian motion $\{W_i(t),t\geq 0\}$ with standard deviation $\sigma$, and the arrival process is governed by the Brownian motion $\{W_A(t),t\geq 0\}$ with standard deviation $\sigma_A$. The queue length at server $i$ at time $t>0$ equals
\begin{align}\label{eq: transient queue length}
Q_{i}(t,\beta):=\sup_{0<s<t}((W_i(t)+W_A(t)-\beta t)-(W_i(s)+W_A(s)-\beta s)),
\end{align}
with $Q_{i}(0,\beta)=0$.
For $i,j\leq N$ with $i\neq j$ the Brownian motions $\{W_i(t),t\geq 0\}$ and $\{W_j(t),t\geq 0\}$ are i.i.d.
\end{defn}
Formally, the Brownian motions $\{W_i(t),t\geq 0\}$ and  $\{W_A(t),t\geq 0\}$ represent fluctuations in the service and arrival processes, as they have zero mean. The controllable parameter $\beta$ represents the excess capacity
in each individual queue. 

\paragraph{\textbf{\emph{Base-stock level and capacity.}}}
To buffer against uncertainties in the supply and demand processes, we introduce a base-stock level $I_i$ for each component $i\leq N$.
We define $\beta_i>0$ as the net capacity for component $i$, i.e. the difference between the production rate and arrival rate, in other words, $\beta_i$ captures the capacity investment of server $i$. As mentioned before, we assume that for all servers, the net capacity and the base-stock levels are the same, thus $\beta_i=\beta_j=\beta$ and $I_i=I_j=I$. The backlog $Q_{i}(t,\beta)$ represents the number of outstanding orders of component $i\leq N$ at time $t$, with $Q_{i}(t,\beta)$ given in Definition \ref{def: brownian fork-join queue}. 
If $\sigma_A^2>0$, $(Q_{i}(t,\beta))_{i\leq N}$ are dependent random variables.


\paragraph{\textbf{\emph{Transient inventory levels and backorders.}}}
We proceed by developing an expression for the total system costs, which requires expressions for the inventory and backorders. 
The inventory of component $i$ consists of two parts: first, the excess supply that works as a buffer against uncertain demand; second, the committed inventory that consists of items that are committed to realized demand but put aside because other components are not yet available.
I.e., the excess supply of component $i$ is given by $(I-Q_{i}(t,\beta))^+$. Moreover, the number of backorders for component $i$ at time $t$ is equal to $(Q_{i}(t,\beta)-I)^+$, since for $Q_{i}(t,\beta)\leq I$ the shortage is compensated by inventory $I$ and only the part of $Q_{i}(t,\beta)$ exceeding $I$ represents actual backorders that cannot be satisfied.
Since all components need to be available to assemble the final product, the number of backorders in the system is equal to the number of backorders of the component with the largest backlog and is thus given by $\max_{j\leq N}\left(Q_{j}(t,\beta)-I\right)^+$. Therefore, the committed inventory of component $i$ equals the number of backorders in the system minus its own backlog and can be expressed as
$\max_{j\leq N}\left(Q_{j}(t,\beta)-I\right)^+-(Q_{i}(t,\beta)-I)^+.$
The total inventory of component $i$ at time $t$ is thus given by 
\begin{align}\label{eq: inventory position}
I_{i}(t)=(I-Q_{i}(t,\beta))^++\max_{j\leq N}\left(Q_{j}(t,\beta)-I\right)^+-(Q_{i}(t,\beta)-I)^+
=I-Q_{i}(t,\beta)+\max_{j\leq N}\left(Q_{j}(t,\beta)-I\right)^+,
\end{align}
with $I_{i}(0)=I$.
Observe that the total inventory $I_{i}(t)$ at time $t$ is a function of the number of outstanding orders at time $t$. The reason why this is true, is that the random variable $Q_{i}(t,\beta)$ does not depend on the total inventory, because the servers always produce when there is an incoming task, irrespective whether there are items in stock, or not. When there are items in stock, the product is immediately assembled, but servers work in order to reach the target inventory. When there are no items in stock, servers work to finish their component. Hence, whether or not a server works, does not depend on the total inventory but only on the demand and their own service speed. This means that the total inventory at time $t$ is described as the function given in Equation \eqref{eq: inventory position}. Thus, in order to know the total inventory on a certain time $t$, one should know the number of outstanding orders on that given time $t$, where the dynamics of these outstanding orders are described as the dynamics of reflected Brownian motions until time $t$. Thus, this describes the \textit{dynamics} of the system.

\paragraph{\textbf{\emph{Steady-state limit.}}}
Because the backlogs are modeled as reflected Brownian motions with negative drift, the backlogs have a steady-state limit. This limit extends to the largest backlog in the system and the total inventory of component $i$. We prove this in Lemma \ref{lem: steady-state limit delay}.
\begin{lem}[Steady-state of backlogs.]\label{lem: steady-state limit delay}
Given $(Q_{i}(t,\beta),i\leq N)$ with $Q_{i}(t,\beta)$ defined in \eqref{eq: transient queue length}, we have that
$(Q_{i}(t,\beta),i\leq N) \stackrel{d}{\rightarrow} (Q_{i}(\infty,\beta),i\leq N)$ with
\begin{align}\label{eq: steady state limit max q}
(Q_{i}(\infty,\beta), i\leq N)\overset{d}{=}(\sup_{s>0}(W_i(s)+W_A(s)-\beta s), i\leq N).
\end{align}
In particular, 
\begin{align}\label{eq: steady state limit max q}
\max_{i\leq N}Q_{i}(\infty,\beta)\overset{d}{=}\max_{i\leq N}\sup_{s>0}(W_i(s)+W_A(s)-\beta s).
\end{align}
\end{lem}
\proof{Proof.} The argument in one dimension is standard (see e.g. Section III.6 of \cite{asmussen2003applied}); we extend it to our setting. 
Given $t>0$, we can define Brownian motions $\{\hat{W}_A(s),s\geq 0\}$ and $\{\hat{W}_i(s),s\geq 0\}$ that satisfy $\hat{W}_A(t-s)=W_A(t)-W_A(s)$ and $\hat{W}_i(t-s)=W_i(t)-W_i(s)$. From this, it follows that for fixed $t>0$, we have that
\begin{align*}
   (Q_{i}(t,\beta), i\leq N) &=(\sup_{0\leq s\leq t}(\hat{W}_i(s)+\hat{W}_A(s)-\beta s), i\leq N)\\
   &\overset{d}{=} (\sup_{0\leq s\leq t}(W_i(s)+W_A(s)-\beta s), i\leq N). 
\end{align*}
Now, we obtain the lemma by letting $t\rightarrow\infty$, using monotone convergence.
\hfill\rlap{\hspace*{-2em}\Halmos}
\endproof 
Combining this result with (\ref{eq: inventory position}), we obtain an analogous result for the steady-state total inventory. In particular,
\begin{align*}
    \sum_{i=1}^N I_i(t) \stackrel{d}{\rightarrow} \sum_{i=1}^N  (I-Q_{i}(\infty,\beta)+\max_{j\leq N}\left(Q_{j}(\infty,\beta)-I\right)^+).
\end{align*}
From now on, we write $Q_i(\beta):=Q_{i}(\infty,\beta)$.

\paragraph{\textbf{\emph{Cost function.}}}
We scale the cost of building net capacity to one and let $h^{(N)}$ and $b^{(N)}$ denote (inventory) holding costs and backorder costs, respectively, which may depend on $N$. Our goal is to minimize the expected total costs of the system in steady state.
\begin{defn}
We define
\begin{align}\label{eq: cost function}
    C_N(I,\beta):=\expect*{\sum_{i\leq N}\left[h^{(N)} \left(I-Q_i(\beta)+\max_{j\leq N}(Q_j(\beta)-I)^+\right)\right]+b^{(N)}\max_{j\leq N}(Q_j(\beta)-I)^+},
\end{align}
with the distribution of $Q_i(\beta)$ given in Equation \eqref{eq: steady state limit max q}.
\end{defn}
Equation \eqref{eq: cost function} simplifies to
$$
    C_N(I,\beta) =\expect*{Nh^{(N)} (I-Q_i(\beta))+(Nh^{(N)}+b^{(N)})\left(\max_{j\leq N}Q_j(\beta)-I\right)^+}.
$$
Then, the expected total costs in the system are equal to $C_N(I,\beta)+\beta N$, where the term $\beta N$ reflects our normalization of unity net capacity costs per queue. If this term would be removed, it would be optimal to choose $\beta=\infty$ and $I=0$. 

Due to the self-similarity of Brownian motion, we can write 
\begin{align*}
   \beta \max_{i\leq N}\sup_{s>0}(W_A(s)+W_i(s)-\beta s)
   &=\beta \max_{i\leq N}\sup_{t>0}\left(W_A\left(\frac{t}{\beta^2}\right)  + W_i\left(\frac{t}{\beta^2}\right)-\beta \frac{t}{\beta^2}\right)\\
   &\overset{d}{=}\max_{i\leq N}\sup_{t>0}(W_A(t)+W_i(t)-t).
\end{align*}
This means that $\max_{i\leq N} Q_i(\beta)\overset{d}{=}\frac{1}{\beta}\max_{i\leq N}Q_i(1)$. Therefore, after rescaling the variable $I$, we can write
\begin{align}\label{eq: split beta and I minimization}
     \min_{(I,\beta)}\bigg(C_N(I,\beta)+\beta N\bigg)=\min_{(I,\beta)}\bigg(\frac{1}{\beta}C_N(I\beta,1)+\beta N\bigg)=\min_{(I,\beta)}\bigg(\frac{1}{\beta}C_N(I,1)+\beta N\bigg). 
\end{align}
In the last part of Equation \eqref{eq: split beta and I minimization}, $I$ has the interpretation of the base-stock level where the net capacity $\beta=1$. Therefore, from now on, the actual number of products on stock at time 0 equals $I/\beta$. Similarly, the actual unsatisfied demands of component $i$ equals $Q_i(1)/\beta$ and we write $Q_i=Q_i(1)$. This allows us to write the cost function $F_N(I,\beta)$ to be optimized  as given in Definition \ref{def: optimalcosts independent}.
\begin{defn}\label{def: optimalcosts independent}
We define 
\begin{align}
    F_N(I,\beta):= C_N(I,\beta) +\beta N=\frac{1}{\beta}C_N(I) +\beta N,
\end{align}
with $C_N(I):=C_N(I,1)$ and $C_N(I,\beta)$ given in Equation \eqref{eq: cost function}. 
\end{defn}

Our goal is to solve $\min_{(I,\beta)}F_N(I,\beta)$, focusing on the case where $N$ is large.

\subsection{Preliminary results}
As we have defined the Brownian fork-join queue and the corresponding cost functions, we now state some general results that are valid regardless whether $\sigma_A=0$ or $\sigma_A>0$.
In the next lemma, 
we show that we can write $\min_{(I,\beta)}F_N(I,\beta)$ as two separate minimization problems.

\begin{lem}\label{lem: inventory minimization}
Let $(b^{(N)})_{N\geq 1}, (h^{(N)})_{N\geq 1}$ be sequences such that $h^{(N)}>0$ and $b^{(N)}>0$ for all $N$. Let $\left(I_N,\beta_N\right)$ minimize $F_N(I,\beta)$. Then the optimal base-stock level $I_N$ minimizes $C_N(I)$ and the optimal $\beta_N$ minimizes $\frac{1}{\beta}C_N(I_N)+\beta N$. Furthermore, the function $C_N(I)$ is convex with respect to $I$, and the function $\frac{1}{\beta}C_N(I)+\beta N$ is convex with respect to $\beta$.
\end{lem}

Using Lemma \ref{lem: inventory minimization}, we can characterize the optimal net capacity and base-stock level. In Lemma \ref{lem: optimal beta general} we provide expressions for the optimal net capacity and costs in terms of the optimal base-stock level, which is given in Lemma \ref{lem: optimal I general}.
\begin{lem}\label{lem: optimal beta general}
Given $I^*_N= \arg\min_I C_N(I)$, minimizing $F_N(I,\beta)$ with respect to $\beta$ yields $\beta^*_N=\sqrt{\frac{C_N(I^*_N)}{N}}$. Furthermore, the corresponding costs are $F_N(I^*_N,\beta^*_N)=2N\beta^*_N=2\sqrt{C_N(I^*_N)N}$.
\end{lem}
The optimal value of $I$ can be expressed as a quantile of the distribution of 
$\max_{i\leq N} Q_i$:
\begin{lem}\label{lem: optimal I general}
$I^*_N$ is the unique solution of 
\begin{align*}
    \probability*{\max_{i\leq N} Q_i \leq I_N^*} = \frac{b^{(N)} }{N h^{(N)}+b^{(N)}}.
\end{align*}
\end{lem}
The main technical issue is that the distribution of this maximum is in general not very tractable, especially when $N$ is large. The main theme of our work is to consider approximations of this distribution using extreme value theory, to analyze their quality if $N$ is large.  

To explain our ideas, we mention the following first-order approximation of $\max_{i\leq N} Q_i$: 
\begin{lem} \label{lem: first order EVT}
$\max_{i\leq N} Q_i$ satisfies the first-order approximation
\begin{align*}
    \frac{\max_{i\leq N} Q_i}{ \log N} \LimitL \frac{\sigma^2}{2},
\end{align*}
as $N\to\infty$.
\end{lem}
The lemma easily follows from more refined results that are proven later on in this paper.

This first-order approximation is valid regardless whether $\sigma_A=0$ or $\sigma_A>0$.
In the subsequent two sections, we consider more refined extreme-value theory approximations covering both cases. It turns out that the second-order behavior of the maximum is qualitatively different when $\sigma_A$ becomes strictly positive. This has, in turn, an impact on the structure of the optimal solution of our cost minimization problem when $N$ grows large. 

To better understand this structure, we heuristically analyze the first-order approximation of the cost minimization problem and apply it to approximate $I_N^*$ and $\beta_N^*$. 
First, we use the approximation $\max_{i\leq N} Q_i \approx\frac{\sigma^2}{2} \log N$
to write
\begin{align*}
    C_N(I)\approx \bar C_N(I) = N h^{(N)}\left(I-\frac{\sigma^2+\sigma_A^2}{2}\right) + (N h^{(N)}+b^{(N)})\left(\frac{\sigma^2}{2} \log N-I\right)^+.
\end{align*}
The optimal value $\bar I_N$ for the associated first-order minimization problem 
$\min_I \bar C_N(I)$ is given by  $\bar I_N = \frac{\sigma^2}{2} \log N$, since $b^{(N)}>0$. 
Using this approximation, we see that $C_N(\bar I_N) \approx \bar C_N(\bar I_N)=(1+o(1)) \frac{\sigma^2}{2} N h^{(N)} \log N$, $\bar \beta_N = \sqrt{\bar C_N(\bar I_N)/N}= (1+o(1)) \sqrt{ \frac{\sigma^2}{2} h^{(N)} \log N}$,
and $F_N(\bar I_N, \bar \beta_N) \approx 2 \sqrt{N} \sqrt{\frac{\sigma^2}{2} N h^{(N)}\log N}$.
These results can be made rigorous and the decision rule $\bar I_N$ can be shown to be asymptotically optimal, i.e.\ that $F_N(\bar I_N, \bar \beta_N ) = F_N(I^*_N, \beta^*_N)(1+o(1))$.
To prove this, we need to specify how the cost parameters $h^{(N)}$ and $b^{(N)}$ scale with $N$. For this, we consider three regimes. These regimes relate to the quantile $b^{(N)}/ (N h^{(N)}+b^{(N)})$ of $\max_i Q_i$ at which $I_N^*$ attains its optimal solution. Assume that $b^{(N)}/ (N h^{(N)}+b^{(N)})$ converges to a constant $1-\gamma$. We classify the three regimes in a similar way as is done
in the analysis of large call centers; cf.\ \cite{borst2004dimensioning}:
\begin{itemize}
\item  We are in the {\em balanced regime} if $\gamma\in (0,1)$.
\item If $\gamma=0$, for large systems, the inventory is always sufficiently high to ensure that the manufacturer can assemble the end-product. We call this the {\em quality-driven regime}.
\item Finally, if $\gamma=1$, inventories are much lower, and we call this the {\em efficiency-driven regime}.
\end{itemize}
When we are in the balanced or efficiency-driven regime we can prove how far the costs under the first order approximation are from the real optimal costs. This is established in Lemma \ref{lem: first order bal reg}:
\begin{lem}\label{lem: first order bal reg}
Assume $\gamma_N=Nh^{(N)}/(Nh^{(N)}+b^{(N)})$, with $\gamma_N=\gamma\in(0,1)$ or $\gamma_N\LimitN 1$. Then
\begin{align*}
    \frac{F_N(I^*_N,\beta^*_N)}{F_N(\bar I_N, \bar \beta_N)}=1-o(1).
\end{align*}
\end{lem}

In the next two sections, we carry out a more elaborate program using more refined extreme value estimates of $\max_{i\leq N} Q_i$. This analysis gives sharper order bounds than those given in Lemma \ref{lem: first order bal reg}.
In particular, in the following sections we consider the minimization in two distinct cases. First, in Section \ref{sec: basic model - independent}, we look at the case where demand is assumed to be deterministic, such that $W_A=0$. Thereafter, in Section \ref{sec: model - dependent}, we consider the stochastic demand case. In the former case, we utilize existing results in extreme value theory, while the latter case requires the development of a novel limit theorem. Furthermore, we use the result given in Corollary \ref{cor: ratio optimal approx}; this corollary shows how the ratio between the optimal costs and approximate costs can be represented, when the approximate base-stock level and net capacity are solutions to a minimization problem as well. This corollary follows trivially from Lemma \ref{lem: optimal beta general}.

\begin{cor}\label{cor: ratio optimal approx}
Assume we have a function $\tilde{F}_N(I,\beta): (0,\infty)\times(0,\infty)\rightarrow\mathbb{R}$. Furthermore, assume that the function $\tilde{F}_N$ has the form  
\begin{align*}
    \Tilde{F}_N(I,\beta)=\frac{1}{\beta}\Tilde{C}_N(I)+\beta N,
\end{align*}
where $\Tilde{C}_N$ is a positive function with domain $(0,\infty)$. Moreover, assume that the minimum value $\Tilde{F}_N(\tilde{I}_N,\tilde{\beta}_N)=2N\tilde{\beta}_N=2\sqrt{\Tilde{C}_N(\tilde{I}_N)N}$, where $\tilde{I}_N$ and $\tilde{\beta}_N$ are minimizers, then
\begin{align*}
    \frac{F(I^*_N,\beta^*_N)}{F(\tilde{I}_N,\tilde{\beta}_N)}=\frac{2\sqrt{C_N(I^*_N)}\sqrt{\Tilde{C}_N(\tilde{I}_N)}}{C_N(\tilde{I}_N)+\Tilde{C}_N(\tilde{I}_N)}.
\end{align*}
\end{cor}

\section{The basic model: deterministic arrival stream}\label{sec: basic model - independent}

In this section, we consider the case that demand is deterministic. From this, it follows that all $N$ queues are mutually independent.

\subsection{Solution and convergence of the minimization problem}\label{subsec: solution independent}
We now analyze the minimization of the cost function described in Definition \ref{def: optimalcosts independent} for the special case with $W_A=0$ representing deterministic demand.
Although we can simplify the minimization problem significantly, by using the self-similarity of Brownian motions and by writing the minimization problem as two separate minimization problems, as shown in Lemma \ref{lem: inventory minimization}, the function $F_N$ still has a difficult form, since we have the expression $\max_{i\leq N}Q_i$ in this function. 
In Lemma \ref{lem: minimization problem case WA=0} we give the optimal base-stock level in order to minimize costs. We assume that the holding and backlog costs $h^{(N)}$ and $b^{(N)}$ are positive sequences, and we distinguish three cases. First of all, we consider the balanced regime $\gamma_N=Nh^{(N)}/(Nh^{(N)}+b^{(N)})=\gamma\in(0,1)$ for all $N>0$. Secondly, we consider the quality driven regime, where $\gamma_N\LimitN 0$. Finally, we investigate the efficiency driven regime, where $\gamma_N\LimitN 1$. All proofs for this section can be found in Appendix \ref{app: proofs basic model - independent}. We present numerical results for the three regimes in Section \ref{subsec: numerical results independent}. 

\begin{lem}\label{lem: minimization problem case WA=0}
Let $Q_i=\sup_{s>0}(W_i(s)-s)$, with $(W_i,1\leq i\leq N)$ independent Brownian motions with mean 0 and variance $\sigma^2$. Let $h^{(N)}$ and $b^{(N)}$ be positive sequences.
In order to minimize $F_N(I,\beta)$, the optimal base-stock level $I^*_N$ satisfies,
\begin{align}\label{eq: optimal inventory case WA=0}
    I^*_N=P^{-1}_N\left(1-\gamma_N\right)=\frac{\sigma^2}{2}\log\left(\frac{1}{1-\left(1-\gamma_N\right)^{\frac{1}{N}}}\right),
\end{align}
with $P^{-1}_N$ the quantile function of  $\probability*{\max_{i\leq N}Q_i<x}$ and $\gamma_N=Nh^{(N)}/(Nh^{(N)}+b^{(N)})$.
\end{lem}

To get a better understanding of the limiting behavior of the solution to $\min_{(I,\beta)}F_N(I,\beta)$, we would like to approximate the function $F_N$. Since $(Q_i,i\leq N)$ are independent and exponentially distributed, we know by standard extreme value theory (cf.\ \cite{de2007extreme}) that $\frac{2}{\sigma^2}\max_{i\leq N} Q_i-\log N\LimitD G$, as $N\to\infty$, with $G\sim\text{Gumbel}$. Therefore, for $N$ large, $\max_{i\leq N}Q_i\overset{d}{\approx}\frac{\sigma^2}{2}G+\frac{\sigma^2}{2}\log N$. We get a new minimization problem when we replace $\max_{i\leq N} Q_i$ with this approximation $\frac{\sigma^2}{2}G+\frac{\sigma^2}{2}\log N$. In Definition \ref{def: approxcosts independent} we give the resulting function $\ApproxCosts(I,\beta)$ that is to be minimized.
\begin{defn}\label{def: approxcosts independent}
\begin{align}
  \ApproxC{I}:=\mathbb{E}\left[Nh^{(N)}\left(I-Q_i\right)+\bigg(Nh^{(N)}+b^{(N)}\bigg)\left(\frac{\sigma^2}{2}G+\frac{\sigma^2}{2}\log N-I\right)^+\right],
\end{align}
and
\begin{align}\label{eq: approximate problem independent}
    \ApproxCosts(I,\beta):=  \frac{1}{\beta}\ApproxC{I}+\beta N.
\end{align}
\end{defn}
In the remainder of this section, we investigate whether minimizing $\ApproxCosts(I,\beta)$ results in costs that are close to those when we minimize $F_N(I,\beta)$. Note that we write $(I^*_N,\beta^*_N)$ for the minimizers of the cost function $F_N$ defined in Definition \ref{def: optimalcosts independent}, and we write $(\ApproxIN,\ApproxBetaN)$ for the minimizers of the cost function $\ApproxCosts$ defined in Definition \ref{def: approxcosts independent}. Throughout this paper, we indicate second-order approximations by the $\wedge$-symbol.

In Proposition \ref{prop: approx case WA=0}, we present the base-stock level that minimizes $\ApproxCosts$. This base-stock level turns out to be a quantile of $\frac{\sigma^2}{2}G$  added to $\frac{\sigma^2}{2}\log N$.

\begin{prop}[Approximation]\label{prop: approx case WA=0}
Minimizing $\ApproxCosts(I,\beta)$ with $G\sim\text{Gumbel}$, gives solution $(\ApproxIN,\ApproxBetaN,\ApproxCosts(\ApproxIN,\ApproxBetaN))$, with 
\begin{align}\label{eq: approx inventory case WA=0}
    \ApproxIN=\frac{\sigma^2}{2}\log N-\frac{\sigma^2}{2}\log\left(-\log\left(1-\gamma_N\right)\right),
\end{align}
and
\begin{align}\label{eq: approx capacity case WA=0}
\ApproxC{\ApproxIN}=&Nh^{(N)}\bigg(\ApproxIN-\frac{\sigma^2}{2}\bigg)
+(Nh^{(N)}+b^{(N)})\frac{\sigma^2}{2}\left( \int_{-\log \left(1-\gamma_N\right)}^{\infty}\frac{e^{-t}}{t}dt+\Gamma   +\log\left(-\log\left(1-\gamma_N\right)\right)\right),
\end{align}
where $\Gamma\approx 0.577$ is Euler's constant and $\gamma_N=Nh^{(N)}/(Nh^{(N)}+b^{(N)})$.
\end{prop}

Combining Equations \eqref{eq: approx inventory case WA=0} and \eqref{eq: approx capacity case WA=0} with the results in Lemma \ref{lem: optimal beta general} gives the solution $(\ApproxIN,\ApproxBetaN,\ApproxCosts(\ApproxIN,\ApproxBetaN))$.

We compare the costs under the optimal base-stock level and net capacity with the costs under the approximate base-stock level and net capacity. We distinguish the balanced regime, quality driven regime and efficiency driven regime. 

By using the results from Lemmas \ref{lem: G and Q same prob space} and \ref{lem: bounds C Ctilde} in Appendix \ref{app: proofs basic model - independent}, we prove the order bounds in the balanced, quality driven and efficiency driven regime in Theorem \ref{thm: order bounds indep}. In the efficiency driven regime, we impose the additional condition $\gamma_N<1-\exp(-N)$ needed to make sure that $\ApproxIN>0$. If we namely choose $\gamma_N>1-\exp(-N)$, we get that $\ApproxIN<0$, which is not feasible, because $\ApproxIN$ has the physical meaning of the number of items that needs to be stored.

\begin{thm}[Order bounds]\label{thm: order bounds indep}
Assume $\gamma_N=Nh^{(N)}/(Nh^{(N)}+b^{(N)})$, if $\gamma_N=\gamma\in(0,1)$, in the \textit{balanced regime}, then
\begin{align}
    \frac{F_N(I^*_N,\beta^*_N)}{F_N(\ApproxIN,\ApproxBetaN)}=1-O(1/(N\log N)),
\end{align}
if $\gamma_N\LimitN 0$, in the \textit{quality driven regime}, then
\begin{align}
     \frac{F_N(I^*_N,\beta^*_N)}{F_N(\ApproxIN,\ApproxBetaN)}=1-O(\gamma_N/(N\log(N/\gamma_N))),
\end{align}
and if $\gamma_N\LimitN 1$ and $\gamma_N<1-\exp(-N)$, in the \textit{efficiency driven regime}, then
\begin{align}
  \frac{F_N(I^*_N,\beta^*_N)}{F_N(\ApproxIN,\ApproxBetaN)}=1-O(1/\log N).
\end{align}
\end{thm}
Using the order bounds given in Theorem \ref{thm: order bounds indep}, we can establish for the three different regimes how $F_N(I^*_N,\beta^*_N)$ scales with $N$ as $N$ becomes large.
\begin{lem}\label{lem: magnitude F}
Assume $\gamma_N=Nh^{(N)}/(Nh^{(N)}+b^{(N)})$, if $\gamma_N=\gamma\in(0,1)$ in the \textit{balanced regime}, then
\begin{align}
      &F_N(I^*_N,\beta^*_N)\nonumber\\
      =&2\sqrt{N}\sqrt{Nh^{(N)}\frac{\sigma^2}{2}\left(\log N-\log(-\log(1-\gamma))-1\right)+(Nh^{(N)}+b^{(N)})\frac{\sigma^2}{2}\expect*{\left(G+\log(-\log(1-\gamma))\right)^+}}\nonumber\\
      +&O(\sqrt{h^{(N)}}/\sqrt{\log N}),
\end{align}
if $\gamma_N\LimitN 0$ in the \textit{quality driven regime}, then
\begin{align}
     F_N(I^*_N,\beta^*_N)
      =&2\sqrt{N}\sqrt{Nh^{(N)}\frac{\sigma^2}{2}\left(\log (N/\gamma_N)-1\right)+(Nh^{(N)}+b^{(N)})\frac{\sigma^2}{2}\gamma_N}\nonumber\\
      +&O(\gamma_N\sqrt{h^{(N)}}/\sqrt{\log (N/\gamma_N)}),
\end{align}
and if $\gamma_N\LimitN 1$ and $\gamma_N<1-\exp(-N)$ in the \textit{efficiency driven regime}, then
\begin{align}
F_N(I^*_N,\beta^*_N)
      =2\sqrt{N}\sqrt{Nh^{(N)}\frac{\sigma^2}{2}\left(\log N-1\right)+b^{(N)}\frac{\sigma^2}{2}\log(-\log(1-\gamma_N))}
      +O(N\sqrt{h^{(N)}}/\sqrt{\log N}).
\end{align}
\end{lem}
The results given in Theorem \ref{thm: order bounds indep} and Lemma \ref{lem: magnitude F} are obtained by using the properties stated in Lemmas \ref{lem: G and Q same prob space} and \ref{lem: bounds C Ctilde}. In Lemma \ref{lem: G and Q same prob space} we show that we can write a Gumbel distributed random variable that is on the same probability space as $\max_{i\leq N}Q_i$. This gives us a very powerful result; namely that $\max_{i\leq N}Q_i$ and $G_N$ are ordered and that their difference decreases as $\max_{i\leq N}Q_i$ becomes large. Consequently, we obtain very sharp bounds on $|C_N(I^*_N)-C_N(\ApproxIN)|$ and $|\ApproxC{\ApproxIN}-C_N(\ApproxIN)|$ in Lemma \ref{lem: bounds C Ctilde} which leads to sharp results in Theorem \ref{thm: order bounds indep} and Lemma \ref{lem: magnitude F}.

\subsection{Numerical experiments}\label{subsec: numerical results independent}
We now provide some numerical results to illustrate the solutions to the minimization problem and their characteristics discussed in Section \ref{subsec: solution independent}. In all experiments, we let $\sigma=1$ and let $N$ vary from 10 to 1000. The results for the balanced regime, quality driven regime and efficiency driven regime are given in Tables \ref{tab: numerical balanced regime}, \ref{tab: numerical quality driven regime} and \ref{tab: numerical efficiency driven regime without condition}, respectively. We can observe that in all regimes the approximate solutions are close to the optimal solutions. Most importantly, already for small $N$, the fraction of the costs corresponding to the optimal solution over the costs corresponding to the approximate solution nearly equals 1.

\begin{table}[H]
\caption{Balanced Regime, $h^{(N)}=1,b^{(N)}=N$ such that $\gamma_N=\frac{1}{2}$.}
\label{tab: numerical balanced regime}
\begin{tabular}{c|c|c|c|c|c|c|c}
$N$ & $I^*_N$ & $\beta^*_N$ & $F_N(I^*_N,\beta^*_N)$ & $\ApproxIN$ & $\ApproxBetaN$ & $F_N(\ApproxIN,\ApproxBetaN)$ & $\left(1-\frac{F_N(I^*_N,\beta^*_N)}{F_N(\ApproxIN,\ApproxBetaN)}\right)N\log N$\\
\hline
 10 & 1.35178 & 1.19648 & 23.9296 & 1.33455 & 1.19328 & 23.9315  & 0.001807\\
 50 & 2.14273 & 1.49338 & 149.338 & 2.13927 & 1.49286 & 149.338  & 0.000379\\
 100 & 2.48757 & 1.60499 & 320.997 & 2.48584 & 1.60475 & 320.997  & 0.000192\\
 200 & 2.83328 & 1.70944 & 683.775 & 2.83242 & 1.70932 & 683.775  & $9.68\cdot 10^{-5}$\\
 500 & 3.29091 & 1.8385 & 1838.5 & 3.29056 & 1.83846 & 1838.5  & $3.91\cdot 10^{-5}$\\
 1000 & 3.63731 & 1.93044 & 3860.87 & 3.63713 & 1.93042 & 3860.87  & $1.97\cdot 10^{-5}$\\
\end{tabular}
\end{table}

\begin{table}[H]
\caption{Quality Driven Regime, $h^{(N)}=1,b^{(N)}=N^2$ such that $\gamma_N=\frac{1}{1+N}$.}
\label{tab: numerical quality driven regime}
\begin{tabular}{c|c|c|c|c|c|c|c}
$N$ & $I^*_N$ & $\beta^*_N$ & $F_N(I^*_N,\beta^*_N)$ & $\ApproxIN$ & $\ApproxBetaN$ & $F_N(\ApproxIN,\ApproxBetaN)$& $\left(1-\frac{F_N(I^*_N,\beta^*_N)}{F_N(\ApproxIN,\ApproxBetaN)}\right)\frac{N}{\gamma_N}\log \frac{N}{\gamma_N}$\\
\hline
10	&	2.32898	&	1.52962	&	30.5925	&	2.3266	&	1.52924	&	30.5925	 & 0.000617\\
50	&	3.91708	&	1.97978	&	197.978	&	3.91698	&	1.97976	&	197.978	 & $2.52\cdot 10^{-5}$\\
100	&	4.60768	&	2.14684	&	429.368	&	4.60766	&	2.14684	&	429.368	 & $6.31162\cdot 10^{-6}$\\
200	&	5.29957	&	2.30221	&	920.886	&	5.29956	&	2.30221	&	920.886	 & $1.21801\cdot 10^{-6}$\\
500	&	6.21511	&	2.49306	&	2493.06	&	6.21511	&	2.49306	&	2493.06	 & $5.51467\cdot 10^{-6}$\\
1000	&	6.90801	&	2.62833	&	5256.66	&	6.90801	&	2.62833	&	5256.66	 & 0.000176\\
\end{tabular}
\end{table}

\begin{table}[H]
\caption{Efficiency Driven Regime, $h^{(N)}=N,b^{(N)}=1$ such that $\gamma_N=\frac{N^2}{N^2+1}$.}
\label{tab: numerical efficiency driven regime without condition}
\begin{tabular}{c|c|c|c|c|c|c|c}
$N$ & $I^*_N$ & $\beta^*_N$ & $F_N(I^*_N,\beta^*_N)$ & $\ApproxIN$ & $\ApproxBetaN$ & $F_N(\ApproxIN,\ApproxBetaN)$ &$\left(1-\frac{F_N(I^*_N,\beta^*_N)}{F_N(\ApproxIN,\ApproxBetaN)}\right)\log N$ \\
\hline
10	&	0.497572	&	3.12224	&	62.4448	&	0.386624	&	3.08439	&	62.4616	 & 0.000797\\
50	&	0.965997	&	9.35451	&	935.451	&	0.927385	&	9.34122	&	935.452	 & $8.65678\cdot 10^{-6}$ \\
100	&	1.21527	&	14.4701	&	2894.02	&	1.19242	&	14.4615	&	2894.02	 & $1.30518\cdot 10^{-6}$\\
200	&	1.48208	&	22.0864	&	8834.57	&	1.46889	&	22.0808	&	8834.57	 & $2.20863\cdot 10^{-7}$\\
500	&	1.85348	&	38.0553	&	38055.3	&	1.84728	&	38.0521	&	38055.3	 & $2.51171\cdot 10^{-8}$\\
1000	&	2.14443	&	56.945	&	113890	&	2.14098	&	56.9428	&	113890	 & $5.30189\cdot 10^{-9}$\\
\end{tabular}
\end{table}
\section{Stochastic Demand}\label{sec: model - dependent}
We now extend our framework to the case where demand is stochastic. This means that stochasticity not only arises from the production process of the individual components, but also results from uncertain demands. Consequently, delays may no longer only be caused by low production of a specific component, but may also occur when there is a sudden peak in demand. Since all components need to be available to assemble the end-product and satisfy demand, delays of the different components are now correlated. We use the same strategy when demand is stochastic as in the basic model with deterministic demand. However, we can no longer approximate the maximum queue length distribution with the Gumbel distribution. In Section \ref{subsec: second order convergence} we show that for $N$ large, $\max_{i\leq N}Q_i\approx \frac{\sigma^2}{2}\log N+\frac{\sigma\sigma_A}{\sqrt{2}}\sqrt{\log N}X$ with $X$ a standard normal random variable. Using this approximation, we obtain a new minimization problem, in which we minimize $\ApproxWACosts(I,\beta)$ as given in Definition \ref{def: approxcosts dependent} with respect to $I$ and $\beta$.
\begin{defn}\label{def: approxcosts dependent}
\begin{align*}
    \ApproxWAC{I}=\mathbb{E}\left[Nh^{(N)}\left(I-Q_i\right)+\bigg(Nh^{(N)}+b^{(N)}\bigg)\left(\frac{\sigma^2}{2}\log N+\frac{\sigma\sigma_A}{\sqrt{2}}\sqrt{\log N}X-I\right)^+\right],
\end{align*}
and
\begin{align*}
  \ApproxWACosts(I,\beta)=  \frac{1}{\beta}\ApproxWAC{I}+\beta N.
\end{align*}
\end{defn}
In Section \ref{subsec: solution dependent} we elaborate on the solution and convergence of the minimization problem.

\subsection{Extreme value limit}\label{subsec: second order convergence}
In this section, we focus on the maximum of $N$ dependent random variables. In Theorem \ref{thm: second order convergence process} we prove that a scaled version of $\max_{i\leq N}Q_i(\beta)$ converges in distribution to a normally distributed random variable, as $N$ goes to infinity.

\begin{thm}\label{thm: second order convergence process}
Let $(W_i,1\leq i\leq N)$ be independent Brownian motions with mean 0 and variance $\sigma^2$, and $W_A$ be a Brownian motion with mean 0 and variance $\sigma_A^2$. Then
\begin{align}\label{eq: limit approximate problem dependent}
\frac{\max_{i\leq N}\sup_{s>0}\left(W_i(s)+W_A(s)-\beta s\right)-\frac{\sigma^2}{2\beta}\log N}{\sqrt{\log N}}\LimitD\frac{\sigma\sigma_A}{\sqrt{2}\beta}X,
\end{align}
with $X\sim\mathcal{N}(0,1)$. In other words, for all $x\in \mathbb{R}$
\begin{align*}
\probability*{\frac{\max_{i\leq N}\sup_{s>0}\left(W_i(s)+W_A(s)-\beta s\right)-\frac{\sigma^2}{2\beta}\log N}{\sqrt{\log N}}> x}\LimitN 1-\Phi\left(\frac{x\sqrt{2}\beta}{\sigma\sigma_A}\right),
\end{align*}
with $\Phi$ the cumulative distribution function of a standard normal random variable.
\end{thm}

A heuristic explanation of the result in Theorem \ref{thm: second order convergence process} is as follows: though $(Q_i,i\leq N)$ are dependent random variables, since we are adding the same Brownian motion $W_A$, $\max_{i\leq N} W_i(s)$ will dominate more and more over $W_A$ as $N$ becomes larger. Consequently, $W_A$ does not affect the time at which the supremum of $\max_{i\leq N} W_i(s)+W_A(s)-\beta s$ is attained. 
Hence, for $N$ large $\max_{i\leq N} Q_i(\beta)\approx \max_{i\leq N} \sup_{s>0}(W_i(s)-\beta s)+W_A(\tau)$, with $\tau$ the hitting time of the supremum of $\max_{i\leq N} (W_i(s)-\beta s)$.
Based on theory on conditional expectations of L\'evy processes we know that the conditional expectation of the hitting time $\tau(x)$ to reach a point $x$ is linear with $x$, to be precise, for $N=1$, it is known that $\expect*{\tau(x)\mid \tau(x)<\infty}=x/\beta$. Combining this with the fact that $\max_{i\leq N}\sup_{s>0}(W_i(s)-\beta s)\sim \frac{\sigma^2}{2\beta }\log N$, we expect that the supremum of $\max_{i\leq N} (W_i(s)-\beta s)$ is reached at $\tau\approx\frac{1}{\beta}\cdot\frac{\sigma^2}{2\beta}\log N=\frac{\sigma^2}{2\beta^2}\log N$.
Therefore, $W_A(\tau)\overset{d}{\approx} \frac{\sigma\sigma_A}{\sqrt{2}\beta}\sqrt{\log N}X$, with $X$ standard normally distributed, which results in Equation \eqref{eq: limit approximate problem dependent}. 

The proof of Theorem \ref{thm: second order convergence process} consists of four parts, which are stated in Lemmas \ref{lem: lower bound second order convergence process}, \ref{lem: upper bound second order convergence process <d-e}, \ref{lem: upper bound second order convergence process >d+e} and \ref{lem: upper bound second order convergence process d-e<s<d+e} for which the proofs are provided in Appendix \ref{app: proofs second order convergence}. For a process $X$ we have for all $t>0$ that
\begin{align*}
    \probability*{\sup_{s>0}X(s)>x}\geq \probability*{X(t)>x}.
\end{align*}
Furthermore, for every $0<t_1<t_2$,
\begin{align*}
    \probability*{\sup_{s>0}X(s)>x}
    \leq &\probability*{\sup_{0<s<t_1}X(s)>x}+\probability*{\sup_{t_1\leq s<t_2}X(s)>x}+\probability*{\sup_{s\geq t_2}X(s)>x}.
\end{align*}
We prove that these lower and upper bounds are tight for the process given in Theorem \ref{thm: second order convergence process} for appropriately chosen $t,t_1,t_2$.
More specifically, in Lemma \ref{lem: lower bound second order convergence process} we prove the asymptotic behavior at the critical time $d\log N$ where $d=\frac{\sigma^2}{2\beta^2}$, resulting in the tight lower bound. We show that times before and after this critical time have no influence in Lemmas \ref{lem: upper bound second order convergence process <d-e} and \ref{lem: upper bound second order convergence process >d+e}, respectively, leading up to Lemma \ref{lem: upper bound second order convergence process d-e<s<d+e} that shows the concentration around the critical time $d\log N$, proving a tight upper bound. 

\begin{lem}\label{lem: lower bound second order convergence process}
For $d=\frac{\sigma^2}{2\beta^2}$, 
\begin{align}
\frac{\max_{i\leq N}(W_i(d\log N)+W_A(d\log N))-\beta d\log N-\frac{\sigma^2}{2\beta}\log N}{\sqrt{\log N}}\LimitD\frac{\sigma\sigma_A}{\sqrt{2}\beta}X,
\end{align}
with $X\sim\mathcal{N}(0,1)$, as $N\to\infty$. 
\end{lem}

\begin{lem}\label{lem: upper bound second order convergence process <d-e}
For $d=\frac{\sigma^2}{2\beta^2}$ and $0<\epsilon<d$, and for all $x$, 
\begin{align}
\mathbb{P}\left(\frac{\max_{i\leq N}\sup_{0<s<(d-\epsilon)\log N}\left(W_i(s)+W_A(s)-\beta s\right)-\frac{\sigma^2}{2\beta}\log N}{\sqrt{\log N}}\geq x\right)\LimitN 0.
\end{align}
\end{lem}

\begin{lem}\label{lem: upper bound second order convergence process >d+e}
For $d=\frac{\sigma^2}{2\beta^2}$ and all $\epsilon>0$, and $x\in \mathbb{R}$, 
\begin{align}\label{eq: upper bound second order convergence process >d+e}
\mathbb{P}\left(\frac{\max_{i\leq N}\sup_{s\geq(d+\epsilon)\log N}\left(W_i(s)+W_A(s)-\beta s\right)-\frac{\sigma^2}{2\beta}\log N}{\sqrt{\log N}}\geq x\right)\LimitN 0.
\end{align}
\end{lem}

\begin{lem}\label{lem: upper bound second order convergence process d-e<s<d+e}
For $d=\frac{\sigma^2}{2\beta^2}$ and $\epsilon>0$ and for all $x$, 
\begin{align}
\limsup_{N\to\infty}&\mathbb{P}\left(\frac{\max_{i\leq N}\sup_{(d-\epsilon)\log N\leq s <(d+\epsilon)\log N}\left(W_i(s)+W_A(s)-\beta s\right)-\frac{\sigma^2}{2\beta}\log N}{\sqrt{\log N}}\geq x\right)\nonumber\\
\leq &\mathbb{P}\left(\sigma_A\sqrt{\frac{\sigma^2}{2\beta^2}-\epsilon}X_1+\sqrt{2\epsilon}\sigma_A|X_2|>x\right),
\end{align}
with $X_1,X_2\sim\mathcal{N}(0,1)$ and independent.
\end{lem}
In Appendix \ref{app: proofs second order convergence} we show how these lemmas can be used to prove Theorem \ref{thm: second order convergence process}. In Lemma \ref{lem: L1 convergence}, we prove that convergence holds even in $L_1$, when $X$ is chosen approprately.
\begin{lem}\label{lem: L1 convergence}
Define $X_N:=\frac{\sqrt{2}\beta}{\sigma\sigma_A}\frac{W_A\left(\frac{\sigma^2}{2\beta^2}\log N\right)}{\sqrt{\log N}}$. Then,
\begin{align*}
    \expect*{\left|\frac{\max_{i\leq N}\sup_{s>0}\left(W_i(s)+W_A(s)-\beta s\right)-\frac{\sigma^2}{2\beta}\log N}{\sqrt{\log N}}-\frac{\sigma\sigma_A}{\sqrt{2}\beta}X_N\right|}\LimitN 0.
\end{align*}
\end{lem}

The proof of Lemma \ref{lem: L1 convergence} is also given in Appendix \ref{app: proofs second order convergence}. In the next section, we apply Theorem \ref{thm: second order convergence process} and Lemma \ref{lem: L1 convergence} to solve and approximate the minimization problem. Specifically, Lemma \ref{lem: L1 convergence} gives us an order bound between the optimal base-stock level and the approximate base-stock level.

\subsection{Solution and Convergence of the Minimization Problem}\label{subsec: solution dependent}
We can use the convergence result proven in Theorem \ref{thm: second order convergence process} to prove asymptotics of the minimization of the function $F_N$. Since $ \frac{\sqrt{2}\beta}{\sigma\sigma_A}  \frac{\max_{i\leq N}Q_i(\beta)-\frac{\sigma^2}{2\beta}\log N}{\sqrt{\log N}}$ is a continuous random variable, we know that its quantile function converges to the quantile function of a standard normal random variable; cf. \citet[p.~305, Lem.~21.2]{vaart_1998}. So we can use this to derive asymptotics of the minimization problem of $F_N$. 

Using $P_N^{A}(z)$ as described in Definition \ref{def: P_N^A}, we can solve the minimization problem, which yields the optimal base-stock level and net capacity given in Lemma \ref{lem: solution minimization problem}. The proofs concerning the solution and subsequent convergence results are provided in Appendix \ref{app: proofs dependent model}.
\begin{defn}\label{def: P_N^A}
We define 
\begin{align*}
    P_N^A(z)=\probability*{\frac{\sqrt{2}}{\sigma\sigma_A}\frac{\max_{i\leq N}Q_i-\frac{\sigma^2}{2}\log N}{\sqrt{\log N}}\leq z}.
\end{align*}
\end{defn}
\begin{lem}\label{lem: solution minimization problem}
Let $(b^{(N)})_{N\geq 1}, (h^{(N)})_{N\geq 1}$ be sequences such that $h^{(N)}>0$ and $b^{(N)}>0$ for all $N$, and $\gamma_N=Nh^{(N)}/(Nh^{(N)}+b^{(N)})$. Let $\left(\OptimalWABetaN,\OptimalWAIN\right)$ minimize $F_N(I,\beta)$. Then 
\begin{align}\label{eq: OptimalWAIN}
    \OptimalWAIN=\frac{\sigma^2}{2}\log N+\frac{\sigma\sigma_A}{\sqrt{2}}{P_N^A}^{-1}\left(1-\gamma_N\right)\sqrt{\log N}.
\end{align}
\end{lem}
When we are in the balanced regime, we can approximate the minimization problem given in Definition \ref{def: approxcosts dependent}, using the convergence result in Theorem \ref{thm: second order convergence process}, and prove how far the approximate solution is from the optimal solution. This is done in Proposition \ref{prop: asymptotic approximation dependent} and Theorem \ref{thm: order bounds dep}. In Lemma \ref{lem: bal reg dep} we show how the optimal costs scale with $N$ when we are in the balanced regime. The proofs are given in Appendix \ref{app: proofs dependent model}.

\begin{prop}\label{prop: asymptotic approximation dependent}
For $(b^{(N)})_{N\geq 1}, (h^{(N)})_{N\geq 1}$ and $\gamma_N=Nh^{(N)}/(Nh^{(N)}+b^{(N)})$,
\begin{align}\label{eq: asymptotic I approximation dependent}
    \ApproxWAIN=\frac{\sigma^2}{2}\log N+\frac{\sigma  \sigma _A }{\sqrt{2} }\sqrt{\log N} \Phi^{-1}\left(1-\gamma_N\right),
\end{align}
and
\begin{align}\label{eq: C tilde dep}
    \ApproxWAC{\ApproxWAIN}=Nh^{(N)} \left(\frac{\sigma ^2}{2} \log N-\frac{\sigma ^2+\sigma_A ^2}{2}\right)+(Nh^{(N)}+b^{(N)})\frac{ \sigma  \sigma _A \sqrt{\log N} e^{-\frac{1}{2}\Phi^{-1}\left(1-\gamma_N\right)^2}}{2\sqrt{\pi }}.
\end{align}

\end{prop}
\begin{thm}[Order bound]\label{thm: order bounds dep}
Assume $\gamma_N=Nh^{(N)}/(Nh^{(N)}+b^{(N)})$, with $\gamma_N=\gamma\in(0,1)$. Then
\begin{align*}
    \left|\frac{F_N(\OptimalWAIN,\OptimalWABetaN)   }{F_N(\ApproxWAIN,\ApproxWABetaN)}-1\right|=o\left(\frac{1}{\sqrt{\log N}}\right).
\end{align*}
\end{thm}

\begin{lem}[Balanced regime]\label{lem: bal reg dep}
Assume $\gamma_N=Nh^{(N)}/(Nh^{(N)}+b^{(N)})$, with $\gamma_N=\gamma\in(0,1)$. Then
\begin{align}\label{eq: bal reg dep inventory}
    \OptimalWAIN= \frac{\sigma^2}{2}\log N+\frac{\sigma  \sigma _A }{\sqrt{2} }\sqrt{\log N} \Phi^{-1}\left(1-\gamma\right)+o(\sqrt{\log N}),
\end{align}
and
\begin{align}\label{eq: bal reg dep costs}
  F_N(\OptimalWAIN,\OptimalWABetaN)= 2\sqrt{N}\sqrt{\ApproxWAC{\ApproxWAIN}}+o(N\sqrt{h^{(N)}}).
  \end{align}
\end{lem}

The result in Lemma \ref{lem: bal reg dep} only holds for the balanced regime, so a natural question is what we can say about the efficiency and the quality driven regime. As is shown in Lemma \ref{lem: first order bal reg}, in the efficiency driven regime, the first order approximation $\bar I_N=\frac{\sigma^2}{2}\log N$ gives that the ratio of the approximate costs and the optimal costs converge to 1. Thus we expect the approximation given in \eqref{eq: asymptotic I approximation dependent} will also satisfy this convergence result. In order to determine whether this approximation also satisfies the order bound given in Theorem \ref{thm: order bounds dep}, a further analysis is needed. The analysis we provide for the balanced regime heavily relies on \citet[p.~305, Lem.~21.2]{vaart_1998}, which says that if $Y_N\LimitD Y$, then for $\gamma\in(0,1)$, $P^{-1}_{Y_N}(\gamma)\LimitN P^{-1}_{Y}(\gamma)$. This gives us the convergence result \eqref{eq: bal reg dep inventory} of the inventory in the balanced regime. In order to be able to prove a similar result for the efficiency driven regime, we need an improvement of \citet[p.~305, Lem.~21.2]{vaart_1998} which also holds when $\gamma_N\LimitN 1$.

However, for the quality driven regime, this convergence result does not hold, because we see in Lemma \ref{lem: magnitude F} that $\OptimalWAIN\approx \frac{\sigma^2}{2}\log(N/\gamma_N)$. In order to find a sharp order bound such as given in Theorem \ref{thm: order bounds dep} we should resort to the analysis of tail asymptotics, which is beyond the scope of this study.

\subsection{Numerical Experiments}\label{subsec: numerical results dependent}
In Section \ref{subsec: solution dependent}, we provided expressions to calculate the asymptotically optimal net capacity and base-stock level. The question remains how large the number of components has to be for these approximations to be of use. Therefore, we now examine the expected costs under both the optimal net capacity and base-stock level and under these asymptotic approximations. Since it is not straightforward to calculate $\mathbb{E}\left[ \left(\max_{i\leq N}Q_i-I\right)^+ \right]$ for dependent $Q_i$, to evaluate the cost function given in Definition \ref{def: optimalcosts independent} we resort to simulation. First, we explain the details of our simulation experiment, after which we discuss the numerical results.

In our simulation, we aim to determine the maximum delay over all components, so $\max_{i\leq N} Q_i$. For this, we use the algorithm proposed by \citet[\S4.5]{asmussen1995}, who describe an exact algorithm for simulating a reflected Brownian motion at the grid points. At every grid point, we draw normal random variables with the required drift and variance for the supply and demand processes and update the maximum. We use a step size of 0.001 for the grid points. Since we cannot simulate over an infinite horizon, we have to determine when to terminate the simulation. The maximum value is expected to be attained at a time which is smaller than $\hat{t}=\frac{\sigma^2+\sigma_A^2}{2}\sum_{j=1}^N \frac{1}{j}$. To simulate well beyond this point, we run the simulation until $t=2\hat{t}$. 

Using the above method to simulate $\max_{i\leq N} Q_i$, we can estimate ${P_N^A}^{-1}(1-\gamma_N)$ with $P_N^A(z)$ as described in Definition \ref{def: P_N^A}. 
To obtain a median-unbiased estimate of the quantile, we use the approach suggested by \citet[p.~982-983]{zielinski2009}. For this, we sample $\max_{i\leq N} Q_i$ 100 times and randomly choose between the observations $(1-\gamma_N)\cdot 100$ and $(1-\gamma_N)\cdot 100+1$, with weights depending on the value of the fractile. Our estimate is equal to the median over 100 iterations. 
Once we have our estimate of ${P_N^A}^{-1}(1-\gamma_N)$, we determine the value of the optimal base-stock level as given in Equation \eqref{eq: OptimalWAIN}.
Using the optimal base-stock level we determine the optimal net capacity given in Lemma \ref{lem: optimal beta general}. Since this also requires the expectation of $\left(\max_{i\leq N}Q_i-I\right)^+$, we determine this value by taking the average based on 10,000 simulations.

Next, we compare the costs under our asymptotic approximations of the net capacity and base-stock level (provided in Proposition \ref{prop: asymptotic approximation dependent}) to the costs under the optimal net capacity and base-stock level obtained from the simulation. We again sample $\left(\max_{i\leq N}Q_i-I\right)^+$ based on 10,000 new simulations and determine the costs of the different policies using cost function $F_N(I,\beta)$. 

The procedure described above is applicable for $N$ in the order of hundreds, however, it is close to impossible to provide a fast simulation for $N$ in the order of thousands. Hence, to give a useful approximation of the optimal capacity and base-stock level in these cases, we need to use the limit we derived in Theorem \ref{thm: second order convergence process}.

In order to assess the performance of the approximations and its sensitivity to various model parameters, we perform a full factorial experiment. In our experiment, we vary the number of components, demand variability and backorder costs. The setup of the experiment is given in Table \ref{tab: experiment setup}. We set $h^{(N)}=1$ and $\sigma=1$ in all experiments. In total we have 24 instances. The results are given in Tables \ref{tab: costs asymptotic} and \ref{tab: costs asymptotic b=3N} for $b^{(N)}=N$ and $b^{(N)}=3N$, respectively.

\begin{table}[H]
    \centering
    \caption{Parameter settings for experiments}
    \begin{tabular}{ll}
        \hline
        Parameter & Values  \\
        \hline
        $N$ & 10, 50, 100\\
        $\sigma_A$ & 0.1, 0.5, 0.75, 1\\
        $b^{(N)}$ & $N$, $3N$\\
        \hline
    \end{tabular}
    \label{tab: experiment setup}
\end{table}

\begin{table}[h]
\caption{Comparison of costs approximate solution for $h^{(N)}=1$, $b^{(N)}=N$}
\label{tab: costs asymptotic}
\begin{tabular}{c|c|c|c|c|c|c|c|c}
$N$ & $\sigma_A$& $ I^A_N $ & $\beta^A_N$ & $F_N(I^A_N,\beta^A_N)$  &  $\ApproxWAIN$ & $\ApproxWABetaN$ & $F_N(\ApproxWAIN,\ApproxWABetaN)$&  $\left(1-\frac{F_N(I^A_N,\beta^A_N)}{F_N(\ApproxWAIN,\ApproxWABetaN)}\right)\sqrt{\log N}$\\
\hline
10 & 0.1 & 1.327 & 1.1583 & 23.1894 & 1.151 & 0.855514 & 24.5143 &  0.0820\\
50 & 0.1 & 2.122 & 1.47611 & 147.534 & 1.956 & 1.25004 & 150.337 &  0.0369\\
100 & 0.1 & 2.455 & 1.58865 & 318.588 & 2.303 & 1.38516 & 322.994 &  0.0293\\ \hline
10 & 0.5 & 1.486 & 1.25448 & 25.333 & 1.151 & 0.976909 & 26.9363 &  0.0903\\
50 & 0.5 & 2.338 & 1.59412 & 159.934 & 1.956 & 1.3744 & 164.689 &  0.0571\\
100 & 0.5 & 2.715 & 1.71664 & 343.937 & 2.303 & 1.51094 & 352.91 &  0.0546 \\ \hline
10 & 0.75 & 1.714 & 1.36908 & 27.191 & 1.151 & 1.00605 & 29.7614 &  0.1311\\
50 & 0.75 & 2.638 & 1.70591 & 171.443 & 1.956 & 1.41834 & 180.556 &  0.0998\\
100 & 0.75 & 2.980 & 1.83438 & 367.348 & 2.303 & 1.55865 & 383.319 &  0.0894\\ \hline
10 & 1 & 1.990 & 1.47358 & 29.8393 & 1.151 & 1.0037 & 34.6552 &  0.2109\\
50 & 1 & 3.006 & 1.84276 & 185.25 & 1.956 & 1.43941 & 201.314 &  0.1578\\
100 & 1 & 3.394 & 1.97602 & 393.668 & 2.303 & 1.58534 & 421.505 &  0.1417\\
\end{tabular}
\end{table}

\begin{table}[h]
\caption{Comparison of costs approximate solution for $h^{(N)}=1$, $b^{(N)}=3N$}
\label{tab: costs asymptotic b=3N}
\begin{tabular}{c|c|c|c|c|c|c|c|c}
$N$ & $\sigma_A$& $ I^A_N $ & $\beta^A_N$ & $F_N(I^A_N,\beta^A_N)$  &  $\ApproxWAIN$ & $\ApproxWABetaN$ & $F_N(\ApproxWAIN,\ApproxWABetaN)$&  $\left(1-\frac{F_N(I^A_N,\beta^A_N)}{F_N(\ApproxWAIN,\ApproxWABetaN)}\right)\sqrt{\log N}$\\
\hline
10 & 0.1 & 1.726 & 1.31058 & 25.9539 & 1.224 & 0.884692 & 31.2239 &  0.2561\\
50 & 0.1 & 2.533 & 1.5931 & 159.026 & 2.050 & 1.27624 & 173.141 &  0.1612\\
100 & 0.1 & 2.883 & 1.69656 & 341.44 & 2.405 & 1.41084 & 367.575 &  0.1526\\ \hline
10 & 0.5 & 2.067 & 1.43331 & 28.3311 & 1.513 & 1.0992 & 31.2606 &  0.1422\\
50 & 0.5 & 2.987 & 1.74381 & 173.875 & 2.428 & 1.48993 & 183.166 &  0.1003\\
100 & 0.5 & 3.370 & 1.86469 & 371.779 & 2.814 & 1.62542 & 387.809 &  0.0887\\ \hline
10 & 0.75 & 2.449 & 1.57036 & 31.4004 & 1.694 & 1.18023 & 35.5139 &  0.1758\\
50 & 0.75 & 3.418 & 1.89842 & 190.571 & 2.664 & 1.58369 & 205.174 &  0.1408\\
100 & 0.75 & 3.899 & 2.01955 & 404.306 & 3.070 & 1.72277 & 429.58 &  0.1263\\ \hline
10 & 1 & 2.913 & 1.72878 & 34.6096 & 1.875 & 1.23092 & 40.7704 &  0.2293\\
50 & 1 & 4.158 & 2.06968 & 207.553 & 2.899 & 1.65341 & 230.281 &  0.1952\\
100 & 1 & 4.567 & 2.20696 & 439.681 & 3.326 & 1.79761 & 479.663 &  0.1789\\
\end{tabular}
\end{table}

There are several important observations to be made from Table \ref{tab: costs asymptotic}. First of all, we can observe that for $N=10$ the difference in costs between the simulated optimal solution and the asymptotic solution is around 10\% for most cases, the case $N=10$ and $\sigma_A=1$ is an outlier, where the difference is around 15\%. As $N$ increases to 50, the difference decreases. Furthermore, the difference becomes larger when $\sigma$ increases.  In the last column, we verify the convergence result from Theorem \ref{thm: order bounds dep}. We observe that the difference decreases as $N$ increases, and that increasing $\sigma_A$ causes the difference to increase. 

When we consider the results for $b^{(N)}=3N$ given in Table \ref{tab: costs asymptotic b=3N}, we observe that the difference between the asymptotic and optimal costs is considerably higher than for $b^{(N)}=N$. Especially for $N=10$, the difference is around 15\% of the optimum, except for $N=10$ and $\sigma_A=0.1$, where the difference is around 20\%. However, for a larger number of components, the difference is around 10\% of the optimum. Interestingly, for the case $\sigma_A=1$, the difference between $b^{(N)}=N$ and $b^{(N)}=3N$ is relatively small.

Overall, in most of our experiments the difference between the costs under the optimal base-stock level and net capacity and the costs under the approximations are around 10\%. Furthermore, we can conclude that for small variations in demand and low backorder costs, the asymptotic approach performs well in terms of costs already for a reasonable number of components. Also, the performance improves by increasing $N$. Finally, the performance of the approximations highly depends on the backorder costs relative to the holding costs.

\section{Mixed-behavior approximations}\label{sec: master formula}
The numerical results in Section \ref{subsec: numerical results dependent} show that the approximations are in most of the cases around 10-15\% off the optimal value. In this section, we show how we can further improve the approximations.

Under deterministic demand and stochastic demand, the approximate problems are given in Definition \ref{def: approxcosts independent} and Definition \ref{def: approxcosts dependent}. If $\sigma_A$ is small, then we know that on the one hand,
\begin{align*}
     \max_{i\leq N} Q_i\overset{d}{\approx} \frac{\sigma^2}{2}G+\frac{\sigma^2}{2}\log N,
\end{align*}
because $Q_i$ and $Q_j$ are only slightly correlated. But on the other hand,
\begin{align*}
     \max_{i\leq N} Q_i\overset{d}{\approx} \frac{\sigma\sigma_A}{\sqrt{2}}\sqrt{\log N}X+\frac{\sigma^2}{2}\log N\approx \frac{\sigma^2}{2}\log N. 
\end{align*}
Since the Gumbel term is missing here, this could be the reason that this approximation is not working well for small $N$. Thus, it could be beneficial to look at the combination of these two approximations. Then, we have
\begin{align}\label{eq: max Q_i master}
\max_{i\leq N}Q_i\overset{d}{\approx} \frac{\sigma^2}{2}\log N+\frac{\sigma\sigma_A}{\sqrt{2}}\sqrt{\log N}X+\frac{\sigma^2}{2}G.  
\end{align}
When we replace $\max_{i\leq N}Q_i$ with Equation \eqref{eq: max Q_i master} in the minimization problem, we get
\begin{align*}
    \min_{I,\beta}\left(\frac{1}{\beta}\expect*{Nh^{(N)}(I-Q_i)+(Nh^{(N)}+b^{(N)})\left(\frac{\sigma^2}{2}\log N+\frac{\sigma\sigma_A}{\sqrt{2}}\sqrt{\log N}X+\frac{\sigma^2}{2}G-I\right)^+}+\beta N\right).
\end{align*}
The optimal $I^{M}_N$ satisfies $\probability*{\frac{\sigma^2}{2}\log N+\frac{\sigma\sigma_A}{\sqrt{2}}\sqrt{\log N}X+\frac{\sigma^2}{2}G<I^{M}_N}=1-\gamma_N$.
Thus,
\begin{align}\label{eq: computation I master}
    \int_{-\infty}^{\infty} \exp\bigg(-\exp\bigg(-\frac{2}{\sigma^2}\bigg(I^{M}_N-\frac{\sigma^2}{2}\log N-\frac{\sigma\sigma_A}{\sqrt{2}}\sqrt{\log N}x\bigg)\bigg)\bigg)\phi(x)dx=1-\gamma_N.
\end{align}
Now, $I^{M}_N$ can be computed through standard numerical methods such as the bisection method. Furthermore, the optimal net capacity $\beta^{M}_N$ satisfies
\begin{align}\label{eq: symb expr beta mixed behavior}
   \beta^{M}_N= \frac{\sqrt{\expect*{Nh^{(N)}(I^{M}_N-Q_i)+(Nh^{(N)}+b^{(N)})\left(\frac{\sigma^2}{2}\log N+\frac{\sigma\sigma_A}{\sqrt{2}}\sqrt{\log N}X+\frac{\sigma^2}{2}G-I^{M}_N\right)^+}}}{\sqrt{N}}.
\end{align}
The relevant expectations in this symbolic expression can be computed numerically; see Appendix \ref{app: mixed behavior} for details.

\subsection{Numerical results mixed-behavior approximations}
Using the same simulation procedure as described in Section \ref{subsec: numerical results dependent}, we evaluate the performance of these adjusted approximations. The results for the cases of $h^{(N)}=1,\:b^{(N)}=N$ and $h^{(N)}=1,\:b^{(N)}=3N$ are given in Tables \ref{tab: costs master} and \ref{tab: costs master b=3N}, respectively. 

From the simulation results we can conclude that these adjusted approximations result in costs that are much closer to the optimal costs, already for small $N$. When comparing the last two columns, where the last column repeats the results from Section \ref{subsec: numerical results dependent}, we observe that the mixed-behavior approximations show better convergence, also when $\sigma_A$ is larger.  
Furthermore, where we saw in Section \ref{subsec: numerical results dependent} that the cost difference increased considerably with the change in $b^{(N)}$, we now do see a slight increase, but the difference is still small for a larger value of $b^{(N)}$. Therefore, we can conclude that these mixed-behavior approximations perform well especially when demand variations are no more than 75\% of the variations in component production, even with a small number of components. 

\begin{table}[H]
\caption{Comparison of costs master solution for $h^{(N)}=1$, $b^{(N)}=N$}
\label{tab: costs master}
\begin{tabular}{c|c|c|c|c|c|c}

$N$& $\sigma_A$& $I^{M}_N$ & $\beta^{M}_N$ & $F_N(I^{M}_N,\beta^{M}_N)$& $\left(1-\frac{F_N(I^A_N,\beta^A_N)}{F_N(I^{M}_N,\beta^{M}_N)}\right)\sqrt{\log N}$ & $\left(1-\frac{F_N(I^A_N,\beta^A_N)}{F_N(\ApproxWAIN,\ApproxWABetaN)}\right)\sqrt{\log N}$\\
\hline
10	&	0.1	&	1.33785	&	1.1945	&	23.2022	&	0.000837	&	0.082011	\\	
50	&	0.1	&	2.14487	&	1.49567	&	147.567	&	0.000442	&	0.036877	\\	
100	&	0.1	&	2.49244	&	1.60808	&	318.638	&	0.000337	&	0.029273	\\	\hline
10	&	0.5	&	1.38072	&	1.21129	&	25.4342	&	0.006038	&	0.090320	\\	
50	&	0.5	&	2.19829	&	1.53814	&	160.497	&	0.006938	&	0.057107	\\	
100	&	0.5	&	2.54871	&	1.65808	&	345.247	&	0.008143	&	0.054563	\\	\hline
10	&	0.75	&	1.40013	&	1.2128	&	27.6956	&	0.027647	&	0.131055	\\	
50	&	0.75	&	2.216	&	1.56166	&	174.269	&	0.032074	&	0.099827	\\	
100	&	0.75	&	2.5656	&	1.68745	&	372.643	&	0.030493	&	0.089412	\\	\hline
10	&	1	&	1.41255	&	1.19665	&	31.5428	&	0.081950	&	0.210871	\\	
50	&	1	&	2.22627	&	1.57136	&	192.722	&	0.076684	&	0.157827	\\	
100	&	1	&	2.57434	&	1.70384	&	407.343	&	0.072043	&	0.141724	\\	
\end{tabular}
\end{table}


\begin{table}[H]
\caption{Comparison of costs master solution for $h^{(N)}=1$, $b^{(N)}=3N$}
\label{tab: costs master b=3N}
\begin{tabular}{c|c|c|c|c|c|c}
$N$& $\sigma_A$& $I^{M}_N$ & $\beta^{M}_N$ & $F_N(I^{M}_N,\beta^{M}_N)$& $\left(1-\frac{F_N(I^A_N,\beta^A_N)}{F_N(I^{M}_N,\beta^{M}_N)}\right)\sqrt{\log N}$ & $\left(1-\frac{F_N(I^A_N,\beta^A_N)}{F_N(\ApproxWAIN,\ApproxWABetaN)}\right)\sqrt{\log N}$\\
\hline
10	&	0.1	&	1.78238	&	1.34746	&	25.9965	&	0.002487	&	0.256113	\\	
50	&	0.1	&	2.59271	&	1.62088	&	159.162	&	0.001690	&	0.161243	\\	
100	&	0.1	&	2.94168	&	1.72533	&	341.49	&	0.000314	&	0.152581	\\	\hline
10	&	0.5	&	1.94345	&	1.38309	&	28.3671	&	0.001926	&	0.142201	\\	
50	&	0.5	&	2.83775	&	1.68955	&	174.284	&	0.004642	&	0.100327	\\	
100	&	0.5	&	3.21861	&	1.8044	&	372.617	&	0.004826	&	0.088703	\\	\hline
10	&	0.75	&	2.09429	&	1.41142	&	32.0055	&	0.028689	&	0.175760	\\	
50	&	0.75	&	3.04648	&	1.74512	&	193.854	&	0.033496	&	0.140773	\\	
100	&	0.75	&	3.44819	&	1.86761	&	410.624	&	0.033019	&	0.126256	\\	\hline
10	&	1	&	2.25658	&	1.43095	&	36.5165	&	0.079240	&	0.229298	\\	
50	&	1	&	3.26538	&	1.79271	&	216.91	&	0.085321	&	0.195211	\\	
100	&	1	&	3.68765	&	1.92281	&	456.859	&	0.080689	&	0.178876	\\	
\end{tabular}
\end{table}

\section{Analyzing asymmetric systems}\label{sec: asymmetry}
This paper derived several new, analytical results for joint capacity and inventory optimization for large-scale, symmetric assembly systems. In this section we provide an informal discussion of the application of such results in asymmetric settings. 

For ease of exposition, consider a case where different components have different holding costs. For other parameters, our assumptions remain in place. In practical settings, component prices might range from a few thousand euros to hundreds of thousands of euros. Companies seeking to apply advanced methods for optimizing capacity and inventory investments would focus on the most expensive components: For inexpensive components some coarse heuristics would suffice. 

Suppose the company seeks to derive separate inventory buffer and capacity rules for two groups of components: Expensive and very expensive components. This yields $k=2$ groups of components. We seek to apply our results on extremes as the total number of components $N$ in these two groups grows large; we keep $k$ and the ratio of components in the two groups fixed. Also, since we seek to derive rules at the group-level, it makes sense to assume symmetry within groups, i.e.\ by averaging cost parameters within the groups. For example, consider the following: $N/2$ servers have a holding cost $h_1^{(N)}$ and $N/2$ servers have a holding cost $h_2^{(N)}$. Then we need to minimize 

\begin{multline}\label{eq: cost function asymm}
\frac{N}{2}\left(h_1^{(N)}\frac{1}{\beta_1}\left(I_1-\frac{\sigma^2}{2}\right)+\beta_1\right)
+\frac{N}{2}\left(h_2^{(N)}\frac{1}{\beta_2}\left(I_2-\frac{\sigma^2}{2}\right)+\beta_2\right)\\
+\left(\frac{N}{2}h_1^{(N)}+\frac{N}{2}h_2^{(N)}+b^{(N)}\right)\expect*{\max\left(\frac{1}{\beta_1}\max_{i\leq N/2}(Q_i(1)-I_1),\frac{1}{\beta_2}\max_{N/2+1\leq i\leq N}(Q_i(1)-I_2)\right)^+}.
\end{multline}
over $(I_1,I_2,\beta_1,\beta_2)$.
Obviously,
\begin{multline*}
    \expect*{\max\left(\frac{1}{\beta_1}\max_{i\leq N/2}(Q_i(1)-I_1),\frac{1}{\beta_2}\max_{N/2+1\leq i\leq N}(Q_i(1)-I_2)\right)^+}\\
\leq \expect*{\frac{1}{\beta_1}\max_{i\leq N/2}(Q_i(1)-I_1)^+}+\expect*{\frac{1}{\beta_2}\max_{i\leq N/2}(Q_i(1)-I_2)^+}.
\end{multline*}
The cost function in Equation \eqref{eq: cost function asymm} can therefore be bounded from above by
\begin{multline*}
   \frac{N}{2}\left(h_1^{(N)}\frac{1}{\beta_1}\left(I_1-\frac{\sigma^2}{2}\right)+\beta_1\right)
+\frac{N}{2}\left(h_2^{(N)}\frac{1}{\beta_2}\left(I_2-\frac{\sigma^2}{2}\right)+\beta_2\right)\\
+\left(\frac{N}{2}h_1^{(N)}+\frac{N}{2}h_2^{(N)}+b^{(N)}\right)\left(\expect*{\frac{1}{\beta_1}\max_{i\leq N/2}(Q_i(1)-I_1)^+}+\expect*{\frac{1}{\beta_2}\max_{i\leq N/2}(Q_i(1)-I_2)^+}\right).
\end{multline*}
Our analytical results enable us to minimize this upper bound; for instance choosing $\tilde{h}_{1,2}^{(N)}=h_{1,2}^{(N)}$ and $\tilde{b}_{1,2}^{(N)}=\frac{N}{2}h_{2,1}^{(N)}+b^{(N)}$ yields
 \begin{align*}
   &\frac{N}{2}\left(h_1^{(N)}\frac{1}{\beta_1}\left(I_1-\frac{\sigma^2}{2}\right)+\beta_1\right)
+\frac{N}{2}\left(h_2^{(N)}\frac{1}{\beta_2}\left(I_2-\frac{\sigma^2}{2}\right)+\beta_2\right)\\
&\qquad+\left(\frac{N}{2}h_1^{(N)}+\frac{N}{2}h_2^{(N)}+b^{(N)}\right)\left(\expect*{\frac{1}{\beta_1}\max_{i\leq N/2}(Q_i(1)-I_1)^+}+\expect*{\frac{1}{\beta_2}\max_{i\leq N/2}(Q_i(1)-I_2)^+}\right)\\
&\quad=\frac{N}{2}\left(\tilde{h}_1^{(N)}\frac{1}{\beta_1}\left(I_1-\frac{\sigma^2}{2}\right)+\beta_1\right)+\left(\frac{N}{2}\tilde{h}_1^{(N)}+\tilde{b}_1^{(N)}\right)\expect*{\frac{1}{\beta_1}\max_{i\leq N/2}(Q_i(1)-I_1)^+}\\
&\qquad+\frac{N}{2}\left(\tilde{h}_2^{(N)}\frac{1}{\beta_2}\left(I_2-\frac{\sigma^2}{2}\right)+\beta_2\right)+\left(\frac{N}{2}\tilde{h}_2^{(N)}+\tilde{b}_2^{(N)}\right)\expect*{\frac{1}{\beta_2}\max_{i\leq N/2}(Q_i(1)-I_2)^+}.
\end{align*}
This is the sum of two functions that can be minimized using the exact solutions that we derived. In Table \ref{tab: numerical results asymmetric} we compare numerically the actual costs under the capacity and base-stock level that are obtained by minimizing this upper bound, with the costs under the optimal capacity and base-stock level. In this table, the ratio indicates how many servers have a holding cost $h_1^{(N)}$, and how many servers have a holding cost $h_2^{(N)}$; the 1:1 ratio corresponds to the above example while the 1:3 ratio can be treated similarly. The table demonstrates that our asymptotic results may be useful when optimizing asymmetric systems as well as symmetric systems. 
\begin{table}[H]
\centering
\caption{Comparison of optimal costs and costs under upper bound heuristic, $\sigma=1$, $\sigma_A=0$.}
\label{tab: numerical results asymmetric}
\begin{tabular}{c|c|c|c|c|c|c|c}
$N$ & $h_1^{(N)}$ & $h_2^{(N)}$& Ratio &$b^{(N)}$& \text{Optimal}  & \text{Heuristic} & Diff. \\
\hline
10	&	1	&	10	& 1:1	&  10	&	42.3 $\pm$ 0.1	&	 42.9	$\pm$ 0.1 & 0.14 \% \\
100	&	1	&	10	& 1:1	&  100	&	615.6	$\pm$ 1.2 &	617.4 $\pm$ 1.0 &	0.3 \%\\
1000	&	1	&	10	& 1:1	&  1000	&	7597.9 $\pm$ 8.2 &	 7643.0 $\pm$ 7.8&  0.6 \%
\\
10	&	10	&	100	& 1:1	&  1	&	126.0 $\pm$ 0.4	&	127.0	$\pm$ 0.4&0.7 \%\\
100	&	100	&	1000	& 1:1	&  1 &	5967$\pm$10.9	&	6002	$\pm$	9.6&0.6 \%\\
1000	&	1000	&	10000	& 1:1	&  1	&	 236063	$\pm$	256&236402$\pm$ 233	& 0.1 \%
\\
10	&	1	&	10	& 1:3	&  10	&	53.1$\pm$0.2	&	53.2 $\pm$ 0.2&0.2 \%\\
100	&	1	&	10	& 1:3	&  100	&	770.5$\pm$1.3	&	772.9$\pm$	1.2 &0.3 \%\\
1000	&	1	&	10	& 1:3	&  1000	& 9551.1$\pm$10.7	 &  9581.6$\pm$ 9.5&0.3 \%	
\end{tabular}
\end{table}

\section{Conclusions}\label{sec: conclusions}
In this study, we defined a large scale assembly system in which $N$ components are assembled into a final product. We studied an assembly system with linear demand and production, subjected to some random noise. Thus, we imposed the natural assumption that this noise is normally distributed. Hence, delays per component are written as an all-time supremum of a Brownian motion minus a drift term. We aimed to minimize the total costs in the system with respect to the inventory and net capacity per component. The costs in the system consist of inventory holding costs for each component and penalty costs for delays in assembly of the final product, which is equal to the delay of the slowest produced component.  Before attempting to solve the minimization problem, we simplified the minimization problem, using the self-similarity property of a Brownian motion, into two separate minimization problems. We distinguished two cases, first of all we covered the case of deterministic demand, resulting in all delays being independent. Secondly, we investigated the case that demand is stochastic and consequently delays of the components are dependent.

For the deterministic demand scenario, we proved order bounds for three different regimes: balanced, quality driven and efficiency driven. Additionally, we verified numerically that already for a limited number of components, our approximations result in costs that are very close to the costs corresponding to the optimal solution.
For the stochastic demand scenario, we developed a limit theorem that we use to obtain approximate solutions. We showed numerically that even though theoretically these approximations perform well, for practical situations there is still room for improvement. However, this limit theorem is still necessary for systems with $N$ of the order of thousands, because it is close to impossible to simulate these systems fast. Therefore, we provided additional approximations for a mixed-behavior regime, where we use a combination of the approximations for the deterministic and stochastic demand scenarios. We demonstrated numerically that these approximations perform very well already for a practical number of components.

Future work could extend the model to a decentralized minimization problem, where the components are not produced in-house by the manufacturer but are sourced at outside suppliers that have their own objectives, which results in an asymptotic analysis of a game theoretical equilibrium, cf.\ \cite{nair2016provisioning, gopalakrishnan2016routing} and \cite{kumar2010exploiting}. Additionally, we expect that we can extend the result in Theorem \ref{thm: second order convergence process} to general L\'evy processes. However, the cost minimization problem relies heavily on the self-similarity property of Brownian motions. Thus, to solve the minimization problem for L\'evy processes, other techniques are needed.

\ACKNOWLEDGMENT{This work is part of the research program Complexity in high-tech manufacturing, (partly) financed by the Dutch Research Council (NWO) through contract 438.16.121. The research is also supported by the NWO programs MEERVOUD [Vlasiou: 632.003.002] and Talent VICI [Zwart: 639.033.413].}

\DeclareRobustCommand{\VAN}[3]{#3}

\bibliographystyle{informs2014} 

\bibliography{refs} 

\begin{thebibliography}{72}
\providecommand{\natexlab}[1]{#1}
\providecommand{\url}[1]{\texttt{#1}}
\providecommand{\urlprefix}{URL }

\bibitem[{Abate \protect\BIBand{} Whitt(1987)}]{abate1987transient}
Abate J, Whitt W (1987) {Transient behavior of regulated Brownian motion, I:
  starting at the origin}. \emph{Advances in Applied Probability}
  19(3):560--598.

\bibitem[{Adler \protect\BIBand{} Taylor(2007)}]{adler2007random}
Adler RJ, Taylor JE (2007) \emph{Random fields and geometry}, volume~80
  (Springer).

\bibitem[{Ak{\c{c}}ay \protect\BIBand{} Xu(2004)}]{akccay2004joint}
Ak{\c{c}}ay Y, Xu SH (2004) Joint inventory replenishment and component
  allocation optimization in an assemble-to-order system. \emph{Management
  Science} 50(1):99--116.

\bibitem[{Altendorfer \protect\BIBand{} Minner(2011)}]{altendorfer2011}
Altendorfer K, Minner S (2011) Simultaneous optimization of capacity and
  planned lead time in a two-stage production system with different customer
  due dates. \emph{European Journal of Operational Research} 213(1):134--146.

\bibitem[{{ASML Holding N.V.}(2021)}]{asmlreport}
{ASML Holding NV} (2021) {ASML} annual report 2020.
  \url{https://www.asml.com/en/investors/annual-report/2020}.

\bibitem[{Asmussen(2003)}]{asmussen2003applied}
Asmussen S (2003) \emph{Applied probability and queues}, volume~2 (Springer).

\bibitem[{Asmussen et~al.(1995)Asmussen, Glynn, \protect\BIBand{}
  Pitman}]{asmussen1995}
Asmussen S, Glynn PW, Pitman J (1995) {Discretization error in simulation of
  one-dimensional reflecting Brownian motion}. \emph{The Annals of Applied
  Probability} 875--896.

\bibitem[{Atan et~al.(2017)Atan, Ahmadi, Stegehuis, de~Kok, \protect\BIBand{}
  Adan}]{atan2017}
Atan Z, Ahmadi T, Stegehuis C, de~Kok T, Adan I (2017) Assemble-to-order
  systems: A review. \emph{European Journal of Operational Research}
  261(3):866--879.

\bibitem[{Atan \protect\BIBand{} Rousseau(2016)}]{atan2016}
Atan Z, Rousseau M (2016) Inventory optimization for perishables subject to
  supply disruptions. \emph{Optimization Letters} 10(1):89--108.

\bibitem[{Atar et~al.(2012)Atar, Mandelbaum, \protect\BIBand{}
  Zviran}]{atar2012control}
Atar R, Mandelbaum A, Zviran A (2012) Control of fork-join networks in heavy
  traffic. \emph{2012 50th Annual Allerton Conference on Communication,
  Control, and Computing (Allerton)}, 823--830 (IEEE).

\bibitem[{Baccelli(1985)}]{baccelli1985two}
Baccelli F (1985) {Two parallel queues created by arrivals with two demands:
  The M/G/2 symmetrical case}. \emph{RR-0426, INRIA. ffinria-00076130} .

\bibitem[{Baccelli \protect\BIBand{} Makowski(1989)}]{baccelli1989queueing}
Baccelli F, Makowski AM (1989) Queueing models for systems with synchronization
  constraints. \emph{Proceedings of the IEEE} 77(1):138--161.

\bibitem[{Bijvank et~al.(2014)Bijvank, Huh, Janakiraman, \protect\BIBand{}
  Kang}]{bijvank2014}
Bijvank M, Huh WT, Janakiraman G, Kang W (2014) Robustness of order-up-to
  policies in lost-sales inventory systems. \emph{Operations Research}
  62(5):1040--1047.

\bibitem[{Bollapragada et~al.(2004)Bollapragada, Rao, \protect\BIBand{}
  Zhang}]{bollapragada2004managing}
Bollapragada R, Rao US, Zhang J (2004) Managing two-stage serial inventory
  systems under demand and supply uncertainty and customer service level
  requirements. \emph{IIE transactions} 36(1):73--85.

\bibitem[{Borst et~al.(2004)Borst, Mandelbaum, \protect\BIBand{}
  Reiman}]{borst2004dimensioning}
Borst S, Mandelbaum A, Reiman MI (2004) Dimensioning large call centers.
  \emph{Operations Research} 52(1):17--34.

\bibitem[{Bradley \protect\BIBand{} Glynn(2002)}]{bradley2002}
Bradley JR, Glynn PW (2002) Managing capacity and inventory jointly in
  manufacturing systems. \emph{Management Science} 48(2):273--288.

\bibitem[{Brown \protect\BIBand{} Resnick(1977)}]{brown1977extreme}
Brown BM, Resnick SI (1977) Extreme values of independent stochastic processes.
  \emph{Journal of Applied Probability} 732--739.

\bibitem[{Chaturvedi \protect\BIBand{} Mart{\'\i}nez-de
  Alb{\'e}niz(2016)}]{chaturvedi2016}
Chaturvedi A, Mart{\'\i}nez-de Alb{\'e}niz V (2016) Safety stock, excess
  capacity or diversification: Trade-offs under supply and demand uncertainty.
  \emph{Production andOperations Management} 25(1):77--95.

\bibitem[{Dai \protect\BIBand{} Harrison(1992)}]{dai1992reflected}
Dai J, Harrison JM (1992) {Reflected Brownian motion in an orthant: numerical
  methods for steady-state analysis}. \emph{The Annals of Applied Probability}
  65--86.

\bibitem[{Denton(2021)}]{barrons}
Denton J (2021) {ASML} cuts guidance in the face of supply-chain issues. the
  chip stock is falling.
  \url{https://www.barrons.com/articles/asml-cuts-guidance-supply-chain-issues-51634728074}.

\bibitem[{D\k{e}bicki et~al.(2015)D\k{e}bicki, Hashorva, Ji, \protect\BIBand{}
  Tabi{\'s}}]{debicki2015extremes}
D\k{e}bicki K, Hashorva E, Ji L, Tabi{\'s} K (2015) {Extremes of vector-valued
  Gaussian processes: Exact asymptotics}. \emph{Stochastic Processes and their
  Applications} 125(11):4039--4065.

\bibitem[{D\k{e}bicki et~al.(2020)D\k{e}bicki, Ji, \protect\BIBand{}
  Rolski}]{debicki2020exact}
D\k{e}bicki K, Ji L, Rolski T (2020) {Exact asymptotics of component-wise
  extrema of two-dimensional Brownian motion}. \emph{Extremes} 23:569–602.

\bibitem[{Do{\u{g}}ru et~al.(2017)Do{\u{g}}ru, Reiman, \protect\BIBand{}
  Wang}]{dougru2017}
Do{\u{g}}ru MK, Reiman MI, Wang Q (2017) Assemble-to-order inventory management
  via stochastic programming: Chained boms and the m-system. \emph{Production
  and Operations Management} 26(3):446--468.

\bibitem[{Ewing \protect\BIBand{} Clark(2021)}]{ewing2021}
Ewing J, Clark D (2021) Lack of tiny parts disrupts auto factories worldwide.
  \emph{New York Times}
  \urlprefix\url{https://www.nytimes.com/2021/01/13/business/auto-factories-semiconductor-chips.html}.

\bibitem[{Flatto \protect\BIBand{} Hahn(1984)}]{flatto1984two}
Flatto L, Hahn S (1984) {Two parallel queues created by arrivals with two
  demands I}. \emph{SIAM Journal on Applied Mathematics} 44(5):1041--1053.

\bibitem[{Gans et~al.(2003)Gans, Koole, \protect\BIBand{}
  Mandelbaum}]{gans2003telephone}
Gans N, Koole G, Mandelbaum A (2003) Telephone call centers: Tutorial, review,
  and research prospects. \emph{Manufacturing \& Service Operations Management}
  5(2):79--141.

\bibitem[{Glasserman(1997)}]{glasserman1997}
Glasserman P (1997) Bounds and asymptotics for planning critical safety stocks.
  \emph{Operations Research} 45(2):244--257.

\bibitem[{Goldberg et~al.(2016)Goldberg, Katz-Rogozhnikov, Lu, Sharma,
  \protect\BIBand{} Squillante}]{goldberg2016asymptotic}
Goldberg DA, Katz-Rogozhnikov DA, Lu Y, Sharma M, Squillante MS (2016)
  Asymptotic optimality of constant-order policies for lost sales inventory
  models with large lead times. \emph{Mathematics of Operations Research}
  41(3):898--913.

\bibitem[{Goldberg et~al.(2021)Goldberg, Reiman, \protect\BIBand{}
  Wang}]{goldberg2021}
Goldberg DA, Reiman MI, Wang Q (2021) A survey of recent progress in the
  asymptotic analysis of inventory systems. \emph{Production and Operations
  Management} 30(6):1718--1750.

\bibitem[{Gopalakrishnan et~al.(2016)Gopalakrishnan, Doroudi, Ward,
  \protect\BIBand{} Wierman}]{gopalakrishnan2016routing}
Gopalakrishnan R, Doroudi S, Ward AR, Wierman A (2016) Routing and staffing
  when servers are strategic. \emph{Operations Research} 64(4):1033--1050.

\bibitem[{{\VAN{Haan}{De}{de}}~Haan \protect\BIBand{}
  Ferreira(2006)}]{de2007extreme}
{\VAN{Haan}{De}{de}}~Haan L, Ferreira A (2006) \emph{Extreme value theory: an
  introduction} (Springer Science \& Business Media).

\bibitem[{Halfin \protect\BIBand{} Whitt(1981)}]{halfin1981heavy}
Halfin S, Whitt W (1981) Heavy-traffic limits for queues with many exponential
  servers. \emph{Operations Research} 29(3):567--588.

\bibitem[{Harrison(1985)}]{harrison1985brownian}
Harrison JM (1985) \emph{Brownian motion and stochastic flow systems} (Wiley
  New York).

\bibitem[{Harrison(2013)}]{harrison2013}
Harrison JM (2013) \emph{Brownian Models of Performance and Control} (Cambridge
  University Press),
  \urlprefix\url{http://dx.doi.org/10.1017/CBO9781139087698}.

\bibitem[{Huh et~al.(2009)Huh, Janakiraman, Muckstadt, \protect\BIBand{}
  Rusmevichientong}]{huh2009}
Huh WT, Janakiraman G, Muckstadt JA, Rusmevichientong P (2009) Asymptotic
  optimality of order-up-to policies in lost sales inventory systems.
  \emph{Management Science} 55(3):404--420.

\bibitem[{Karsten et~al.(2012)Karsten, Slikker, \protect\BIBand{} van
  Houtum}]{karsten2012inventory}
Karsten F, Slikker M, van Houtum GJ (2012) Inventory pooling games for
  expensive, low-demand spare parts. \emph{Naval Research Logistics (NRL)}
  59(5):311--324.

\bibitem[{Klein(1988)}]{de1988fredholm}
Klein SJd (1988) \emph{Fredholm integral equations in queueing analysis}. Ph.D.
  thesis, Rijksuniversiteit Utrecht.

\bibitem[{Klosterhalfen et~al.(2014)Klosterhalfen, Minner, \protect\BIBand{}
  Willems}]{klosterhalfen2014}
Klosterhalfen ST, Minner S, Willems SP (2014) Strategic safety stock placement
  in supply networks with static dual supply. \emph{Manufacturing and Service
  Operations Management} 16(2):204--219.

\bibitem[{Ko \protect\BIBand{} Serfozo(2004)}]{ko2004response}
Ko SS, Serfozo RF (2004) {Response times in M/M/s fork-join networks}.
  \emph{Advances in Applied Probability} 36(3):854--871.

\bibitem[{Kou et~al.(2016)Kou, Zhong et~al.}]{kou2016first}
Kou S, Zhong H, et~al. (2016) {First-passage times of two-dimensional Brownian
  motion}. \emph{Advances in Applied Probability} 48(4):1045--1060.

\bibitem[{Kumar \protect\BIBand{} Randhawa(2010)}]{kumar2010exploiting}
Kumar S, Randhawa RS (2010) Exploiting market size in service systems.
  \emph{Manufacturing \& Service Operations Management} 12(3):511--526.

\bibitem[{Leadbetter et~al.(1983)Leadbetter, Lindgren, \protect\BIBand{}
  Rootz{\'e}n}]{leadbetter2012extremes}
Leadbetter MR, Lindgren G, Rootz{\'e}n H (1983) \emph{Extremes and related
  properties of random sequences and processes} (Springer Science \& Business
  Media).

\bibitem[{{\VAN{Leeuwaarden}{Van}{van}}~Leeuwaarden
  et~al.(2019){\VAN{Leeuwaarden}{Van}{van}}~Leeuwaarden, Mathijsen,
  \protect\BIBand{} Zwart}]{van2019economies}
{\VAN{Leeuwaarden}{Van}{van}}~Leeuwaarden JS, Mathijsen BW, Zwart B (2019)
  Economies-of-scale in many-server queueing systems: Tutorial and partial
  review of the qed halfin--whitt heavy-traffic regime. \emph{SIAM Review}
  61(3):403--440.

\bibitem[{Lu \protect\BIBand{} Pang(2015)}]{lu2015gaussian}
Lu H, Pang G (2015) Gaussian limits for a fork-join network with
  nonexchangeable synchronization in heavy traffic. \emph{Mathematics of
  Operations Research} 41(2):560--595.

\bibitem[{Lu \protect\BIBand{} Pang(2017{\natexlab{a}})}]{lu2017heavy}
Lu H, Pang G (2017{\natexlab{a}}) {Heavy-traffic limits for a fork-join network
  in the Halfin-Whitt regime}. \emph{Stochastic Systems} 6(2):519--600.

\bibitem[{Lu \protect\BIBand{} Pang(2017{\natexlab{b}})}]{lu2017heavy2}
Lu H, Pang G (2017{\natexlab{b}}) Heavy-traffic limits for an infinite-server
  fork--join queueing system with dependent and disruptive services.
  \emph{Queueing Systems} 85(1-2):67--115.

\bibitem[{Lu \protect\BIBand{} Song(2005)}]{lu2005}
Lu Y, Song JS (2005) Order-based cost optimization in assemble-to-order
  systems. \emph{Operations Research} 53(1):151--169.

\bibitem[{Mayorga \protect\BIBand{} Ahn(2011)}]{mayorga2011}
Mayorga ME, Ahn HS (2011) Joint management of capacity and inventory in
  make-to-stock production systems with multi-class demand. \emph{European
  Journal of Operational Research} 212(2):312--324.

\bibitem[{Nair et~al.(2016)Nair, Wierman, \protect\BIBand{}
  Zwart}]{nair2016provisioning}
Nair J, Wierman A, Zwart B (2016) Provisioning of large-scale systems: The
  interplay between network effects and strategic behavior in the user base.
  \emph{Management Science} 62(6):1830--1841.

\bibitem[{Nelson \protect\BIBand{} Tantawi(1988)}]{nelson1988approximate}
Nelson R, Tantawi AN (1988) Approximate analysis of fork/join synchronization
  in parallel queues. \emph{IEEE Transactions on Computers} 37(6):739--743.

\bibitem[{Nguyen(1993)}]{nguyen1993processing}
Nguyen V (1993) Processing networks with parallel and sequential tasks: Heavy
  traffic analysis and brownian limits. \emph{The Annals of Applied
  Probability} 28--55.

\bibitem[{Nguyen(1994)}]{nguyen1994trouble}
Nguyen V (1994) The trouble with diversity: Fork-join networks with
  heterogeneous customer population. \emph{The Annals of Applied Probability}
  1--25.

\bibitem[{Pan \protect\BIBand{} So(2016)}]{pan2016}
Pan W, So KC (2016) Component procurement strategies in decentralized assembly
  systems under supply uncertainty. \emph{IIE Transactions} 48(3):267--282.

\bibitem[{Pickands~III(1968)}]{pickands1968moment}
Pickands~III J (1968) Moment convergence of sample extremes. \emph{The Annals
  of Mathematical Statistics} 39(3):881--889.

\bibitem[{Plambeck(2008)}]{plambeck2008}
Plambeck EL (2008) Asymptotically optimal control for an assemble-to-order
  system with capacitated component production and fixed transport costs.
  \emph{Operations Research} 56(5):1158--1171.

\bibitem[{Plambeck \protect\BIBand{} Ward(2008)}]{plambeckward2008}
Plambeck EL, Ward AR (2008) Optimal control of a high-volume assemble-to-order
  system with maximum leadtime quotation and expediting. \emph{Queueing
  Systems} 60(1):1--69.

\bibitem[{Reddy \protect\BIBand{} Kumar(2020)}]{reddy2020}
Reddy KN, Kumar A (2020) Capacity investment and inventory planning for a
  hybrid manufacturing-remanufacturing system in the circular economy.
  \emph{International Journal of Production Research} 1--29.

\bibitem[{Reed \protect\BIBand{} Zhang(2017)}]{reed2017}
Reed J, Zhang B (2017) Managing capacity and inventory jointly for multi-server
  make-to-stock queues. \emph{Queueing Systems} 86(1-2):61--94.

\bibitem[{Reiman \protect\BIBand{} Wang(2015)}]{reiman2015}
Reiman MI, Wang Q (2015) Asymptotically optimal inventory control for
  assemble-to-order systems with identical lead times. \emph{Operations
  Research} 63(3):716--732.

\bibitem[{R{\'e}nyi(1953)}]{renyi1953theory}
R{\'e}nyi A (1953) On the theory of order statistics. \emph{Acta Math. Acad.
  Sci. Hung} 4(2).

\bibitem[{Resnick(1987)}]{resnick2013extreme}
Resnick SI (1987) \emph{Extreme values, regular variation and point processes}
  (Springer).

\bibitem[{Sleptchenko et~al.(2003)Sleptchenko, van~der Heijden,
  \protect\BIBand{} van Harten}]{sleptchenko2003}
Sleptchenko A, van~der Heijden MC, van Harten A (2003) Trade-off between
  inventory and repair capacity in spare part networks. \emph{Journal of the
  Operational Research Society} 54(3):263--272.

\bibitem[{Song(1998)}]{song1998}
Song JS (1998) On the order fill rate in a multi-item, base-stock inventory
  system. \emph{Operations research} 46(6):831--845.

\bibitem[{{\VAN{Vaart}{Van der}{van der}}~Vaart(1998)}]{vaart_1998}
{\VAN{Vaart}{Van der}{van der}}~Vaart AW (1998) \emph{Asymptotic Statistics}.
  Cambridge Series in Statistical and Probabilistic Mathematics (Cambridge
  University Press),
  \urlprefix\url{http://dx.doi.org/10.1017/CBO9780511802256}.

\bibitem[{Varma(1990)}]{varma1990heavy}
Varma S (1990) \emph{Heavy and light traffic approximations for queues with
  synchronization constraints}. Ph.D. thesis, University of Maryland.

\bibitem[{Wright(1992)}]{wright1992two}
Wright PE (1992) Two parallel processors with coupled inputs. \emph{Advances in
  Applied Probability} 24(4):986--1007.

\bibitem[{Wu \protect\BIBand{} Chao(2014)}]{wu2014}
Wu J, Chao X (2014) {Optimal control of a Brownian production/inventory system
  with average cost criterion}. \emph{Mathematics of Operations Research}
  39(1):163--189.

\bibitem[{Xin \protect\BIBand{} Goldberg(2016)}]{xin2016}
Xin L, Goldberg DA (2016) Optimality gap of constant-order policies decays
  exponentially in the lead time for lost sales models. \emph{Operations
  Research} 64(6):1556--1565.

\bibitem[{Xin \protect\BIBand{} Goldberg(2018)}]{xin2018asymptotic}
Xin L, Goldberg DA (2018) Asymptotic optimality of tailored base-surge policies
  in dual-sourcing inventory systems. \emph{Management Science} 64(1):437--452.

\bibitem[{Zhang et~al.(2020)Zhang, Zhang, \protect\BIBand{}
  Zhang}]{zhang2020simple}
Zhang H, Zhang J, Zhang RQ (2020) Simple policies with provable bounds for
  managing perishable inventory. \emph{Production and Operations Management}
  29(11):2637--2650.

\bibitem[{Zieli{\'n}ski(2009)}]{zielinski2009}
Zieli{\'n}ski R (2009) {Optimal nonparametric quantile estimators. Towards a
  general theory. A survey}. \emph{Communications in Statistics-Theory and
  Methods} 38(7):980--992.

\bibitem[{Zou et~al.(2004)Zou, Pokharel, \protect\BIBand{} Piplani}]{zou2004}
Zou X, Pokharel S, Piplani R (2004) Channel coordination in an assembly system
  facing uncertain demand with synchronized processing time and delivery
  quantity. \emph{International Journal of Production Research}
  42(22):4673--4689.

\end{thebibliography}




\newpage
\begin{APPENDIX}{E-companion}
\section{Proofs}
\label{sec: appendix}
\subsection{Proofs of Section \ref{sec: problem formulation}}\label{app: proofs problem formulation}

\proof{Proof of Lemma \ref{lem: inventory minimization}}
$F_N(I,\beta)>0$, hence $F_N$ has a global infimum, and since $\lim_{\beta\downarrow 0}F_N(I,\beta)=\infty$, $\lim_{\beta\to \infty}F_N(I,\beta)=\infty$ and $ \lim_{I\to \infty}F_N(I,\beta)=\infty$, $F_N$ has a global minimum. Now, assume $F_N(I_N,\beta_N)=\min_{(I,\beta)}F_N(I,\beta)$. Assume that there exists an $\hat{I}_N$ such that 
\begin{multline*}
    \InventoryExpCosts{\hat{I}_N}\\
    < \InventoryExpCosts{I_N}.
\end{multline*}
Then $F_N(\hat{I}_N,\beta_N)< F_N(I_N,\beta_N)$. This contradicts the statement that $(I_N,\beta_N)$ gives the minimum of $F_N$. Hence, the optimal base-stock level minimizes $C_N(I)$. The proof that $\beta_N$ minimizes $\frac{1}{\beta}C_N(I_N)+\beta N$ goes analogously. 

To prove that $C_N(I)$ is convex with respect to $I$, we observe that
\begin{align*}
    \frac{d^2}{dI^2}C_N(I)=\left(b^{(N)}+Nh^{(N)}\right)\frac{d^2}{dI^2}\expect*{\left(\max_{i\leq N}Q_i-I\right)^+}=&\left(b^{(N)}+Nh^{(N)}\right)\frac{d^2}{dI^2}\int_{I}^{\infty} \probability*{\max_{i\leq N}Q_i>x}dx\nonumber\\
    =&\left(b^{(N)}+Nh^{(N)}\right)f(I)\geq 0,
\end{align*}
because $f$ is the probability density function of $\max_{i\leq N}Q_i$. This density exists; cf.\ \citet[Prop. 2a]{dai1992reflected}. In conclusion, we have a convex minimization problem. Moreover, $\frac{d^2}{d\beta^2}\left(\frac{1}{\beta}C_N(I_N)+\beta N\right)=\frac{2}{\beta^3}C_N(I_N)>0$. Thus $\frac{1}{\beta}C_N(I_N)+\beta N$ is also convex with respect to $\beta$.
\hfill\rlap{\hspace*{-2em}\Halmos}
\endproof 

\proof{Proof of Lemma \ref{lem: optimal beta general}}
$F_N(I,\beta)$ has the form $F_N(I,\beta)=\frac{1}{\beta}C_N(I)+\beta N$, thus in order to minimize $F_N(I^*_N,\beta)$, we know by Lemma \ref{lem: inventory minimization} that we need to solve $\frac{d}{d\beta}F_N(I^*_N,\beta)=-\frac{1}{\beta^2}C_N(I^*_N)+N=0$. Thus, $\beta^*_N=\frac{\sqrt{C_N(I^*_N)}}{\sqrt{N}}$, and $F_N(I^*_N,\beta^*_N)=2\sqrt{NC_N(I^*_N)}=2N\beta^*_N$.
\hfill\rlap{\hspace*{-2em}\Halmos}
\endproof 

\proof{Proof of Lemma \ref{lem: optimal I general}}
To solve $\min_I C_N(I)$ we have to solve $\frac{d}{dI}C_N(I)=0$, this gives for the optimal base-stock level $I^*_N$ that
\begin{align*}
    Nh^{(N)}-\left(Nh^{(N)}+b^{(N)}\right)\probability*{\max_{i\leq N} Q_i>I^*_N}=0.
\end{align*}
Hence $I^*_N=P^{-1}_N\left(\frac{b^{(N)}}{Nh^{(N)}+b^{(N)}}\right)$, with $P^{-1}_N$ the quantile function of $\max_{i\leq N}Q_i$.
\hfill\rlap{\hspace*{-2em}\Halmos}
\endproof 
\proof{Proof of Lemma \ref{lem: first order bal reg}}
Following Corollary \ref{cor: ratio optimal approx}, we have 
\begin{align*}
   \frac{F_N(I^*_N,\beta^*_N)}{F_N(\bar I_N, \bar \beta_N)}=\frac{2\sqrt{C_N(I^*_N)}\sqrt{\bar{C}_N(\bar{I}_N)}}{C_N(\bar{I}_N)+\bar{C}_N(\bar{I}_N)}.
\end{align*}
Furthermore, observe that 
\begin{align*}
    \expect*{\max_{i\leq N}Q_i}\geq \expect*{\max_{i\leq N}\sup_{s>0}(W_i(s)-s)+W_A(\tau)}=\frac{\sigma^2}{2}\sum_{i=1}^N\frac{1}{i}\geq \frac{\sigma^2}{2}\log N,
\end{align*}
where $\tau$ is the first hitting time of the supremum of $\max_{i\leq N}(W_i(t)-t)$. From this it follows that for $I<\frac{\sigma^2}{2}\log N$, $ \frac{\sigma^2}{2}\log N-I<\expect*{\max_{i\leq N}Q_i-I}<\expect*{\left(\max_{i\leq N}Q_i-I\right)^+}.$
For $I>\frac{\sigma^2}{2}\log N$, $(\frac{\sigma^2}{2}\log N-I)^+=0<\expect*{(\max_{i\leq N}Q_i-I)^+}$. In conclusion, $C_N(I)>\bar{C}_N(I)$. Therefore,
\begin{align*}
     \frac{F_N(I^*_N,\beta^*_N)}{F_N(\bar I_N, \bar \beta_N)}=\frac{2\sqrt{C_N(I^*_N)}\sqrt{\bar{C}_N(\bar{I}_N)}}{C_N(\bar{I}_N)+\bar{C}_N(\bar{I}_N)}\geq \frac{\sqrt{C_N(I^*_N)}\sqrt{\bar{C}_N(\bar{I}_N)}}{C_N(\bar{I}_N)}. 
\end{align*}
We have $|C_N(I^*_N)-C_N(\bar I_N)|\leq (2Nh^{(N)}+b^{(N)})|I^*_N-\bar I_N|,$ and
\begin{align*}
    |\bar C_N(\bar I_N)-C_N(\bar I_N)|\leq (Nh^{(N)}+b^{(N)})\expect*{\left|\max_{i\leq N}Q_i-\frac{\sigma^2}{2}\log N\right|}.
\end{align*}
In the case that $\gamma_N=\gamma\in(0,1)$, we have by applying Lemma \ref{lem: first order EVT} that $|\bar C_N(\bar I_N)-C_N(\bar I_N)|=o((Nh^{(N)}+b^{(N)})\log N)$. Furthermore, $C_N(\bar I_N)\sim Nh^{(N)}\frac{\sigma^2}{2}\log N$, and since $\max_{i\leq N}Q_i/\log N\LimitP \sigma^2/2$, as $N\to\infty$, we also have that $I^*_N/\log N\LimitN \sigma^2/2$. Thus $|C_N(I^*_N)-C_N(\bar I_N)|=o((Nh^{(N)}+b^{(N)})\log N)$, and  the lemma follows.

In the case that $\gamma_N\LimitN 1$, we first observe that $\bar C_N(\bar I_N)=Nh^{(N)}\left(\frac{\sigma^2}{2}\log N-\frac{\sigma^2+\sigma_A^2}{2}\right)\sim Nh^{(N)}\frac{\sigma^2}{2}\log N $. Furthermore,
\begin{align*}
    C_N(\bar I_N)=&Nh^{(N)}\left(\frac{\sigma^2}{2}\log N-\frac{\sigma^2+\sigma_A^2}{2}\right)+(Nh^{(N)}+b^{(N)})\expect*{\left(\max_{i\leq N}Q_i-\frac{\sigma^2}{2}\log N\right)^+}\nonumber\\
    \leq&Nh^{(N)}\left(\frac{\sigma^2}{2}\log N-\frac{\sigma^2+\sigma_A^2}{2}\right)+(Nh^{(N)}+b^{(N)})\expect*{\left|\max_{i\leq N}Q_i-\frac{\sigma^2}{2}\log N\right|}.
\end{align*}
Thus,
\begin{align*}
    \frac{C_N(\bar I_N)}{Nh^{(N)}\log N}\leq \frac{\sigma^2}{2}+o(1)+\frac{1}{\gamma_N}\frac{\expect*{\left|\max_{i\leq N}Q_i-\frac{\sigma^2}{2}\log N\right|}}{\log N}.
\end{align*}
By Lemma \ref{lem: first order EVT}, we know that $\expect*{\left|\max_{i\leq N}Q_i-\frac{\sigma^2}{2}\log N\right|}/\log N\LimitN 0$. Thus 
\begin{align*}
\limsup_{N\to\infty}C_N(\bar I_N)/(Nh^{(N)}\log N)\leq \sigma^2/2.
\end{align*}
Finally,
\begin{align*}
    C_N(I^*_N)=&Nh^{(N)}\left(I^*_N-\frac{\sigma^2+\sigma_A^2}{2}\right)+(Nh^{(N)}+b^{(N)})\expect*{\left(\max_{i\leq N}Q_i-I^*_N\right)^+}\nonumber\\
    \geq&Nh^{(N)}\left(I^*_N-\frac{\sigma^2+\sigma_A^2}{2}\right)+(Nh^{(N)}+b^{(N)})\expect*{\max_{i\leq N}Q_i-I^*_N}\nonumber\\
    \geq&-Nh^{(N)}\frac{\sigma^2+\sigma_A^2}{2}+(Nh^{(N)}+b^{(N)})\frac{\sigma^2}{2}\log N-b^{(N)}I^*_N.
\end{align*}
$I^*_N=O(\log N)$, and $b^{(N)}/(Nh^{(N)})\LimitN 0$, therefore, $\liminf_{N\to\infty}C_N(I^*_N)/(Nh^{(N)}\log N)\geq \sigma^2/2$. Combining these results gives 
\begin{align*}
    \liminf_{N\to\infty}\frac{F_N(I^*_N,\beta^*_N)}{F_N(\bar I_N, \bar \beta_N)}\geq\liminf_{N\to\infty} \frac{\sqrt{C_N(I^*_N)}\sqrt{\bar{C}_N(\bar{I}_N)}}{C_N(\bar{I}_N)}=1.
\end{align*}
\hfill\rlap{\hspace*{-2em}\Halmos}
\endproof 
\subsection{Proofs of Section \ref{sec: basic model - independent}}\label{app: proofs basic model - independent}
\proof{Proof of Lemma \ref{lem: minimization problem case WA=0}}
In Lemma \ref{lem: optimal I general}, it is shown that $I^*_N=P^{-1}_N(1-\gamma_N)$, with $P^{-1}_N$ the quantile function of $\max_{i\leq N}Q_i$. Because $(Q_i,i\leq N)$ are independent and exponentially distributed,
\begin{align*}
    \probability*{\max_{i\leq N}Q_i\leq P^{-1}_N(x)}=x=\left(1-e^{-\frac{2}{\sigma^2}P^{-1}_N(x)}\right)^N.
\end{align*}
From this it follows that $P^{-1}_N(x)=\frac{\sigma^2}{2}\log\Big(1\Big/\Big(1-x^{\frac{1}{N}}\Big)\Big)$.
\hfill\rlap{\hspace*{-2em}\Halmos}
\endproof 

\proof{Proof of Proposition \ref{prop: approx case WA=0}}
Minimizing $\ApproxCosts(\ApproxIN,\ApproxBetaN)$ goes analogously as minimizing $F_N(I_N,\beta_N)$ in Lemma \ref{lem: minimization problem case WA=0}.
Hence $ \ApproxIN=\ApproxQuantile(1-\gamma_N).$ Thus, we have to solve 
\begin{align*}
    \probability*{\frac{\sigma^2}{2}G+\frac{\sigma^2}{2}\log N\leq \ApproxQuantile(x)}=\probability*{G\leq \frac{2}{\sigma^2}\ApproxQuantile(x)-\log N}=e^{-e^{-\left(\frac{2}{\sigma^2}\ApproxQuantile(x)-\log N\right)}}=x.
\end{align*}
Therefore, $\ApproxQuantile(x)=\frac{\sigma^2}{2}\log N-\frac{\sigma^2}{2}\log(-\log x).$ Hence, the optimal base-stock level is given in Equation \eqref{eq: approx inventory case WA=0}. Furthermore,
\begin{align*}
    \expect*{\left(\frac{\sigma^2}{2}G+\frac{\sigma^2}{2}\log N-\ApproxIN\right)^+}
    =& \expect*{\left(\frac{\sigma^2}{2}G+\frac{\sigma^2}{2}\log\left(-\log(1-\gamma_N)\right)\right)^+}\nonumber\\
    =&\frac{\sigma^2}{2}\int_{-\log\left(-\log(1-\gamma_N)\right)}^{\infty}1-e^{-e^{-x}}dx.
    \end{align*}
By using partial integration and substitution we can write  
    \begin{align*}
   &\frac{\sigma^2}{2}\int_{-\log\left(-\log(1-\gamma_N)\right)}^{\infty}1-e^{-e^{-x}}dx
   =\frac{\sigma^2}{2}\left( \int_{-\log (1-\gamma_N)}^{\infty}\frac{e^{-t}}{t}dt+\Gamma   +\log\left(-\log(1-\gamma_N)\right)\right).
\end{align*}
Hence, this gives us the expression of $\ApproxC{\ApproxIN}$ in \eqref{eq: approx capacity case WA=0}.
\hfill\rlap{\hspace*{-2em}\Halmos}
\endproof
\begin{lem}\label{lem: G and Q same prob space}
Define 
\begin{align}
    G_N:=-\log\left(-\log\left(\left(1-\exp\left(-\frac{2}{\sigma^2}\max_{i\leq N}Q_i\right)\right)^N\right)\right),
\end{align}
then $\probability*{G_N<x}=e^{-e^{-x}}$, for all $N$. Moreover,
\begin{align}\label{eq: max Qi> G}
    \max_{i\leq N}Q_i>\frac{\sigma^2}{2}G_N+\frac{\sigma^2}{2}\log N,
\end{align}
and $\max_{i\leq N}Q_i-\frac{\sigma^2}{2}G_N-\frac{\sigma^2}{2}\log N$ strictly decreases as a function of $\max_{i\leq N} Q_i$ with limit 0.
\end{lem}
\proof{Proof} 
To prove that $G_N$ follows a Gumbel distribution, we first observe that $\probability*{\max_{i\leq N}Q_i<x}=\left(1-\exp\left(-\frac{2}{\sigma^2}x\right)\right)^N.$ Therefore, $\left(1-\exp\left(-\frac{2}{\sigma^2}\max_{i\leq N}Q_i\right)\right)^N\sim\text{Unif}[0,1]$.
Then,
\begin{align*}
    \probability*{G_N<x}=&\probability*{-\log\left(-\log\left(\left(1-\exp\left(-\frac{2}{\sigma^2}\max_{i\leq N}Q_i\right)\right)^N\right)\right)<x}\nonumber\\
    =&\probability*{-\log\left(\left(1-\exp\left(-\frac{2}{\sigma^2}\max_{i\leq N}Q_i\right)\right)^N\right)>e^{-x}}\nonumber\\
    =&\probability*{\left(1-\exp\left(-\frac{2}{\sigma^2}\max_{i\leq N}Q_i\right)\right)^N<e^{-e^{-x}}}=e^{-e^{-x}}.
\end{align*}
To prove \eqref{eq: max Qi> G}, we need to show that for all $x>0$ and $N$ 
\begin{align*}
    x>&-\frac{\sigma^2}{2}\log\left(-\log\left(\left(1-\exp\left(-\frac{2}{\sigma^2}x\right)\right)^N\right)\right)+\frac{\sigma^2}{2}\log N.
\end{align*}
This is equivalent to the inequality $x>-\frac{\sigma^2}{2}\log\left(-\log\left(1-\exp\left(-\frac{2}{\sigma^2}x\right)\right)\right)$, which is equivalent to $  1-e^{-\frac{2}{\sigma^2}x}<e^{-e^{-\frac{2}{\sigma^2}x}}$, with $x>0$. This is equivalent to $e^{-y}>1-y$ for $y\in(0,e^{-1}]$. Observe that for $y=0$, we have equality, and we have for $y>0$ that $(e^{-y})'>-1=(1-y)'$. The statement follows.
To prove that the larger $\max_{i\leq N}Q_i$ becomes, the smaller the difference between $\max_{i\leq N}Q_i$ and $\frac{\sigma^2}{2}G_N+\frac{\sigma^2}{2}\log N$ becomes, we first observe that
\begin{align*}
    \frac{\sigma^2}{2}G_N+\frac{\sigma^2}{2}\log N=&-\frac{\sigma^2}{2}\log\left(-\log\left(\left(1-\exp\left(-\frac{2}{\sigma^2}\max_{i\leq N}Q_i\right)\right)^N\right)\right)+\frac{\sigma^2}{2}\log N\nonumber\\
    =&-\frac{\sigma^2}{2}\log\left(-\log\left(1-e^{-\frac{2}{\sigma^2}\max_{i\leq N}Q_i}\right)\right).
\end{align*}
Thus we need to obtain that $x+\frac{\sigma^2}{2}\log(-\log(1-e^{-\frac{2}{\sigma^2}x}))$ is strictly decreasing in $x$ for $x>0$. Taking the first derivative gives the inequality
\begin{align*}
    \frac{e^{-\frac{2 x}{\sigma ^2}}}{\left(1-e^{-\frac{2 x}{\sigma ^2}}\right) \log \left(1-e^{-\frac{2 x}{\sigma ^2}}\right)}+1<0.
\end{align*}
This is equivalent to the inequality $-y/((1-y)\log(1-y))>1$ for $y\in(0,1)$, which can be rewritten to $\log y>1-1/y$, which is a basic logarithm inequality. Finally, $\lim_{x\to\infty}x+\frac{\sigma^2}{2}\log(-\log(1-e^{-\frac{2}{\sigma^2}x}))=0$.
\hfill\rlap{\hspace*{-2em}\Halmos}
\endproof
\begin{lem}\label{lem: bounds C Ctilde}
Let $\gamma_N=Nh^{(N)}/(Nh^{(N)}+b^{(N)})$, then
\begin{align}\label{eq: ineq C(I) C(Itilde)}
    \left|C_N(I^*_N)-C_N(\ApproxIN)\right|\leq (I^*_N-\ApproxIN)(Nh^{(N)}+b^{(N)})\left(1-\gamma_N-\left(1+\frac{\log(1-\gamma_N)}{N}\right)^N\right),
\end{align}
\begin{align}\label{eq: ineq C Ctilde}
   \left|\ApproxC{\ApproxIN}-C_N(\ApproxIN)\right|\leq (I^*_N-\ApproxIN)Nh^{(N)}\left(1-\left(1+\frac{\log(1-\gamma_N)}{N}\right)^N\right).
\end{align}
\end{lem}
\proof{Proof}
Due to the inequality in \eqref{eq: max Qi> G}, $I^*_N>\ApproxIN$, then, we have 
\begin{align*}
    C_N(I^*_N)-C_N(\ApproxIN)=&Nh^{(N)}(I^*_N-\ApproxIN)+(Nh^{(N)}+b^{(N)})\expect*{\left(\max_{i\leq N}Q_i-I^*_N\right)^+-\left(\max_{i\leq N}Q_i-\ApproxIN\right)^+}\nonumber\\
    =&Nh^{(N)}(I^*_N-\ApproxIN)+(Nh^{(N)}+b^{(N)})\expect*{\bigg(\ApproxIN-I^*_N\bigg)\mathbbm{1}\left(\max_{i\leq N}Q_i>I^*_N\right)}\nonumber\\
    -&(Nh^{(N)}+b^{(N)})\expect*{\left(\max_{i\leq N}Q_i-\ApproxIN\right)^+\mathbbm{1}\left(\ApproxIN<\max_{i\leq N}Q_i<I^*_N\right)}.
\end{align*}
We have $\probability*{\max_{i\leq N}Q_i>I^*_N}=\gamma_N=Nh^{(N)}/(Nh^{(N)}+b^{(N)})$, thus
\begin{align*}
    Nh^{(N)}(I^*_N-\ApproxIN)+(Nh^{(N)}+b^{(N)})\expect*{(\ApproxIN-I^*_N)\mathbbm{1}\left(\max_{i\leq N}Q_i>I^*_N\right)}=0.
\end{align*}
Furthermore,
\begin{align*}
   \expect*{\left(\max_{i\leq N}Q_i-\ApproxIN\right)^+\mathbbm{1}\left(\ApproxIN<\max_{i\leq N}Q_i<I^*_N\right)}\leq& (I^*_N-\ApproxIN)\probability*{\ApproxIN<\max_{i\leq N}Q_i<I^*_N}\nonumber\\
    =&(I^*_N-\ApproxIN)\left(1-\gamma_N-\left(1+\frac{\log(1-\gamma_N)}{N}\right)^N\right).
\end{align*}
Equation \eqref{eq: ineq C(I) C(Itilde)} follows. To prove Equation \eqref{eq: ineq C Ctilde}, we observe that 
\begin{align}
     &|\ApproxC{\ApproxIN}-C_N(\ApproxIN)|\nonumber\\
     &\quad=(Nh^{(N)}+b^{(N)})\expect*{\left(\max_{i\leq N}Q_i-\ApproxIN\right)^+-\left(\frac{\sigma^2}{2}G_N+\frac{\sigma^2}{2}\log N-\ApproxIN\right)^+}\nonumber\\
     &\quad=(Nh^{(N)}+b^{(N)})\expect*{\left(\max_{i\leq N}Q_i-\frac{\sigma^2}{2}G_N-\frac{\sigma^2}{2}\log N\right)\mathbbm{1}\left(\frac{\sigma^2}{2}G_N+\frac{\sigma^2}{2}\log N>\ApproxIN\right)}\\
     &\qquad+(Nh^{(N)}+b^{(N)})\expect*{\left(\max_{i\leq N}Q_i-\ApproxIN\right)\mathbbm{1}\left(\frac{\sigma^2}{2}G_N+\frac{\sigma^2}{2}\log N<\ApproxIN<\max_{i\leq N}Q_i\right)}\label{subeq: Ctilde - C}.
\end{align}
Because $G_N$ and $\max_{i\leq N}Q_i$ are on the same probability space, we have $\mathbb{P}\Big(\max_{i\leq N}Q_i=I^*_N\Big| \frac{\sigma^2}{2}G_N+\frac{\sigma^2}{2}\log N=\ApproxIN\Big)=1.$ Furthermore, $x+\frac{\sigma^2}{2}\log(-\log(1-e^{-\frac{2}{\sigma^2}x}))$ is decreasing in $x$. Thus, we can bound 
\begin{align}
    &\expect*{\bigg(\max_{i\leq N}Q_i-\frac{\sigma^2}{2}G_N-\frac{\sigma^2}{2}\log N\bigg)\mathbbm{1}\left(\frac{\sigma^2}{2}G_N+\frac{\sigma^2}{2}\log N>\ApproxIN\right)}\nonumber\\
    &\quad\leq (I^*_N-\ApproxIN)\probability*{\frac{\sigma^2}{2}G_N+\frac{\sigma^2}{2}\log N>\ApproxIN}\nonumber\\
    &\quad=(I^*_N-\ApproxIN)\gamma_N.\label{subeq: Ctilde - C bound 1}
\end{align}
Similarly, for \eqref{subeq: Ctilde - C}, we observe that if $\frac{\sigma^2}{2}G_N+\frac{\sigma^2}{2}\log N<\ApproxIN$, then $\max_{i\leq N}Q_i<I^*_N$, thus,
\begin{align}
    &\expect*{\bigg(\max_{i\leq N}Q_i-\ApproxIN\bigg)\mathbbm{1}\left(\frac{\sigma^2}{2}G_N+\frac{\sigma^2}{2}\log N<\ApproxIN<\max_{i\leq N}Q_i\right)}\nonumber\\
    &\quad\leq(I^*_N-\ApproxIN)\probability*{\frac{\sigma^2}{2}G_N+\frac{\sigma^2}{2}\log N<\ApproxIN<\max_{i\leq N}Q_i}\nonumber\\
    &\quad\leq(I^*_N-\ApproxIN)\left(1-\left(1+\frac{\log(1-\gamma_N)}{N}\right)^N-\gamma_N\right).\label{subeq: Ctilde - C bound 2}
\end{align}
Adding the bounds in \eqref{subeq: Ctilde - C bound 1} and \eqref{subeq: Ctilde - C bound 2} gives the result.
\hfill\rlap{\hspace*{-2em}\Halmos}
\endproof
\proof{Proof of Theorem \ref{thm: order bounds indep}}
First of all, we assume that $\gamma_N=\gamma\in(0,1)$. Using Corollary \ref{cor: ratio optimal approx}, we have
\begin{align*}
    \frac{F_N(I^*_N,\beta^*_N)}{F_N(\ApproxIN,\ApproxBetaN)}=\frac{2\sqrt{C_N(I^*_N)}\sqrt{\ApproxC{\ApproxIN}}}{C_N(\ApproxIN)+\ApproxC{\ApproxIN}}.
\end{align*}
Because of the inequality in \eqref{eq: max Qi> G}, we have for all $I$ that $C_N(I)>\ApproxC{I}$, thus
\begin{align*}
    \frac{F_N(I^*_N,\beta^*_N)}{F_N(\ApproxIN,\ApproxBetaN)}>\frac{2\sqrt{C_N(I^*_N)}\sqrt{\ApproxC{\ApproxIN}}}{2C_N(\ApproxIN)}.
\end{align*}
We write $f(x):= I^*_{1/x}-\hat{I}_{1/x}$ for $x>0$. Then, we have that
\begin{align*}
   f(x)=&\frac{\sigma^2}{2}\log \left(\frac{1}{1-(1-\gamma )^x}\right)+\frac{\sigma^2}{2}\log x +\frac{\sigma^2}{2}\log(-\log(1-\gamma)))\\
   =&\frac{\sigma^2}{2}\log \left(\frac{x}{1-(1-\gamma )^x}\right)+\frac{\sigma^2}{2}\log(-\log(1-\gamma))).
\end{align*}
By first noting that $x/(1-(1-\gamma)^x)=1/(1-e^{-x \log(1-\gamma))})\rightarrow 1/(- \log (1-\gamma))>0$, we see that $\log (x/(1-(1-\gamma )^x))\longrightarrow -\log(-\log(1-\gamma))$ as $x\downarrow 0$. From this, it follows that $f(x)\longrightarrow 0$ as $x\downarrow 0$ and we can extend the domain of the function $f$ such that $f(0):=0$ and that $f$ is twice differentiable at $x=0$.
By computing the Taylor series of the function $f$ at $x=0$, we get 
\begin{align*}
   f(x)=-\frac{\sigma^2}{4} x \log (1-\gamma )+O(x^2).
\end{align*}
Thus, $(I^*_N-\ApproxIN)\sim -\sigma^2\log(1-\gamma)/(4N),$ as $N\to\infty$.
Following \eqref{eq: ineq C Ctilde}, we can conclude that $|\ApproxC{\ApproxIN}-C_N(\ApproxIN)|/(Nh^{(N)})=O(1/N).$
We can do the same for $\probability*{\ApproxIN<\max_{i\leq N}Q_i<I^*_N}$, and get
\begin{align*}
    \left(1-\gamma-\left(1+\frac{\log(1-\gamma)}{N}\right)^N\right)\sim\frac{1}{2N}(1-\gamma)\log(1-\gamma)^2.
\end{align*}
Thus, after applying the inequality in \eqref{eq: ineq C(I) C(Itilde)}, we get $|C_N(I^*_N)-C_N(\ApproxIN)|/(Nh^{(N)}+b^{(N)})=O(1/N^2)$. We have
\begin{align*}
    \ApproxC{\ApproxIN}=&Nh^{(N)}\frac{\sigma^2}{2}(\log N-\log(-\log(1-\gamma))-1)+(Nh^{(N)}+b^{(N)})\frac{\sigma^2}{2}\expect*{(G+\log(-\log(1-\gamma)))^+}\nonumber\\
    \sim& Nh^{(N)}\frac{\sigma^2}{2}\log N,
\end{align*}
because $(Nh^{(N)}+b^{(N)})/(Nh^{(N)})=1/\gamma$, and $-\log(-\log(1-\gamma))$ and $\expect*{(G_N+\log(-\log(1-\gamma)))^+}$ are of $O(1)$.
In conclusion, we have
\begin{align*}
    \frac{F_N(I^*_N,\beta^*_N)}{F_N(\ApproxIN,\ApproxBetaN)}>&\frac{\sqrt{C_N(I^*_N)}}{\sqrt{C_N(\ApproxIN)}}\frac{\sqrt{\ApproxC{\ApproxIN}}}{\sqrt{C_N(\ApproxIN)}}\nonumber\\
    =&\frac{\sqrt{C_N(\ApproxIN)-O((Nh^{(N)}+b^{(N)})/N^2)}}{\sqrt{C_N(\ApproxIN)}}\frac{\sqrt{C_N(\ApproxIN)-O(Nh^{(N)}/N)}}{\sqrt{C_N(\ApproxIN)}}\nonumber\\
    =&\sqrt{1-O(1/(N^2\log N))}\sqrt{1-O(1/(N\log N))}\nonumber\\
    =&1-O(1/(N\log N)).
\end{align*}

Now, we assume that $\gamma_N\LimitN 0$, then we have that $-\log(-\log(1-\gamma_N))\sim -\log(\gamma_N)$, thus $\ApproxIN\sim\frac{\sigma^2}{2} \log (N/\gamma_N)$. Also, \begin{align*}
    \expect*{(G_N+\log(-\log(1-\gamma_N)))^+}\sim\expect*{(G_N+\log(\gamma_N))^+}\sim\gamma_N.
\end{align*}
From this it follows that $\ApproxC{\ApproxIN}\sim Nh^{(N)}\frac{\sigma^2}{2}\log (N/\gamma_N)$. Furthermore,
\begin{align*}
   \probability*{\max_{i\leq N}Q_i>\ApproxIN}= 1-\left(1+\frac{\log(1-\gamma_N)}{N}\right)^N\leq N\probability*{Q_i>\ApproxIN}= -\log(1-\gamma_N)=\gamma_N(1+O(\gamma_N/2)).
\end{align*}
From this it follows that
\begin{align*}
    \left(1-\gamma_N-\left(1+\frac{\log(1-\gamma_N)}{N}\right)^N\right)\leq -\log(1-\gamma_N)-\gamma_N=\frac{\gamma_N^2}{2}(1+o(1)).
\end{align*}
Also 
\begin{align*}
    \probability*{\max_{i\leq N}Q_i<I^*_N}=\probability*{\frac{\sigma^2}{2}G_N+\frac{\sigma^2}{2}\log N<\ApproxIN}=1-\gamma_N\LimitN 1.
\end{align*}
Earlier, we showed that when $\gamma_N=\gamma$, $(I^*_N-\ApproxIN)=O(1/N)$, now $I^*_N$ is larger, because $\probability*{\max_{i\leq N}Q_i<I^*_N}=1-\gamma_N\LimitN 1$. Following the statement in Lemma \ref{lem: G and Q same prob space} that the difference between $\max_{i\leq N} Q_i$ and $\frac{\sigma^2}{2}G_N+\frac{\sigma^2}{2}\log N$ decreases as $\max_{i\leq N}Q_i$ increases, we can conclude that $(I^*_N-\ApproxIN)=O(1/N).$ Following the proof before, and by using the order bounds in \eqref{eq: ineq C(I) C(Itilde)} and \eqref{eq: ineq C Ctilde}, we have that 
\begin{align*}
      \frac{F_N(I^*_N,\beta^*_N)}{F_N(\ApproxIN,\ApproxBetaN)}=1-O(\gamma_N/(N\log(N/\gamma_N))).
\end{align*}

Finally, we consider the case that $\gamma_N\LimitN 1$ and $\gamma_N\leq 1-\exp(-N)$. Then, $\ApproxIN\geq 0$. Furthermore, when $\gamma_N\LimitN 1$, we have $\log(-\log(1-\gamma_N))\LimitN \infty$, from this it follows that
\begin{align*}
    \expect*{(G_N+\log(-\log(1-\gamma_N)))^+}\sim \log(-\log(1-\gamma_N)).
\end{align*}
Thus
\begin{align*}
    \ApproxC{\ApproxIN}\sim&\frac{\sigma^2}{2}Nh^{(N)}(\log N-\log(-\log(1-\gamma_N)))+\frac{\sigma^2}{2}(Nh^{(N)}+b^{(N)})\log(-\log(1-\gamma_N))\nonumber\\
    =&\frac{\sigma^2}{2}Nh^{(N)}\log N+\frac{\sigma^2}{2}b^{(N)}\log(-\log(1-\gamma_N)).
\end{align*}
Since we consider the efficiency driven regime, we have $b^{(N)}/(Nh^{(N)})\LimitN 0$. Also, it is easy to deduce that when $\gamma_N<1-\exp(-N)$, we have $\log(-\log(1-\gamma_N))<\log N$. Thus $\ApproxC{\ApproxIN}\sim\frac{\sigma^2}{2}Nh^{(N)}\log N$. Furthermore, $I^*_N-\ApproxIN=O(1)$, thus the bounds in \eqref{eq: ineq C(I) C(Itilde)} and \eqref{eq: ineq C Ctilde} are of $O(Nh^{(N)})$. By using the same argument as in the proof for the balanced regime,
\begin{align*}
    \frac{F_N(I^*_N,\beta^*_N)}{F_N(\ApproxIN,\ApproxBetaN)}=1-O(1/\log N).
\end{align*}
\hfill\rlap{\hspace*{-2em}\Halmos}
\endproof

\proof{Proof of Lemma \ref{lem: magnitude F}}
Following Equations \eqref{eq: ineq C(I) C(Itilde)} and \eqref{eq: ineq C Ctilde} and using the same arguments as in the proof of Theorem \ref{thm: order bounds indep}, we can find the same order bound for $F_N(I^*_N,\beta^*_N)/\ApproxCosts(\ApproxIN,\ApproxBetaN)=\sqrt{C_N(I^*_N)}/\sqrt{\ApproxC{\ApproxIN}}$.

In the  case that $\gamma_N=\gamma\in(0,1)$, we have 
\begin{multline*}
    \ApproxC{\ApproxIN}
    =Nh^{(N)}\frac{\sigma^2}{2}\left(\log N-\log(-\log(1-\gamma))-1\right)\\
    +(Nh^{(N)}+b^{(N)})\frac{\sigma^2}{2}\expect*{\left(G+\log(-\log(1-\gamma))\right)^+}.
\end{multline*}
Thus $\ApproxCosts(\ApproxIN,\ApproxBetaN)/(N\log N)=2\sqrt{N}\sqrt{\ApproxC{\ApproxIN}}/(N\log N)=O(\sqrt{h^{(N)}}/\sqrt{\log N})$.

When $\gamma_N\LimitN 0$, we have that $-\log(-\log(1-\gamma_N))\sim -\log(\gamma_N)$, thus $\ApproxIN\sim\frac{\sigma^2}{2} \log (N/\gamma_N)$. Also,
\begin{align*}
    \expect*{(G_N+\log(-\log(1-\gamma_N)))^+}\sim\expect*{(G_N+\log(\gamma_N))^+}\sim\gamma_N.
\end{align*}
From this it follows that 
\begin{align*}
    \ApproxC{\ApproxIN}\sim Nh^{(N)}\frac{\sigma^2}{2}\left(\log (N/\gamma_N)-1\right)+(Nh^{(N)}+b^{(N)})\frac{\sigma^2}{2}\gamma_N.
\end{align*}
Therefore, $2\sqrt{N}\sqrt{\ApproxC{\ApproxIN}}\gamma_N/(N\log (N/\gamma_N))=O(\gamma_N\sqrt{h^{(N)}}/\sqrt{\log (N/\gamma_N)})$.

When $\gamma_N\LimitN 1$, we have 
\begin{align*}
    \ApproxC{\ApproxIN}\sim&\frac{\sigma^2}{2}Nh^{(N)}(\log N-\log(-\log(1-\gamma_N)))+\frac{\sigma^2}{2}(Nh^{(N)}+b^{(N)})\log(-\log(1-\gamma_N))\nonumber\\
    =&\frac{\sigma^2}{2}Nh^{(N)}\log N+\frac{\sigma^2}{2}b^{(N)}\log(-\log(1-\gamma_N)).
\end{align*}
Thus, $2\sqrt{N}\sqrt{\ApproxC{\ApproxIN}}/\log N=O(N\sqrt{h^{(N)}}/\sqrt{\log N})$.
\hfill\rlap{\hspace*{-2em}\Halmos}
\endproof 
\subsection{Proofs of Section \ref{subsec: second order convergence}}\label{app: proofs second order convergence}
\proof{Proof of Lemma \ref{lem: lower bound second order convergence process}}
Let  $b_N=\sqrt{2\log N}-\log(4\pi\log N)/(2\sqrt{2\log N})$. Then
\begin{align*}
b_N\left(\frac{\max_{i\leq N}W_i(d\log N)}{\sigma\sqrt{d\log N}}-b_N\right)\LimitD G,
\end{align*}
with $G\sim\text{Gumbel}$, as $N\to\infty$, cf.\ \citet[p.~11, Ex.~1.1.7]{de2007extreme} for a proof. Observe that 
\begin{align*}
  &b_N\left(\frac{\max_{i\leq N}W_i(d\log N)}{\sigma\sqrt{d\log N}}-b_N\right)\nonumber\\
 = &\frac{1}{\sigma\sqrt{d}}\left(\sqrt{2\log N}-\frac{\log(4\pi\log N)}{2\sqrt{2\log N}}\right)\frac{\max_{i\leq N}W_i(d\log N)-\sigma\sqrt{2d}\log N+\frac{\sigma\sqrt{d}\log(4\pi\log N)}{2\sqrt{2}}}{\sqrt{\log N}}.
\end{align*}
Furthermore, $\beta d+\frac{\sigma^2}{2\beta}=\sigma\sqrt{2d}=\frac{\sigma^2}{\beta}$.
From this it follows that
\begin{align*}
\frac{\max_{i\leq N}W_i(d\log N)-\beta d\log N-\frac{\sigma^2}{2\beta}\log N}{\sqrt{\log N}}\LimitP 0,
\end{align*}
as $N\to\infty$.
Moreover, $\frac{W_A(d\log N)}{\sqrt{\log N}}\overset{d}{=}\frac{\sigma\sigma_A}{\sqrt{2}\beta}X$, with $X\sim\mathcal{N}(0,1)$. The statement follows.
\hfill\rlap{\hspace*{-2em}\Halmos}
\endproof 
\vspace{0.5cm}

\proof{Proof of Lemma \ref{lem: upper bound second order convergence process <d-e}}
To prove Lemma \ref{lem: upper bound second order convergence process <d-e}, we first observe that
\begin{multline}\label{eq: upper bound refl bm}
\frac{\max_{i\leq N}\left(\sup_{0<s<(d-\epsilon)\log N}(W_i(s)+W_A(s)-\beta s)\right)-\frac{\sigma^2}{2\beta}\log N}{\sqrt{\log N}}\\
\leq \frac{\max_{i\leq N}\left(\sup_{0<s<(d-\epsilon)\log N}\left(W_i(s)-\beta s\right)\right)-\frac{\sigma^2}{2\beta}\log N}{\sqrt{\log N}}
+\frac{\sup_{0<s<(d-\epsilon)\log N}W_A(s)}{\sqrt{\log N}}.
\end{multline}
We first focus on the first term on the right-hand side of \eqref{eq: upper bound refl bm}. We know that $\sup_{0<s<(d-\epsilon)\log N}\left(W_i(s)-\beta s\right)$ is a reflected Brownian motion, so we can write down its cumulative distribution function explicitly:
\begin{multline}\label{eq: cdf refl bm}
\probability*{\sup_{0<s<(d-\epsilon)\log N}\left(W_i(s)-\beta s\right)\leq x}\\
=1-\Phi\left(\frac{-x-\beta(d-\epsilon)\log N}{\sigma\sqrt{(d-\epsilon)\log N}}\right)-\text{exp}\left(-\frac{2\beta}{\sigma^2}x\right)\Phi\left(\frac{-x+\beta(d-\epsilon)\log N}{\sigma\sqrt{(d-\epsilon)\log N}}\right);
\end{multline}
see \cite[Eq.\ (1.1)]{abate1987transient}. From this together with the union bound, it follows that 
\begin{align}
    &\probability*{\max_{i\leq N}\sup_{0<s<(d-\epsilon)\log N}\left(W_i(s)-\beta s\right)\geq \frac{\sigma^2}{2\beta}\log N+x\sqrt{\log N}}\\
    &\quad\leq N     \probability*{\sup_{0<s<(d-\epsilon)\log N}\left(W_i(s)-\beta s\right)\geq \frac{\sigma^2}{2\beta}\log N+x\sqrt{\log N}}\nonumber\\
    &\quad=N\Phi\left(\frac{-\beta(2 d-\epsilon)\log N-x\sqrt{\log N}}{\sigma\sqrt{(d-\epsilon)\log N}}\right)+\text{exp}\left(-\frac{2\beta}{\sigma^2}x\sqrt{\log N}\right)\Phi\left(\frac{-\epsilon\beta\log N-x\sqrt{\log N}}{\sigma\sqrt{(d-\epsilon)\log N}}\right).\label{eq: upper bound second order convergence process <d-e}
\end{align}
The cumulative distribution of the normal distribution $\Phi$ satisfies $\Phi(-x)=1-\Phi(x)$. Furthermore, we have that $1-\Phi(x)\sim \exp(-x^2/2)/(\sqrt{2\pi}x)$ as $x\to\infty$; see \cite[Eq.\ (2.1.1), p.\ 49]{adler2007random}. This asymptotic equivalence gives us that the first term in \eqref{eq: upper bound second order convergence process <d-e} satisfies
\begin{align*}
N\Phi\left(\frac{-\beta(2 d-\epsilon)\log N-x\sqrt{\log N}}{\sigma\sqrt{(d-\epsilon)\log N}}\right)=&N\exp\left(-\frac{\beta^2(2d-\epsilon)^2}{2\sigma^2(d-\epsilon)}\log N(1+o(1))\right)\\
=&N\exp\left(-\frac{(2d-\epsilon)^2}{4d(d-\epsilon)}\log N(1+o(1))\right).
\end{align*}
For all $\epsilon\in(0,d)$, we have that $\frac{(2d-\epsilon)^2}{4d(d-\epsilon)}=\frac{4d^2-4d\epsilon+\epsilon^2}{4d(d-\epsilon)}>\frac{4d^2-4d\epsilon}{4d(d-\epsilon)}=1$. Thus, we can conclude that
$$
N\exp\left(-\frac{(2d-\epsilon)^2}{4d(d-\epsilon)}\log N(1+o(1))\right)\LimitN 0.
$$
With the asymptotic equivalence from \cite[Eq.\ (2.1.1), p.\ 49]{adler2007random} we get for the second term in \eqref{eq: upper bound second order convergence process <d-e} that
$$
\text{exp}\left(-\frac{2\beta}{\sigma^2}x\sqrt{\log N}\right)\Phi\left(\frac{-\epsilon\beta\log N-x\sqrt{\log N}}{\sigma\sqrt{(d-\epsilon)\log N}}\right)=\text{exp}\left(-\frac{\epsilon^2\beta^2}{2\sigma^2(d-\epsilon)}\log N(1+o(1))\right)\LimitN 0.
$$

For the second term on the right-hand side of \eqref{eq: upper bound refl bm}, we argue as follows: by filling in $\beta=0$, and replacing $\sigma$ with $\sigma_A$ in Equation \eqref{eq: cdf refl bm}, one can easily see that 
$$
\sup_{0<s<(d-\epsilon)\log N}W_A(s)\overset{d}{=}|W_A((d-\epsilon)\log N)|\overset{d}{=}\sqrt{(d-\epsilon)\log N}|X|,
$$
with $X\sim\mathcal{N}(0,1)$. Thus, we can use the upper bound in \eqref{eq: upper bound refl bm} and conclude that
\begin{align*}
   &\probability*{\frac{\max_{i\leq N}\left(\sup_{0<s<(d-\epsilon)\log N}(W_i(s)+W_A(s)-\beta s)\right)-\frac{\sigma^2}{2\beta}\log N}{\sqrt{\log N}}\geq x}\\
      &\quad\leq \probability*{\max_{i\leq N}\frac{\sup_{0<s<(d-\epsilon)\log N}(W_i(s)-\beta s)-\frac{\sigma^2}{2\beta}\log N}{\sqrt{\log N}}\geq x-y}+\probability*{\frac{\sup_{0<s<(d-\epsilon)\log N}W_A(s)}{\sqrt{\log N}}\geq y}\\
   &\quad\leq N\probability*{\frac{\sup_{0<s<(d-\epsilon)\log N}(W_i(s)-\beta s)-\frac{\sigma^2}{2\beta}\log N}{\sqrt{\log N}}\geq x-y}+\probability*{\frac{\sup_{0<s<(d-\epsilon)\log N}W_A(s)}{\sqrt{\log N}}\geq y}\\
   &\quad\LimitN \probability*{\lvert X\rvert>\frac{y}{\sqrt{d-\epsilon}}}.
\end{align*}
This last expression converges to 0 as $y\to\infty$, the lemma follows.

\hfill\rlap{\hspace*{-2em}\Halmos}
\endproof 

\proof{Proof of Lemma \ref{lem: upper bound second order convergence process >d+e}}
Let $\epsilon>0$ be given. Choose $\delta<\min\left(\frac{2 \left(\beta ^3 \epsilon +\beta  \sigma ^2\right)}{2 \beta ^2 \epsilon +\sigma ^2}-2    \sqrt{\frac{\beta ^2 \sigma ^2}{2 \beta ^2 \epsilon +\sigma ^2}},\frac{2 \beta ^3 \epsilon }{2 \beta ^2 \epsilon +\sigma ^2},\beta\right)$ and positive. Then
\begin{align*}
&\frac{\max_{i\leq N}\left(\sup_{s\geq(d+\epsilon)\log N}(W_i(s)+W_A(s)-\beta s)\right)-\frac{\sigma^2}{2\beta}\log N}{\sqrt{\log N}}\nonumber\\
&\quad\leq \frac{\max_{i\leq N}\left(\sup_{s\geq(d+\epsilon)\log N}(W_i(s)-(\beta-\delta) s)\right)-\frac{\sigma^2}{2\beta}\log N}{\sqrt{\log N}}+\frac{\sup_{s\geq(d+\epsilon)\log N}(W_A(s)-\delta s)}{\sqrt{\log N}}\nonumber\\
&\quad\leq  \frac{\max_{i\leq N}\left(\sup_{s\geq(d+\epsilon)\log N}(W_i(s)-(\beta-\delta) s)\right)-\frac{\sigma^2}{2\beta}\log N}{\sqrt{\log N}}+\frac{\sup_{s>0}(W_A(s)-\delta s)}{\sqrt{\log N}}.
\end{align*}
We have
\begin{align*}
\sup_{s\geq(d+\epsilon)\log N}(W_i(s)-(\beta-\delta) s)\overset{d}{=}&W_i((d+\epsilon)\log N)-(\beta-\delta)(d+\epsilon)\log N+ \sup_{s>0}(W_i'(s)-(\beta-\delta) s),
\end{align*}
with $(W_i',i\leq N)$ independent Brownian motions with mean 0 and variance $\sigma^2$.
We write $E_i= \sup_{s>0}(W_i'(s)-(\beta-\delta) s)$.
Hence, $E_i\sim\text{Exp}\left(\frac{2(\beta-\delta)}{\sigma^2}\right)$.
So
\begin{multline*}
\frac{\max_{i\leq N}\left(\sup_{s\geq(d+\epsilon)\log N}(W_i(s)-(\beta-\delta) s)\right)-\frac{\sigma^2}{2\beta}\log N}{\sqrt{\log N}}\\
\overset{d}{=}\frac{\max_{i\leq N}\left(W_i((d+\epsilon)\log N)+E_i\right)-\left(\frac{\sigma^2}{2\beta}+(\beta-\delta)(d+\epsilon)\right)\log N}{\sqrt{\log N}}.
\end{multline*}
By using the union bound and Chernoff's bound, we get that
\begin{align*}
\mathbb{P}\left(\max_{i\leq N}\left(W_i((d+\epsilon)\log N)+E_i\right)>x\right)\leq &N\mathbb{P}\left(W_i((d+\epsilon)\log N)+E_i>x\right)\nonumber\\
\leq &N\mathbb{E}\left[e^{sW_i((d+\epsilon)\log N)}\right]\mathbb{E}\left[e^{sE_i}\right]e^{-sx},
\end{align*}
for all $s>0$.
$\mathbb{E}\left[e^{sW_i((d+\epsilon)\log N)}\right]=e^{\frac{s^2(\sigma  \sqrt{(d+\epsilon ) \log N})^2}{2}}=N^{\frac{\sigma ^2 (d+\epsilon ) s^2}{2} }$ and $\mathbb{E}\left[e^{sE_i}\right]=\frac{2 (\beta -\delta )}{\sigma ^2}\Big/\left(\frac{2 (\beta -\delta )}{\sigma ^2}-s\right)$.
Hence,
\begin{multline}\label{eq: inequality chernoff bound}
\mathbb{P}\left(\max_{i\leq N}\left(W_i((d+\epsilon)\log N)+E_i\right)>x\sqrt{\log N}+\left(\frac{\sigma^2}{2\beta}+(\beta-\delta)(d+\epsilon)\right)\log N\right)\\
\leq N^{1+\frac{\sigma ^2 (d+\epsilon ) s^2}{2} -s\left(\frac{\sigma^2}{2\beta}+(\beta-\delta)(d+\epsilon)\right)}e^{-sx\sqrt{\log N}}\frac{\frac{2 (\beta -\delta )}{\sigma ^2}}{\frac{2 (\beta -\delta )}{\sigma ^2}-s}.
\end{multline}
Now, we choose $s^{\star}=\frac{\beta }{2 \beta ^2 \epsilon +\sigma ^2}+\frac{\beta -\delta }{\sigma ^2}$.
Because $\delta<\frac{2 \beta ^3 \epsilon }{2 \beta ^2 \epsilon +\sigma ^2}$, $s^{\star}<\frac{2(\beta-\delta)}{\sigma^2}$. Also,
\begin{align*}
1+\frac{\sigma ^2 (d+\epsilon ) s^{\star 2}}{2} -s^{\star}\left(\frac{\sigma^2}{2\beta}+(\beta-\delta)(d+\epsilon)\right)<0,
\end{align*}
because $\delta<\frac{2 \left(\beta ^3 \epsilon +\beta  \sigma ^2\right)}{2 \beta ^2 \epsilon +\sigma ^2}-2    \sqrt{\frac{\beta ^2 \sigma ^2}{2 \beta ^2 \epsilon +\sigma ^2}}$.
Therefore
\begin{align*}
\mathbb{P}\left(\max_{i\leq N}\left(W_i((d+\epsilon)\log N)+E_i\right)>x\sqrt{\log N}+\left(\frac{\sigma^2}{2\beta}+(\beta-\delta)(d+\epsilon)\right)\log N\right)\LimitN 0.
\end{align*}
Moreover, $\sup_{s>0}(W_A(s)-\delta s)\sim\text{Exp}\left(\frac{2\delta}{\sigma_A^2}\right)$.
Therefore, $\frac{\sup_{s>0}(W_A(s)-\delta s)}{\sqrt{\log N}}\LimitP 0$. 
The limit in \eqref{eq: upper bound second order convergence process >d+e} follows.
\hfill\rlap{\hspace*{-2em}\Halmos}
\endproof 

\proof{Proof of Lemma \ref{lem: upper bound second order convergence process d-e<s<d+e}}
First of all, we bound 
\begin{align*}
&\frac{\max_{i\leq N}\sup_{(d-\epsilon)\log N\leq s <(d+\epsilon)\log N}\left(W_i(s)+W_A(s)-\beta s\right)-\frac{\sigma^2}{2\beta}\log N}{\sqrt{\log N}}\nonumber\\
&\quad\leq \sup_{(d-\epsilon)\log N\leq s <(d+\epsilon)\log N}\frac{W_A(s)}{\sqrt{\log N}}+\frac{\max_{i\leq N}\sup_{(d-\epsilon)\log N\leq s <(d+\epsilon)\log N}(W_i(s)-\beta s)-\frac{\sigma^2}{2\beta}\log N}{\sqrt{\log N}}\nonumber\\
&\quad\leq \sup_{(d-\epsilon)\log N\leq s <(d+\epsilon)\log N}\frac{W_A(s)}{\sqrt{\log N}}+\frac{\max_{i\leq N}\sup_{s>0}(W_i(s)-\beta s)-\frac{\sigma^2}{2\beta}\log N}{\sqrt{\log N}}.
\end{align*}
We can write
\begin{align*}
    \sup_{(d-\epsilon)\log N\leq s <(d+\epsilon)\log N}\frac{W_A(s)}{\sqrt{\log N}}=&\frac{W_A((d-\epsilon)\log N)}{\sqrt{\log N}}+\sup_{0\leq s<2\epsilon\log N}\frac{W_A'(s)}{\sqrt{\log N}}\nonumber\\
    \overset{d}{=}&\sigma_A\sqrt{\frac{\sigma^2}{2\beta^2}-\epsilon}X_1+\sqrt{2\epsilon}\sigma_A|X_2|,
\end{align*}
with $X_1,X_2\sim\mathcal{N}(0,1)$ and independent, and $W_A'$ a Brownian motion with mean 0 and variance $\sigma_A^2$.
Furthermore, we have that
\begin{align*}
\frac{2\beta}{\sigma^2}\left(\max_{i\leq N}\sup_{s>0}(W_i(s)-\beta s)-\frac{\sigma^2}{2\beta}\log N\right) \LimitD G,
\end{align*}
as $N\to\infty$, with $G\sim \text{Gumbel}$. Therefore,
\begin{align*}
\frac{\max_{i\leq N}\sup_{s>0}(W_i(s)-\beta s)-\frac{\sigma^2}{2\beta}\log N}{\sqrt{\log N}}\LimitP 0,
\end{align*}
as $N\to\infty$. The statement follows.
\hfill\rlap{\hspace*{-2em}\Halmos}
\endproof 

\proof{Proof of Theorem \ref{thm: second order convergence process}}
We have the following lower bound:
\begin{multline*}
\mathbb{P}\left(\frac{\max_{i\leq N}\sup_{s>0}\left(W_i(s)+W_A(s)-\beta s\right)-\frac{\sigma^2}{2\beta}\log N}{\sqrt{\log N}}\geq x\right)\\
\geq \mathbb{P}\left(\frac{\max_{i\leq N}(W_i(d\log N)+W_A(d\log N))-\beta d\log N-\frac{\sigma^2}{2\beta}\log N}{\sqrt{\log N}}\geq x\right).
\end{multline*}
From this and Lemma \ref{lem: lower bound second order convergence process}, we know that 
\begin{align*}
\liminf_{N\to\infty}\mathbb{P}\left(\frac{\max_{i\leq N}\sup_{s>0}\left(W_i(s)+W_A(s)-\beta s\right)-\frac{\sigma^2}{2\beta}\log N}{\sqrt{\log N}}\geq x\right)\geq 1-\Phi\left(\frac{x\sqrt{2}\beta}{\sigma\sigma_A}\right).
\end{align*}
By using the union bound, we get
\begin{align*}
&\mathbb{P}\left(\frac{\max_{i\leq N}\sup_{s>0}\left(W_i(s)+W_A(s)-\beta s\right)-\frac{\sigma^2}{2\beta}\log N}{\sqrt{\log N}}\geq x\right)\nonumber\\
&\quad\leq\mathbb{P}\left(\frac{\max_{i\leq N}\sup_{0<s<(d-\epsilon)\log N}\left(W_i(s)+W_A(s)-\beta s\right)-\frac{\sigma^2}{2\beta}\log N}{\sqrt{\log N}}\geq x\right)\nonumber\\
&\qquad+\mathbb{P}\left(\frac{\max_{i\leq N}\sup_{(d-\epsilon)\log N\leq s <(d+\epsilon)\log N}\left(W_i(s)+W_A(s)-\beta s\right)-\frac{\sigma^2}{2\beta}\log N}{\sqrt{\log N}}\geq x\right)\nonumber\\
&\qquad+\mathbb{P}\left(\frac{\max_{i\leq N}\sup_{s\geq(d+\epsilon)\log N}\left(W_i(s)+W_A(s)-\beta s\right)-\frac{\sigma^2}{2\beta}\log N}{\sqrt{\log N}}\geq x\right).
\end{align*}
Combining this with the results from Lemmas \ref{lem: upper bound second order convergence process <d-e}, \ref{lem: upper bound second order convergence process >d+e} and \ref{lem: upper bound second order convergence process d-e<s<d+e} gives
\begin{multline*}
\limsup_{N\to\infty}\mathbb{P}\left(\frac{\max_{i\leq N}\sup_{s>0}\left(W_i(s)+W_A(s)-\beta s\right)-\frac{\sigma^2}{2\beta}\log N}{\sqrt{\log N}}\geq x\right)\\
\leq \mathbb{P}\left(\sigma_A\sqrt{\frac{\sigma^2}{2\beta^2}-\epsilon}X_1+\sqrt{2\epsilon}\sigma_A|X_2|>x\right),
\end{multline*}
with $X_1,X_2\sim\mathcal{N}(0,1)$ and independent. This upper bound holds for all $\epsilon>0$, therefore
\begin{align*}
&\limsup_{N\to\infty}\mathbb{P}\left(\frac{\max_{i\leq N}\sup_{s>0}\left(W_i(s)+W_A(s)-\beta s\right)-\frac{\sigma^2}{2\beta}\log N}{\sqrt{\log N}}\geq x\right)\nonumber\\
&\quad\leq \lim_{\epsilon\downarrow 0}\mathbb{P}\left(\sigma_A\sqrt{\frac{\sigma^2}{2\beta^2}-\epsilon}X_1+\sqrt{2\epsilon}\sigma_A|X_2|>x\right)\nonumber\\
&\quad=1-\Phi\left(\frac{x\sqrt{2}\beta}{\sigma\sigma_A}\right).
\end{align*}
Hence, the statement follows.
\hfill\rlap{\hspace*{-2em}\Halmos}
\endproof 
\proof{Proof of Lemma \ref{lem: L1 convergence}}
Because of the self-similarity property, we can assume without loss of generality that $\beta=1$. Let $d=\frac{\sigma^2}{2}$, and $X_N=\frac{\sqrt{2}}{\sigma\sigma_A}\frac{W_A(d\log N)}{\sqrt{\log N}}$. It is easy to see that $X_N\sim\mathcal{N}(0,1)$. Let $0<\epsilon<d$, we write
\begin{align*}
Q_i=\sup_{s>0}(W_i(s)+W_A(s)-s).
\end{align*}
First, observe that
\begin{align}
&\expect*{\left|\frac{\max_{i\leq N}Q_i-\frac{\sigma^2}{2}\log N}{\sqrt{\log N}}-\frac{\sigma\sigma_A}{\sqrt{2}}X_N\right|}\\
&\quad\leq \expect*{\left|\frac{\max_{i\leq N}Q_i-\frac{\sigma^2}{2}\log N}{\sqrt{\log N}}-\frac{\max_{i\leq N}W_i(d\log N)+W_A(d\log N)-\sigma^2\log N}{\sqrt{\log N}}\right|}\label{eq: 1}\\
&\qquad+\expect*{\left|\frac{\max_{i\leq N}W_i(d\log N)+W_A(d\log N)-\sigma^2\log N}{\sqrt{\log N}}-\frac{\sigma\sigma_A}{\sqrt{2}}X_N\right|}.\label{eq: upper bound conv in L1}
\end{align}
Due to \citet[Thm.~3.1]{pickands1968moment}, we obtain for the term in \eqref{eq: upper bound conv in L1} that
\begin{multline}\label{eq: 3}
\expect*{\left|\frac{\max_{i\leq N}W_i(d\log N)+W_A(d\log N)-\sigma^2\log N}{\sqrt{\log N}}-\frac{\sigma\sigma_A}{\sqrt{2}}X_N\right|}\\
=\expect*{\left|\frac{\max_{i\leq N}W_i(d\log N)-\sigma^2\log N}{\sqrt{\log N}}\right|}
\LimitN 0.
\end{multline}

Furthermore, because $Q_i>W_i(d\log N)+W_A(d\log N)-d\log N$, we can rewrite \eqref{eq: 1}:
\begin{align}
&\expect*{\left|\frac{\max_{i\leq N}Q_i-\frac{\sigma^2}{2}\log N}{\sqrt{\log N}}-\frac{\max_{i\leq N}W_i(d\log N)+W_A(d\log N)-\sigma^2\log N}{\sqrt{\log N}}\right|}\nonumber\\
&\quad=\expect*{\frac{\max_{i\leq N}Q_i-\frac{\sigma^2}{2}\log N}{\sqrt{\log N}}-\frac{\max_{i\leq N}W_i(d\log N)+W_A(d\log N)-\sigma^2\log N}{\sqrt{\log N}}}\nonumber\\
&\quad=\expect*{\frac{\max_{i\leq N}Q_i-\frac{\sigma^2}{2}\log N}{\sqrt{\log N}}}-\expect*{\frac{\max_{i\leq N}W_i(d\log N)-\sigma^2\log N}{\sqrt{\log N}}}.\label{eq: 11}
\end{align}
The second term in \eqref{eq: 11} converges to 0 as $N\to\infty$, which follows from the convergence in \eqref{eq: 3}. In order to find a converging upper bound for the first term in \eqref{eq: 11}, we write
\begin{align}
&\expect*{\frac{\max_{i\leq N}Q_i-\frac{\sigma^2}{2}\log N}{\sqrt{\log N}}}\nonumber\\
&\quad\leq \expect*{\frac{\max_{i\leq N}Q_i-\frac{\sigma^2}{2}\log N}{\sqrt{\log N}}\mathbbm{1}\left(-M\leq \frac{\max_{i\leq N}Q_i-\frac{\sigma^2}{2}\log N}{\sqrt{\log N}}\leq M \right)}\label{eq: 12}\\
&\qquad+\expect*{\frac{\max_{i\leq N}Q_i-\frac{\sigma^2}{2}\log N}{\sqrt{\log N}}\mathbbm{1}\left(\frac{\max_{i\leq N}Q_i-\frac{\sigma^2}{2}\log N}{\sqrt{\log N}}> M \right)}.\label{eq: 13}
\end{align}
For the term in \eqref{eq: 12}, we can conclude from Theorem \ref{thm: second order convergence process} together with the dominated convergence theorem that
\begin{multline*}
\expect*{\frac{\max_{i\leq N}Q_i-\frac{\sigma^2}{2}\log N}{\sqrt{\log N}}\mathbbm{1}\left(-M\leq \frac{\max_{i\leq N}Q_i-\frac{\sigma^2}{2}\log N}{\sqrt{\log N}}\leq M \right)}\\ 
\LimitN \expect*{\frac{\sigma\sigma_A}{\sqrt{2}}X\mathbbm{1}\left(-M\leq\frac{\sigma\sigma_A}{\sqrt{2}}X\leq M \right)}=0,    
\end{multline*}
with $X\sim\mathcal{N}(0,1)$.

In order to find a converging upper bound for the term in \eqref{eq: 13}, we bound 
\begin{align*}
\max_{i\leq N}Q_i
\leq \max_{i\leq N}\sup_{s>0}(W_i(s)-(1-1/\sqrt{\log N})s)+\sup_{s>0}(W_A(s)-s/\sqrt{\log N})=:Z_N.
\end{align*}
Then, we have the bound
\begin{align*}
&\expect*{\frac{\max_{i\leq N}Q_i-\frac{\sigma^2}{2}\log N}{\sqrt{\log N}}\mathbbm{1}\left(\frac{\max_{i\leq N}Q_i-\frac{\sigma^2}{2}\log N}{\sqrt{\log N}}\geq M\right)}\\
&\leq \expect*{\frac{Z_N-\frac{\sigma^2}{2}\log N}{\sqrt{\log N}}\mathbbm{1}\left(\frac{\max_{i\leq N}\sup_{s>0}(W_i(s)-(1-1/\sqrt{\log N})s)-\frac{\sigma^2}{2}\log N}{\sqrt{\log N}}\geq M/2\right)}\\
&\quad+\expect*{\frac{Z_N-\frac{\sigma^2}{2}\log N}{\sqrt{\log N}}\mathbbm{1}\left(\frac{\sup_{s>0}(W_A(s)-s/\sqrt{\log N})}{\sqrt{\log N}}\geq M/2\right)}.
\end{align*}
Because $\sup_{s>0}(W_A(s)-s/\sqrt{\log N})$ is exponentially distributed with mean $\sigma_A^2\sqrt{\log N}/2$, we have that
$$
\expect*{\frac{\sup_{s>0}(W_A(s)-s/\sqrt{\log N})}{\sqrt{\log N}}}=\frac{\sigma_A^2}{2}.
$$
Additionally, $\max_{i\leq N}\sup_{s>0}(W_i(s)-(1-1/\sqrt{\log N})s)$ is the maximum of $N$ i.i.d.\ exponentials with mean $\sigma^2/(2(1-1/\sqrt{\log N}))$, it is a standard result that
$$
\expect*{\max_{i\leq N}\sup_{s>0}(W_i(s)-(1-1/\sqrt{\log N})s)}=\frac{\sigma^2}{2(1-1/\sqrt{\log N})}\sum_{i=1}^N\frac{1}{i},
$$
see \cite{renyi1953theory}. From this, it follows that
$$
\expect*{\frac{\max_{i\leq N}\sup_{s>0}(W_i(s)-(1-1/\sqrt{\log N})s)-\frac{\sigma^2}{2}\log N}{\sqrt{\log N}}}\LimitN \frac{\sigma^2}{2}.
$$
Furthermore, due to the memoryless property of exponential random variables, we have that
\begin{multline*}
\expect*{\frac{\sup_{s>0}(W_A(s)-s/\sqrt{\log N})}{\sqrt{\log N}}\mathbbm{1}\left(\frac{\sup_{s>0}(W_A(s)-s/\sqrt{\log N})}{\sqrt{\log N}}\geq M/2\right)}\\
=\exp(-M/\sigma_A^2)\left(\frac{M}{2}+\frac{\sigma_A^2}{2}\right)\overset{M\to\infty}{\longrightarrow} 0,
\end{multline*}
and 
\begin{align*}
&\mathbb{E}\bigg[\frac{\max_{i\leq N}\sup_{s>0}(W_i(s)-(1-1/\sqrt{\log N})s)-\frac{\sigma^2}{2}\log N}{\sqrt{\log N}}\\
&\quad\quad\cdot\mathbbm{1}\bigg(\frac{\max_{i\leq N}\sup_{s>0}(W_i(s)-(1-1/\sqrt{\log N})s)-\frac{\sigma^2}{2}\log N}{\sqrt{\log N}}\geq M/2\bigg)\bigg]\\
&\quad\mathtoolsset{multlined-width=0.98\displaywidth}
\begin{multlined}
=\mathbb{E}\bigg[\max_{i\leq N}\bigg(\frac{\sup_{s>0}(W_i(s)-(1-1/\sqrt{\log N})s)-\frac{\sigma^2}{2}\log N}{\sqrt{\log N}}\\
\cdot\mathbbm{1}\bigg(\frac{\sup_{s>0}(W_i(s)-(1-1/\sqrt{\log N})s)-\frac{\sigma^2}{2}\log N}{\sqrt{\log N}}\geq M/2\bigg)\bigg)\bigg]
\end{multlined}
\\
&\quad\mathtoolsset{multlined-width=0.98\displaywidth}
\begin{multlined}
\leq N \mathbb{E}\bigg[\frac{\sup_{s>0}(W_i(s)-(1-1/\sqrt{\log N})s)-\frac{\sigma^2}{2}\log N}{\sqrt{\log N}}\\
\cdot\mathbbm{1}\bigg(\frac{\sup_{s>0}(W_i(s)-(1-1/\sqrt{\log N})s)-\frac{\sigma^2}{2}\log N}{\sqrt{\log N}}\geq M/2\bigg)\bigg]
\end{multlined}
\\
&\quad=N\exp\left(-\frac{2\left(1-1/\sqrt{\log N}\right) \left(\frac{\sigma ^2}{2} \log N+\frac{M}{2}\sqrt{\log N}\right)}{\sigma ^2}\right)\left(\frac{M}{2}+\frac{\sigma ^2  }{2 (1-1/\sqrt{\log N})}\right)\\
&\quad\LimitN 0,
\end{align*}
for $M>\sigma^2$. From these results, it follows that,
$$
\lim_{M\to\infty}\limsup_{N\to\infty}\expect*{\frac{\max_{i\leq N}Q_i-\frac{\sigma^2}{2}\log N}{\sqrt{\log N}}\mathbbm{1}\left(\frac{\max_{i\leq N}Q_i-\frac{\sigma^2}{2}\log N}{\sqrt{\log N}}\geq M\right)}=0.
$$
The lemma follows.

\hfill\rlap{\hspace*{-2em}\Halmos}
\endproof 
\subsection{Proofs of Section \ref{subsec: solution dependent}}\label{app: proofs dependent model}

\proof{Proof of Lemma \ref{lem: solution minimization problem}}
From Lemma \ref{lem: inventory minimization}, we know that the optimal inventory $\OptimalWAIN$ satisfies
\begin{align*}
    \frac{d}{dI}\InventoryExpCosts{\OptimalWAIN}=0.
\end{align*}
We have
\begin{align*}
    &\frac{d}{dI}\InventoryExpCosts{\OptimalWAIN}\nonumber\\
    &\quad=Nh^{(N)}-(Nh^{(N)}+b^{(N)})\probability*{\max_{i\leq N}Q_i>\OptimalWAIN}\nonumber\\
    &\quad=Nh^{(N)}-(Nh^{(N)}+b^{(N)})\probability*{\frac{\sqrt{2}}{\sigma\sigma_A}\frac{\max_{i\leq N}Q_i-\frac{\sigma^2}{2}\log N}{\sqrt{\log N}}>\frac{\sqrt{2}}{\sigma\sigma_A}\frac{\OptimalWAIN-\frac{\sigma^2}{2}\log N}{\sqrt{\log N}}}.
\end{align*}
Therefore, $\OptimalWAIN$ satisfies $\frac{\sqrt{2}}{\sigma\sigma_A}(\OptimalWAIN-\frac{\sigma^2}{2}\log N)/\sqrt{\log N}={P_N^A}^{-1}(1-\gamma_N)$.
\hfill\rlap{\hspace*{-2em}\Halmos}
\endproof 
\proof{Proof of Proposition \ref{prop: asymptotic approximation dependent}}
We have to find $I$ and $\beta$ such that $F_N(I,\beta)$ is minimized. As before, we know that the optimal $\ApproxWAIN$ should satisfy
\begin{align*}
   Nh^{(N)}- (Nh^{(N)}+b^{(N)})\probability*{\frac{\sigma^2}{2}\log N+\frac{\sigma\sigma_A}{\sqrt{2}}\sqrt{\log N}X>\ApproxWAIN}=0.
\end{align*}
Thus, $\ApproxWAIN$ as given in \eqref{eq: asymptotic I approximation dependent} minimizes $\ApproxWAC{I}$. We know that
\begin{align*}
 \mathbb{E}\left[\left(\frac{\sigma^2}{2}\log N+\frac{\sigma\sigma_A}{\sqrt{2}}\sqrt{\log N} X-\ApproxWAIN\right)^+\right]
    =&\int_{\frac{\ApproxWAIN-\frac{\sigma^2}{2}\log N}{\frac{\sigma\sigma_A}{\sqrt{2}}\sqrt{\log N}}}^{\infty}\left(\frac{\sigma^2}{2}\log N+\frac{\sigma\sigma_A}{\sqrt{2}}\sqrt{\log N} x-\ApproxWAIN\right)\phi(x)dx\nonumber\\
    =&\left(\frac{\sigma^2}{2}\log N-\ApproxWAIN\right)\mathbb{P}\left(\frac{\sigma\sigma_A}{\sqrt{2}}\sqrt{\log N}X\geq \ApproxWAIN-\frac{\sigma^2}{2}\log N\right)\nonumber\\
    &+\quad\frac{\sigma\sigma_A}{\sqrt{2}}\sqrt{\log N}\frac{1}{\sqrt{2\pi}}\text{exp}\left(-\frac{\left(\sigma ^2 \log N-2 \ApproxWAIN\right)^2}{4 \sigma ^2 \sigma _A^2 \log
   N}\right)\nonumber\\
   =&-\frac{\sigma\sigma_A}{\sqrt{2}}\sqrt{\log N}\Phi^{-1}(1-\gamma_N)\gamma_N\nonumber\\
   &+\frac{\sigma\sigma_A}{\sqrt{2}}\sqrt{\log N}\frac{1}{\sqrt{2\pi}}\text{exp}\Big(-\frac{1}{2}\Phi^{-1}(1-\gamma_N)^2\Big).
\end{align*}
The expression in Equation \eqref{eq: C tilde dep} follows.
\hfill\rlap{\hspace*{-2em}\Halmos}
\endproof 

\proof{Proof of Theorem \ref{thm: order bounds dep}}
Using Corollary \ref{cor: ratio optimal approx}, we have 
\begin{align*}
    \frac{F_N(\OptimalWAIN,\OptimalWABetaN)   }{F_N(\ApproxWAIN,\ApproxWABetaN)}=  \frac{2\sqrt{C_N(\OptimalWAIN)}\sqrt{\ApproxWAC{\ApproxWAIN}}}{C_N(\ApproxWAIN)+\ApproxWAC{\ApproxWAIN}}.
\end{align*}
First, assume $\ApproxWAC{\ApproxWAIN}>C_N(\ApproxWAIN)$. Then, $F_N(\OptimalWAIN,\OptimalWABetaN)/F_N(\ApproxWAIN,\ApproxWABetaN)>\sqrt{C_N(\OptimalWAIN)/\ApproxWAC{\ApproxWAIN}}$. We have 
\begin{align*}
    |\ApproxWAC{\ApproxWAIN}-C_N(\OptimalWAIN)|\leq& (2Nh^{(N)}+b^{(N)})|\OptimalWAIN-\ApproxWAIN|+(Nh^{(N)}+b^{(N)})\expect*{\left|\max_{i\leq N}Q_i-\frac{\sigma^2}{2}\log N-\frac{\sigma\sigma_A}{\sqrt{2}}X\right|}.
\end{align*}
We know by \citet[p.~305, Lem.~21.2]{vaart_1998}, that $(\OptimalWAIN-\ApproxWAIN)/\sqrt{\log N}\LimitN 0$. Furthermore, we prove in Lemma \ref{lem: L1 convergence} that $\expect*{\left|\max_{i\leq N}Q_i-\frac{\sigma^2}{2}\log N-\frac{\sigma\sigma_A}{\sqrt{2}}\sqrt{\log N}X\right|/\sqrt{\log N}}\LimitN 0$. From this it follows that $|\ApproxWAC{\ApproxWAIN}-C_N(\OptimalWAIN)|=o((Nh^{(N)}+b^{(N)})\sqrt{\log{N}})$. Since $\ApproxWAC{\ApproxWAIN}\sim \frac{\sigma^2}{2}Nh^{(N)}\log N$, we have $\frac{\sqrt{C_N(\OptimalWAIN)}}{\sqrt{\ApproxWAC{\ApproxWAIN}}}=1-o\left((Nh^{(N)}+b^{(N)})\sqrt{\log N}/(Nh^{(N)}\log N)\right)=1-o\left(1/\sqrt{\log N}\right)$.

Secondly, assume $\ApproxWAC{\ApproxWAIN}<C_N(\ApproxWAIN)$, then 
\begin{align*}
    \frac{F_N(\OptimalWAIN,\OptimalWABetaN)   }{F_N(\ApproxWAIN,\ApproxWABetaN)}>\frac{\sqrt{C_N(\OptimalWAIN)\ApproxWAC{\ApproxWAIN}}}{C_N(\ApproxWAIN)}=\frac{\sqrt{C_N(\OptimalWAIN)}}{\sqrt{C_N(\ApproxWAIN)}}\frac{\sqrt{\ApproxWAC{\ApproxWAIN}}}{\sqrt{C_N(\ApproxWAIN)}}.
\end{align*}
With an analogous derivation, we obtain the same order bound.
\hfill\rlap{\hspace*{-2em}\Halmos}
\endproof 

\proof{Proof of Lemma \ref{lem: bal reg dep}}
We have $\ApproxWAIN=\frac{\sigma^2}{2}\log N+\frac{\sigma\sigma_A}{\sqrt{2}}\sqrt{\log N}\Phi^{-1}(1-\gamma)$. Furthermore, $|\OptimalWAIN-\ApproxWAIN|= o(\sqrt{\log N})$, thus \eqref{eq: bal reg dep inventory} follows. Furthermore, by using the same argument as in Lemma \ref{lem: magnitude F}, \eqref{eq: bal reg dep costs} follows.
\hfill\rlap{\hspace*{-2em}\Halmos}
\endproof 
\subsection{Mixed-behavior approximations}\label{app: mixed behavior}
Though we have a symbolic expression for $\beta^{M}_N$ in \eqref{eq: symb expr beta mixed behavior}, it is not completely clear how to compute the part
\begin{multline*}
    \expect*{\left(\frac{\sigma^2}{2}\log N+\frac{\sigma\sigma_A}{\sqrt{2}}\sqrt{\log N}X+\frac{\sigma^2}{2}G-I^{M}_N\right)^+}\\
    =\int_{I^{M}_N}^{\infty} \probability*{\frac{\sigma^2}{2}\log N+\frac{\sigma\sigma_A}{\sqrt{2}}\sqrt{\log N}X+\frac{\sigma^2}{2}G>x}dx
\end{multline*}
in $\beta^{M}_N$. First, observe that we can write
\begin{align*}
    &\probability*{\frac{\sigma^2}{2}\log N+\frac{\sigma\sigma_A}{\sqrt{2}}\sqrt{\log N}X+\frac{\sigma^2}{2}G>x}\nonumber\\
    &\quad=\probability*{\frac{\sigma_A\sqrt{2}}{\sigma}\sqrt{\log N}X+G>\frac{2}{\sigma^2}x-\log N}\nonumber\\
    &\quad=\int_{-\infty}^{\infty}\probability*{\frac{\sigma_A\sqrt{2}}{\sigma}\sqrt{\log N}X>\frac{2}{\sigma^2}x-\log N-z}\exp(-\exp(-z)-z)dz.
\end{align*}
Now, we write $z=-\log s$. Then,
\begin{multline*}
\int_{-\infty}^{\infty}\probability*{\frac{\sigma_A\sqrt{2}}{\sigma}\sqrt{\log N}X>\frac{2}{\sigma^2}x-\log N-z}\exp(-\exp(-z)-z)dz\\
=  \int_{0}^{\infty}\probability*{\frac{\sigma_A\sqrt{2}}{\sigma}\sqrt{\log N}X>\frac{2}{\sigma^2}x-\log N+\log s}\exp(-s)ds.  
\end{multline*}
Thus,
\begin{align*}
     &\expect*{\left(\frac{\sigma^2}{2}\log N+\frac{\sigma\sigma_A}{\sqrt{2}}\sqrt{\log N}X+\frac{\sigma^2}{2}G-I^{M}_N\right)^+}\nonumber\\
     &\quad=\int_{I^{M}_N}^{\infty}\int_{0}^{\infty}\probability*{\frac{\sigma_A\sqrt{2}}{\sigma}\sqrt{\log N}X>\frac{2}{\sigma^2}x-\log N+\log s}\exp(-s)dsdx\nonumber\\
     &\quad=\int_{0}^{\infty}\int_{I^{M}_N}^{\infty}\probability*{\frac{\sigma_A\sqrt{2}}{\sigma}\sqrt{\log N}X>\frac{2}{\sigma^2}x-\log N+\log s}\exp(-s)dxds.
\end{align*}
It turns out that
\begin{align*}
    \int_{I^{M}_N}^{\infty}\probability*{\frac{\sigma_A\sqrt{2}}{\sigma}\sqrt{\log N}X>\frac{2}{\sigma^2}x-\log N+\log s}\exp(-s)dx
\end{align*}
can be expressed in terms of error functions. Thus, since $I^{M}_N$ can be numerically found by solving Equation \eqref{eq: computation I master}, $\expect*{\left(\frac{\sigma^2}{2}\log N+\frac{\sigma\sigma_A}{\sqrt{2}}\sqrt{\log N}X+\frac{\sigma^2}{2}G-I^{M}_N\right)^+}$ 
can be computed numerically as well. Observe that the procedure to obtain $I^{M}_N$ and $\beta^{M}_N$ is efficient and that its running time is independent of the system size $N$.
\end{APPENDIX}

\end{document}